\documentclass[a4paper,10pt]{amsart}
\usepackage[utf8]{inputenc}
\usepackage[english]{babel}

\usepackage{amsfonts}

\usepackage[english]{babel}
\usepackage{amsmath,amssymb,verbatim,mathrsfs,latexsym,paralist}
\usepackage{graphicx}
\usepackage{epstopdf}
\usepackage{indentfirst}
\usepackage{multirow}
\usepackage{amsthm}
\usepackage{longtable}
\usepackage{caption}

\usepackage{mathdots}
\allowdisplaybreaks[4]

\usepackage{tikz}
\usetikzlibrary{arrows,decorations.markings}
\usetikzlibrary{matrix}
\usetikzlibrary{graphs}
\usetikzlibrary{backgrounds}

\usepackage{capt-of}

	\def\DynkinNodeSize{1.5mm}
	\def\DynkinArrowLength{1.5mm}
	\tikzset{
		dnode/.style={
			circle,
			inner sep=0pt,
			minimum size=\DynkinNodeSize,
			fill=white,
			draw},
		middlearrow/.style={
			decoration={markings,
				mark=at position 0.6 with
				{\draw (0:0mm) -- +(+135:\DynkinArrowLength); \draw (0:0mm) -- +(-135:\DynkinArrowLength);},
			},
			postaction={decorate}
		},
		leftrightarrow/.style={
			decoration={markings,
				mark=at position 0.999 with
				{
					\draw (0:0mm) -- +(+135:\DynkinArrowLength); \draw (0:0mm) -- +(-135:\DynkinArrowLength);
				},
				mark=at position 0.001 with
				{
					\draw (0:0mm) -- +(+45:\DynkinArrowLength); \draw (0:0mm) -- +(-45:\DynkinArrowLength);
				},
			},
			postaction={decorate}
		},
		sedge/.style={
		},
		dedge/.style={
			middlearrow,
			double distance=0.5mm,
		},
		tedge/.style={
			middlearrow,
			double distance=1.0mm+\pgflinewidth,
			postaction={draw}, 
		},
		infedge/.style={
			leftrightarrow,
			double distance=0.5mm,
		},
	}

\newcommand{\Z}{\mathbb{Z}}
\newcommand{\C}{\mathbb{C}}

\newcommand{\wG}{\widehat{G}}
\newcommand{\GL}{\mathrm{GL}}
\newcommand{\SL}{\mathrm{SL}}

\newcommand{\End}{\mathrm{End}}
\newcommand{\Int}{\mathrm{Int}}
\newcommand{\Aut}{\mathrm{Aut}}

\newcommand{\Spin}{\mathrm{Spin}}
\newcommand{\Ad}{\mathrm{Ad}}
\newcommand{\ad}{\mathrm{ad}}
\newcommand{\mf}{\mathfrak}
\newcommand{\g}{\mf{g}}
\newcommand{\h}{\mf{h}}
\renewcommand{\a}{\mf{a}}
\newcommand{\cc}{\mf{c}}

\newcommand{\gl}{\mf{gl}}
\newcommand{\ssl}{\mf{sl}}
\renewcommand{\g}{\mf{g}}

\renewcommand{\u}{\mf{u}}

\newcommand{\z}{\mf{z}}
\newcommand{\ttt}{\mf{t}}

\newcommand{\CC}{\mathcal{C}}

\numberwithin{equation}{section}
\newtheorem{theorem}{Theorem}[section]

\theoremstyle{remark}
 
\theoremstyle{remark}

\usepackage{array,comment,makecell}
\newcolumntype{P}[1]{>{\raggedright\arraybackslash}m{#1}}
\graphicspath{{pics/}}
\newcommand\alb\allowbreak

\begin{document}

\title{A classification of four-tuples of spinors of a ten dimensional space}
\author{Willem A. de Graaf}
\address{
Dipartimento di Matematica\\
Universit\`{a} di Trento\\
Italy}
\author{Alexander Elashvili$^*$}
\thanks{$^*$ Partially supported by the ISF grant 1030/22}
\address{
Razmadze Mathematical Institute\\
Tbilisi State University\\
Georgia}
\author{Mamuka Jibladze}
\address{
Razmadze Mathematical Institute\\
Tbilisi State University\\
Georgia}

\date{}

\begin{abstract}
We use the theory of $\theta$-groups developed by Vinberg, along
with computations in the computer algebra system {\sf GAP}4, to
classify the orbits of  $\Spin(10,\C)\times \SL(4,\C)$ acting on its module
$\Delta_+\otimes \C^4$, where $\Delta_+$ is a half spin module of
$\Spin(10,\C)$. 
\end{abstract}  

\maketitle

\section{Introduction}

The theory of $\theta$-groups has been developed by Vinberg in the 70's
(\cite{vinberg,vinberg2}). They form a class of representations of linear
algebraic groups that arise from a $\Z/m\Z$-grading, or a $\Z$-grading, of a
semisimple complex Lie algebra. These representations share many properties
with the adjoint representation of a semisimple algebraic group on its
Lie algebra. In particular it is possible to use Vinberg's theory to classify
the orbits of a $\theta$-group. The first endeavour in this direction
was the classification by Vinberg and Elashvili of the threevectors of a
9-dimensional complex space (\cite{elashvin}). Subsequently a number of such
classifications have been undertaken (cf., \cite{antonyan2} (which has been
translated in \cite{oeding}), \cite{antelash}, \cite{gatim}, \cite{pervushin},
\cite{nurmiev}).

In this paper we use the same methods to classify the orbits of the
group $\Spin(10,\C)\otimes \SL(4,\C)$ acting on the space $\Delta_+\otimes
\C^4$, where $\Delta_+$ is a half spin module of $\Spin(10,\C)$. However,
in contrast to the above references we also heavily use explicit computation
in the computer algebra system {\sf GAP}4 (\cite{GAP4}) and especially its
package {\sf SLA} (\cite{SLA}). The latter package has a number of
implementations of algorithms for dealing with $\theta$-groups and
with other aspects of simple Lie algebras and their modules. 

This orbit classification has applications in geometry to the study of the
spinor tenfold, which is the orbit of the highest weight vector in
$\mathbb{P}(\Delta_+)$. The linear sections of small codimension of the
spinor tenfold are especially interesting, and those of codimension 4
are closely related to the orbits of $\Spin(10,\C)\times \SL(4,\C)$ on
$\Delta_+\otimes \C^4$. We refer to \cite{liumanivel} for an in depth
investigation that, among many other things, uses the results of the present
paper.

We start with a section with preliminaries
on Vinberg's $\theta$-groups, the particular $\theta$-group that we are
interested in and the construction of the module $\Delta_+$. A particular
property of the natural module of a $\theta$-group is that its elements
have a Jordan decomposition, dividing them into semisimple, nilpotent and
mixed elements (the latter are neither nilpotent nor semisimple). The semisimple
orbits are infinite in number, but be can divided into groups having the same
stabilizer. Section \ref{sec:semsim} is devoted to the classification of these
orbits.
This is followed by Section \ref{sec:stab} in which the methods are explained
that have been used to determine the stabilizers of the semisimple elements. 
Subsequently, in Sections \ref{sec:mixed}, \ref{sec:nilp} the classifications
of respectively the mixed and nilpotent orbits is given. In the last section
we also give the Hasse diagram of the closure relation between the nilpotent
orbits. It shows that the null cone has two irreducible components. 

\noindent{\bf Acknowledgment:} We thank Laurent Manivel for suggesting the
subject of the paper, and for his questions that led us to do many
interesting and fun computations.

\section{Preliminaries on Vinberg's $\theta$-groups}\label{sec:vinb}

In \cite{vinberg} Vinberg introduced and studied a class of representations of
linear algebraic groups that since have become known as $\theta$-groups.
They are constructed from a semisimple complex Lie algebra $\g$ together with
a $\Z/m\Z$-grading
$$\g = \bigoplus_{i\in \Z/m\Z} \g_i \text{ where } [\g_i,\g_j] \subset \g_{i+j}
\text{ for all } i,j.$$
If $m=\infty$ then $\Z/m\Z = \Z$. However, in this paper we restrict to the
case $m< \infty$. Such a grading corresponds to an automorphism $\theta$ of
$\g$ of order $m$ which is constructed as follows. Let $\omega\in\C$ be a fixed
primitive $m$-th root of unity and set $\theta(x) = \omega^i x$ for all
$x\in \g_i$ and extend $\theta$ to $\g$ by linearity. Then $\theta : \g\to\g$
is an automorphism of order $m$. Conversely, any automorphism $\theta$ of order
$m$ yields a $\Z/m\Z$-grading by letting $\g_i$ be the eigenspace of $\theta$
with eigenvalue $\omega^i$. 

We let $G$ be the identity component of the inner automorphism group of $\g$.
We have that the Lie algebra of $G$ is equal to $\ad \g = \{ \ad x \mid
x\in \g\}$ where $\ad x : \g \to \g$ is the adjoint map, $\ad x(y) = [x,y]$.
The subalgebra $\g_0$ is reductive and hence there is a unique connected
subgroup $G_0\subset G$ whose Lie algebra is $\ad \g_0 =\{ \ad x \mid x\in
\g_0\}$. Since $[\g_0,\g_1]\subset \g_1$ we have that $G_0$ acts on
$\g_1$. The representation $G_0 \to \GL(\g_1)$ is called a
$\theta$-representation and the group $G_0$ together with its module $\g_1$ is
called a $\theta$-group.

Results of Vinberg (\cite{vinberg,vinberg2}) and Vinberg and Elashvili
(\cite{elashvin}) make it possible to classify the orbits of a $\theta$-group.
The first observation in this direction is that the space $\g_1$ inherits the
Jordan decomposition of $\g$, that is, if $x\in \g_1$ and $s,n\in \g$ are
such that $x=s+n$ is its Jordan decomposition (with $s$ semisimple and
$n$ nilpotent, see \cite[\S 5.4]{hum}) then $s,n\in \g_1$. This divides the
orbits into three groups: semisimple (whose elements have $n=0$), nilpotent
(whose elements have $s=0$) and mixed (whose elements have $s,n$ both nonzero).

\subsection{Our main example}\label{sec:exa}

The finite order automorphisms of simple Lie algebras $\g$ were classified
by Kac, \cite[Chapter 8]{kac}, see also \cite[\S X.5]{helgason}. One consequence
of this classification is that a finite order automorphism is conveniently
described by its {\em Kac diagram}. This is an affine Dynkin diagram whose
nodes are labeled by non-negative integers. Here we do not go into the details,
but refer to \cite[\S 4.4.7]{onvi}. 
In this paper we consider a particular $\Z/4\Z$-grading 
of the simple Lie algebra $\g$ of type
$E_8$ induced by the automorphism of order 4 whose Kac diagram is

\begin{center}
\begin{tikzpicture}
\node[dnode,label=below:{\small $1$}] (1) at (0,0) {};
\node[dnode,label=left:{\small $2$}] (2) at (2,1) {};
\node[dnode,label=below:{\small $3$}] (3) at (1,0) {};
\node[dnode,label=below:{\small $4$}] (4) at (2,0) {};
\node[dnode,label=below:{\small $5$}] (5) at (3,0) {};
\node[dnode,fill=black,label=below:{\small $6$}] (6) at (4,0) {};
\node[dnode,label=below:{\small $7$}] (7) at (5,0) {};
\node[dnode,label=below:{\small $8$}] (8) at (6,0) {};
\node[dnode,label=below:{\small $0$}] (9) at (7,0) {};
\path (1) edge[sedge] (3)
(3) edge[sedge] (4)
(4) edge[sedge] (2)
(4) edge[sedge] (5)
(5) edge[sedge] (6)
(6) edge[sedge] (7)
(7) edge[sedge] (8)
(8) edge[sedge] (9);
\end{tikzpicture}
\end{center}

This means that all white nodes have label 0 whereas the black node has label 1.
Let $e_1,\ldots,e_8$, $f_1,\ldots,f_8$ be root vectors corresponding to
respectively the simple positive and negative roots. Let $\alpha_0$ be
the lowest root of the root system of $\g$ and let $e_0,f_0$ be root
vectors corresponding to $\alpha_0$ and $-\alpha_0$ respectively.
The automorphism $\theta$ given by this diagram satisfies $\theta(e_i)=e_i$,
$\theta(f_i)=f_i$ for $i\neq 6$ and $\theta(e_6) = i e_6$, $\theta(f_6)=
-i f_6$. Let $\g_k$ be the eigenspace of $\theta$ with eigenvalue
$i^k$. Then $\g = \g_0\oplus \g_1\oplus \g_2\oplus \g_3$ is a $\Z/4\Z$-grading
of $\g$. 

Using \cite[Proposition 8.6]{kac} or \cite[Proposition 17]{vinberg} it follows
from the Kac diagram that $\g_0$ is semisimple of type $D_5+A_3$ and $\g_1$ is
isomorphic, as $\g_0$-module to $\Delta_+\otimes \C^4$, where $\Delta_+$ is
the 16-dimensional semispinor module of the Lie algebra of type $D_5$ and
$\C^4$ is the natural 4-dimensional representation of the Lie algebra of
type $A_3$. 

Let $\wG = \mathrm{Spin}(10,\C)\times \SL(4,\C)$. Then the Lie algebra
$\hat\g$ of
$\wG$ is isomorphic to $\g_0$. Because $\wG$ is simply connected it follows
that there is a surjective homomorphism $\psi : \wG\to G_0$ whose differential
is a fixed isomorphism $\hat\g \to \g_0$. So since $G_0$ acts on $\g_1$, also
$\wG$ acts on that space. In this paper we determine the orbits of $G_0$ acting
on $\g_1$. This is the same as determining the orbits of $\wG$ on that space.
In the sequel we identify the modules $\g_1$ and $\Delta_+\otimes \C^4$.

In order to work with elements of $\Delta_+$ we say some words on the
construction of that module following \cite[\S VIII.13.4]{bou3}. Let $\ell\geq
4$ and set $n=2\ell$. Define the
$n\times n$-matrix 
$$A = \begin{pmatrix} 0 & \dots  & 1\\
  & \iddots &  \\
  1 & \dots & 0 \end{pmatrix}.$$
Define
$\mathfrak{o}(n,\C) = \{ x\in \gl(n,\C) \mid x^TA +Ax = 0\}$.
Then $\mathfrak{o}(n,\C)$ is a simple Lie algebra of type $D_\ell$.
Let $v_1,\ldots, v_n$ be the standard basis of $\C^n$ and consider the
bilinear form $\Psi$ on $\C^n$ defined by $A$, that is, $\Psi(u,v) = u^TAv$. 
Then the Clifford algebra $C$ is the associative algebra generated by
$v_1,\ldots,v_n$ subject to the relations $v_iv_j+v_jv_i= \Psi(v_i,v_j)\cdot 1$.
It is known that the products $v_{i_1}\cdots v_{i_k}$ for $k\geq 0$ and
$i_1<i_2<\cdots <i_k$ form a basis of $C$.

Define the map $f : \mathfrak{o}(n,\C) \to \C$ by $f(a) = \tfrac{1}{2}
\sum_{i=1}^n(av_i) v_{n+1-i}$. Then $f$ is linear and $f([a,b]) = [f(a),f(b)]$
where the latter is the commutator in $C$ (see
\cite[Lemme 1, \S VIII.13 no 4]{bou3}). We partition the basis of $\C^n$ into two
sets $F=\{v_1,\ldots,v_\ell\}$, $F'= \{u_{1},\ldots,v_\ell\}$ where
$u_i = v_{\ell+i}$ for $1\leq i\leq \ell$. Let $U$ denote
the span of $F'$ and let 
$E$ denote the exterior algebra of $U$, that is,
$$E=\bigwedge^0 U +\bigwedge^1 U +\cdots +\bigwedge^\ell U.$$
For $u\in F'$ and $v\in F$ we define the endomorphisms $\lambda(u)$,
$\lambda(v)$ of $E$ by
\begin{align*}
&\lambda(u) \cdot u_{i_1}\wedge \cdots \wedge u_{i_k} = u\wedge u_{i_1}\wedge
\cdots \wedge u_{i_k}\\
&\lambda(v) \cdot u_{i_1}\wedge \cdots \wedge u_{i_k} = \sum_{j=1}^k (-1)^{j-1}
\Psi(u_{i_j},v)
u_{i_1}\wedge \cdots \wedge u_{i_{j-1}} \wedge  u_{i_{j+1}}\wedge \cdots \wedge
u_{i_k}. \end{align*}  
Then the map $\lambda : F\cup F' \to \End(E)$
extends to a homomorphism $\lambda : C\to \End(E)$ and we get a representation
$\rho : \mathfrak{o}(n,\C) \to \gl(E)$ by $\rho(a) = \lambda(f(a))$.

It is known that $E$, as $\mathfrak{o}(n,\C)$-module splits as the direct
sum of two irreducible modules $\Delta_+$ and $\Delta_-$ which are the sum
of the
$\wedge^k U$ with $k$ even, respectively odd. They are called the semispinor
modules of $\mathfrak{o}(n,\C)$. 

In our case we have $\ell=5$ so that $U$ is of dimension 5 and 
$$\Delta_+ = \bigwedge^0 U + \bigwedge^2 U +\bigwedge^4 U.$$

Let $w_1,\ldots,w_4$ be the elements of the
standard basis of $\C^4$. An element of $\Delta_+\otimes \C^4$ is then a linear
combination of elements of the form
$$u_{i_1}\wedge u_{i_2} \wedge \cdots \wedge u_{i_k} \otimes w_j$$
(where $k\in \{0,2,4\}$). 
Throughout we denote this element by $(i_1,i_2,\ldots,i_k)\otimes j$. We
identify the modules $\g_1$ and $\Delta_+\otimes \C^4$. So we will describe
elements of these two modules 
by giving them as linear combinations of basis elements of the
form $(i_1,i_2,\ldots,i_k)\otimes j$.

An element of  $\Delta_+\otimes \C^4$ can be written as $a_1\otimes w_1+\cdots
+a_4\otimes w_4$ where $a_i\in \Delta_+$. For this reason we say that the
elements of  $\Delta_+\otimes \C^4$ are four-tuples of spinors of a
10-dimensional space. 

The above basis vectors are weight vectors for the representation $\Delta_+\otimes \C^4$ as follows.

Weights of $\Delta_+$ can be viewed as quintuples $\frac12(\pm\varepsilon_1\pm\cdots\pm\varepsilon_5)$, with the even number of $+$ signs,
where $\varepsilon_1$, ..., $\varepsilon_5$ form an orthonormal basis of the $5$-dimensional weight space for $D_5$.
The basis vector $u_{i_1}\wedge \cdots \wedge u_{i_k}$ has weight
\[
\varepsilon_{i_1}+\cdots+\varepsilon_{i_k}-\frac12(\varepsilon_1+\cdots+\varepsilon_5).
\]

As for $\C^4$, weights of the standard representation $\C^4$ can be realized as vectors $e_1$, $e_2$, $e_3$, $e_4$ with $e_1+e_2+e_3+e_4=0$ in a $4$-dimensional vector space,
with scalar products $(e_i,e_i)=3/4$ and $(e_i,e_j)=-1/4$ for $i\ne j$ (see e.~g. \cite[p. 292]{onvi}). The basis vector $w_i$ has weight $e_i$. 

Weights of $\Delta_+\otimes \C^4$ are then pairwise sums of weights of $\Delta_+$ and of $\C^4$ in the orthogonal sum of the corresponding weight spaces.

Note that the possible values of scalar products of weights of $\Delta_+$ corresponding to weight vectors $u_{i_1}\wedge \cdots \wedge u_{i_k}$ and $u_{j_1}\wedge \cdots \wedge u_{j_l}$ are $5/4$ (when the weights are equal), $1/4$ (when the symmetric difference of the sets $\{i_1,...,i_k\}$ and $\{j_1,...,j_l\}$ has two elements) and $-3/4$ (when this symmetric difference has four elements).

Thus the possible scalar products for weights of $\Delta_+\otimes \C^4$ turn out to be $2$ (when the weights are equal), $0$, $1$ and $-1$. 
In fact under the identification of $\Delta_+\otimes \C^4$ with $\g_1$, our basis vectors $(i_1,i_2,\ldots,i_k)\otimes j$
correspond to root vectors of $\g$ for certain roots of the $E_8$ root system, and one can check that the scalar product described above coincides with the one inherited from the scalar products of $E_8$ roots.

We will use this scalar product below to display, for a linear combination of weight vectors with nonzero coefficients, its \emph{Dynkin scheme}. The latter is a graph with nodes the corresponding weights,
connected with a solid edge if their scalar product is $-1$ (angle $120^\circ$), with a dashed edge when the scalar product is $1$ (angle $60^\circ$), and with no edge if the scalar product is $0$ (i. e. the weights are orthogonal).

For example,

\includegraphics[scale=.8]{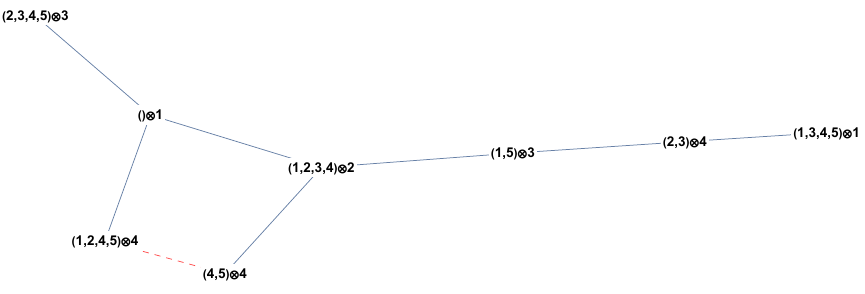}

\noindent is the Dynkin scheme of the vector $()\otimes 1 +(1,3,4,5)\otimes 1 +(1,2,3,4)\otimes2 +{(1,5)\otimes3}+(2,3,4,5)\otimes3+(2,3)\otimes 4+(4,5)\otimes4+(1,2,4,5)\otimes4$.

\section{The semisimple orbits}\label{sec:semsim}

Consider a $\Z/m\Z$-grading $\g = \oplus_{i\in \Z/m\Z} \g_i$ of the semisimple
Lie algebra $\g$. As seen in Section \ref{sec:vinb} this yields the reductive
group $G_0$ which acts on $\g_1$. Here we first describe some general facts
concerning the semisimple orbits in $\g_1$. Then we specialize to the particular
example that we are interested in.

A {\em Cartan subspace} in $\g_1$ is a maximal subspace consisting of commuting
semisimple elements. By \cite[Theorem 1]{vinberg} two Cartan subspaces of
$\g_1$ are $G_0$-conjugate. It follows that every semisimple orbit in $\g_1$
has a point in a fixed Cartan subspace.

Let $\cc\subset \g_1$ be a Cartan subspace and define
\begin{align*}
Z_{G_0}(\cc) &= \{g\in G_0 \mid g(x)=x \text{ for all } x\in \cc\}\\
N_{G_0}(\cc) &= \{ g\in G_0\mid g(x)\in \cc \text{ for all } x\in\cc\}
\end{align*}
and set $W_0 = N_{G_0}(\cc)/Z_{G_0}(\cc)$. The group $W_0$ is called the
{\em little Weyl group} of the grading. A linear transformation $T$ of a
complex vector space is said to be a complex reflection if $T-1$ has rank 1.
Then the kernel of $T-1$ is a hyperplane, which we call the {\em reflection
hyperplane} of $T$.
Vinberg showed that $W_0$ is generated by complex reflections,
\cite[Theorem 8]{vinberg}.

By \cite[Theorem 2]{vinberg}
any two elements of $\cc$ are $G_0$-conjugate if and only if they are
$W_0$-conjugate. It follows that classifying the semisimple orbits in
$\g_1$ reduces to classifying the $W_0$-orbits in $\cc$. 
However, this statement can still be refined in the following way.
For $p\in \cc$ let
$W_p = \{w\in W_0 \mid w(p)=p\}$ be its stabilizer in $W_0$. By
\cite[Theorem 9.44]{lehta} $W_p$ is a reflection subgroup of $W_0$, that is,
it is generated by complex reflections. Now let $M\subset W_0$ be a
reflection subgroup, and define
\begin{align*}
\cc_M &= \{ q\in \cc \mid M\subset W_q\}\\
\cc_M^\circ &= \{ q\in \cc \mid M = W_q\}.
\end{align*}  

Then $\cc_M$ is the intersection of all reflection hyperplanes of all complex
reflections in $M$. So it is a linear subspace of $\cc$. Furthermore,
$\cc_M^\circ$ is the set of $q\in \cc_M$ that do not lie in the reflection
hyperplane of any complex reflection in $W_0$ that does not lie in $M$.
So $\cc_M^\circ$ is Zariski-open in $\cc_M$. Now for $w\in W_0$ and reflection
subgroups $M,M'\subset W_0$ we have 
$$w(\cc_M^\circ) = \cc_{M'}^\circ \text{ if and only if } wM w^{-1} = M'.$$
Let $M_1,\ldots,M_r\subset W_0$ be the reflection subgroups of $W_0$ such
that for each $p\in \cc$ the stabilizer $W_p$ is conjugate to exactly one
$M_i$. Then it follows that every semisimple orbit has a point in a unique
$\cc_{M_i}^\circ$. Moreover, two elements
$p,p'\in \cc_{M_i}^\circ$ are $W_0$-conjugate if and only if they are conjugate
under the group $\Gamma_{i} = N_{W_0}(M_i)/M_i$. We conclude that
the classification of the semisimple $G_0$-orbits also reduces to the
classification of the $\Gamma_{i}$-orbits in $\cc_{M_i}^\circ$ for
$1\leq i\leq r$. 

\subsection{The semisimple orbits of our example}

We let $\g$ and its grading be as in Section \ref{sec:exa}. 

A computation in {\sf GAP}4 shows that the following elements span a
Cartan subspace $\cc$ in $\g_0$ 

\begin{align*}
p_1&= -(3,5)\otimes 1 +(1,2,4,5)\otimes 2-(2,4)\otimes 3-(1,3)\otimes 4,\\
p_2&= -(2,5)\otimes 1+(1,3,4,5)\otimes 2 +(3,4)\otimes 3 +(1,2)\otimes 4,\\
p_3&= (1,2,3,4)\otimes 1+()\otimes 2+(1,2,3,5)\otimes 3-(4,5)\otimes 4,\\
p_4&= (1,4)\otimes 1+(2,3)\otimes 2-(1,5)\otimes 3 +(2,3,4,5)\otimes 4.  
\end{align*} 

Each of the $p_1$, ..., $p_4$ can be viewed as the sum $e_\alpha+e_\beta+e_\gamma+e_\delta$ of four root vectors of $\g$ in $\g_1$.

The Dynkin scheme for each of these quadruples $\alpha$, $\beta$, $\gamma$, $\delta$ is of type extended $A_3$ (a square),
which by \cite{kostant} implies that $e_\alpha+e_\beta+e_\gamma+e_\delta$ is a regular semisimple element in the
subalgebra of type $A_3$ in $\g$ generated by $e_\alpha,e_\beta,e_\gamma,e_\delta$. It follows that each of the $p_1$, ..., $p_4$ is a semisimple element of $\g$,
hence by \cite{vinberg} they are semisimple elements of $\Delta_+\otimes \C^4$.

Note that the obtained four subalgebras of type $A_3$ corresponding to $p_1$, ..., $p_4$ do not centralize each other.
However, the above root vectors can be grouped in such a way that the sums $e_\alpha+e_\beta$, $e_\gamma+e_\delta$ generate a subalgebra of type $2A_1$ in $\g$,
and the obtained four subalgebras of type $2A_1$ corresponding to $p_1$, ..., $p_4$ pairwise centralize each other.

The Dynkin scheme for all 16 weights involved in the $p_1$, ..., $p_4$ looks as follows:

\includegraphics[scale=.4]{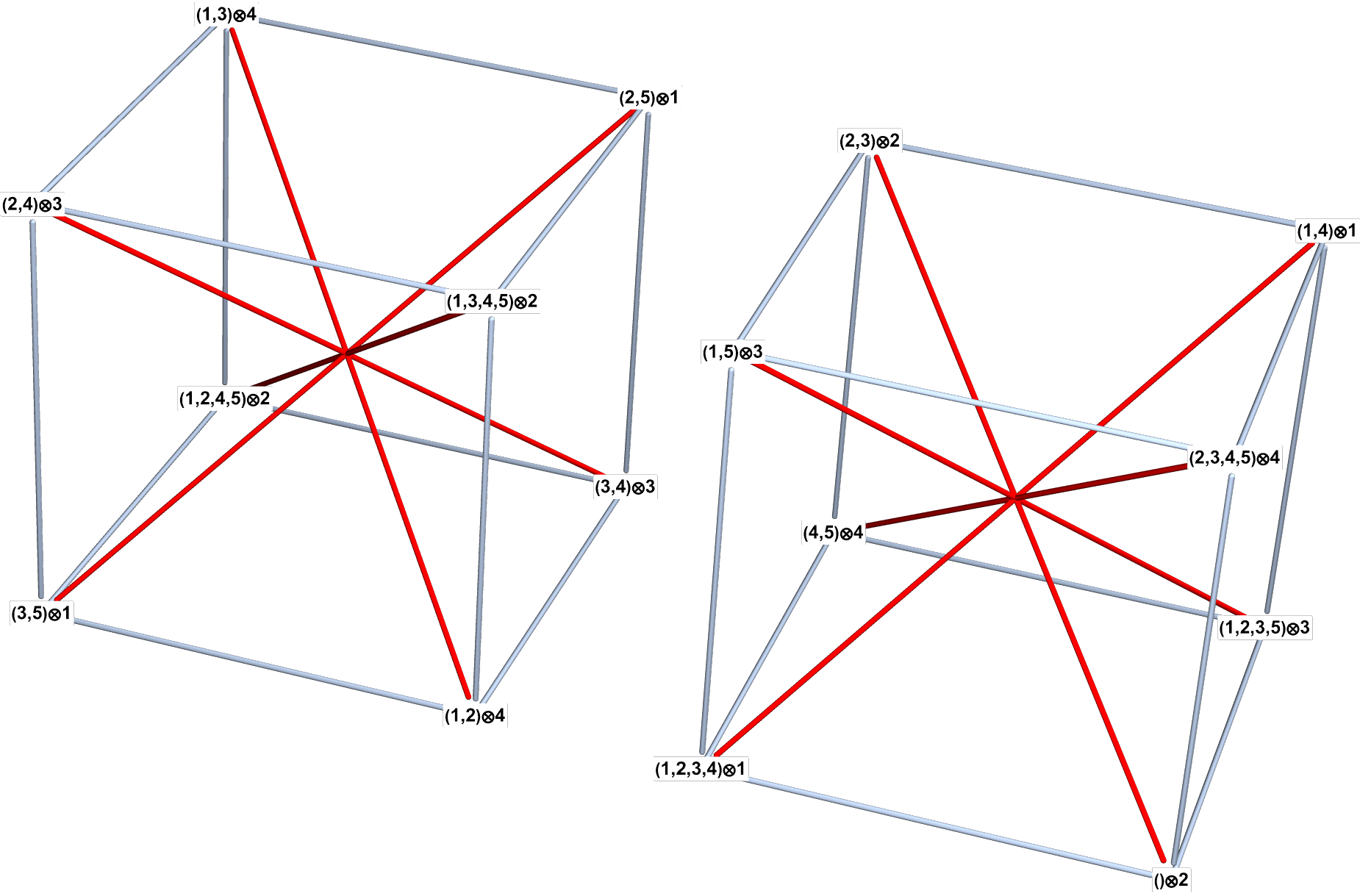}

The linear combinations defining $p_1$, ..., $p_4$ correspond to pairs of opposite sides of the cubes in the picture.

The little Weyl group $W_0$ is generated by the following five complex
reflections (which are given by their matrices with respect to the above basis
of $\cc$):

\begin{align*}
& s_1=\begin{pmatrix} -1&0&0&0\\ 0&1&0&0\\ 0&0&1&0\\ 0&0&0&1 \end{pmatrix},\,\,
s_2=\begin{pmatrix} 0&-1&0&0\\ -1&0&0&0\\ 0&0&1&0\\0&0&0&1\end{pmatrix},\,\,
  s_3=\begin{pmatrix} 0&-i&0&0\\ i&0&0&0\\ 0&0&1&0\\ 0&0&0&1\end{pmatrix},\\
    & s_4=\tfrac{1}{2}\begin{pmatrix} 1&-1&-1&-1\\-1&1&-1&-1\\-1&-1&1&-1\\-1&-1&-1&1
    \end{pmatrix},\,\,
    s_5=\tfrac{1}{2}\begin{pmatrix} 0&0& -1-i & -1+i\\
      0&2&0&0\\ -1+i & 0 & 1 & i\\
     -1-i & 0 & -i & 1\end{pmatrix}.
\end{align*}

There are two ways to compute these generators. Firstly, we have used 
an improved version of the algorithm in \cite{graoriente} to compute them.
Secondly, a computation shows that $\cc$ lies in a unique
Cartan subalgebra $\h$ of $\g$, which is the centralizer of $\cc$ in $\g$.
This Cartan subalgebra is $\theta$-stable. So we can restrict $\theta$ to
$\h$ and consider the centralizer $W^\theta = \{ w\in W \mid w\theta =
\theta w\}$ where $W$ is the Weyl group of the root system of $\g$ with
respect to $\h$.
It is known (see the table in \cite[\S 9]{vinberg}) that $W_0$ is
isomorphic to the group number 31 in the Shephard-Todd
classification of the finite irreducible reflection groups \cite{shephard}. In particular
$|W_0| = 46080$. Furthermore, we have that $W_0\subset \{ w|_\cc \mid w\in
W^\theta\}$. Finally, by computing $W^\theta$ explicitly it is readily checked
that $|W^\theta| =46080$. So the elements of $W^\theta$ all have different
restrictions to $\cc$ and $W_0 = \{ w|_\cc \mid w\in W^\theta\}$.

Including the group itself and the trivial subgroup, $W_0$ has fourty-three
reflection subgroups up to conjugacy, see \cite[Table 9]{taylor}. 
However, only nine of them are equal to the stabilizer of a point. They
are given in Table \ref{tab:semsim}, together with other information.
Here we remark that the elements of a set $\cc_{M_i}^\circ$ all have the
{\em same} stabilizer in $G_0$ (see Theorem \ref{thm:stab} below).
The identity component of such a stabilizer is reductive. We describe it
by giving the type of the root system and the dimension of its centre, where
$T_k$ denotes a $k$-dimensional centre. In all cases the component group
of the stabilizer is a direct product of cyclic groups of order 2. In the
sixth column of the table we indicate a product of $k$ such groups by
 $\CC_2^k$.

\begin{table}[htb]\caption{Stabilizers of points in $\cc$ up to conjugacy. The first column has the index $i$. The second and third columns list generators of the reflection subgroup $M_i$ of $W_0$ and its size. The fourth column has the basis elements of the space $\cc_{M_i}$. The fifth and sixth columns display the type of the identity component and component group of the centralizer in $G_0$ of any $p\in\cc_{M_i}^\circ$. The last column has the size of $\Gamma_i$.}\label{tab:semsim}
\begin{tabular}{|c|c|c|c|c|c|c|}
\hline
$i$ & generators of $M_i$ & size & $\cc_{M_i}$ & $Z_{G_0}(p)^\circ$ & $K$ & $|\Gamma_i|$\\
\hline
1 & & 1 & $p_1,p_2,p_3,p_4$ & 1 & $\CC_2^4$ & 46080\\
\hline
2 & $s_1$ & 2 & $p_2,p_3,p_4$ & $T_1$ & $\CC_2^3$ & 384\\
\hline
3 & $s_1$, $s_4s_2s_3s_5s_4s_5s_3s_2s_4$ & 4 & $p_2,p_3$ & $T_2$ & $\CC_2^2$ & 32\\
\hline
4 &  $s_1$, $s_4$ & 6 & $p_2-p_3,p_3-p_4$ & $A_1$ & $\CC_2^2$ & 24 \\
\hline
5 &  $s_4s_2s_4$, $s_4s_3s_5s_4s_5s_3s_4$, & 16 & $p_1,p_2$ &  $A_1+T_3$ &
$\CC_2^2$ & 96 \\
& $s_4s_2s_1s_5s_4s_5s_1s_2s_4$  & & & & & \\
\hline
6 & $s_1s_2s_1$, $s_1s_2s_4s_3s_5s_3s_4s_2s_1$, & 12 & $p_1+p_2+p_3$ & $A_1+T_1$ &
$\CC_2$ & 4 \\
& $s_4s_2s_3s_5s_4s_5s_3s_2s_4$ & & & & & \\
\hline
7 &  $s_1$, $s_3s_5s_4s_5s_3$, $s_2s_4s_2$ & 24 &  $p_2+p_3$ & $2A_1$ & $\CC_2$
 & 4 \\
\hline
8 & $s_2s_1s_2$, $s_2s_5s_2$, $s_3s_5s_3$, $s_4s_5s_4$ & 192 & $p_1$ &
$2A_1+A_2+T_1$ & $\CC_2$ & 4 \\
\hline
9 & $s_1$, $s_2$, $s_3$, $s_4$, $s_5$ & 46080 & 0 & $D_5+A_3$ & 1 & 1 \\
\hline
\end{tabular}
\end{table}

In the last column of the table we list the sizes of the groups $\Gamma_i =
N_{W_0}(M_i)/M_i$. More explicitly, we have the following descriptions of these
groups, where the matrices are given with respect to the bases of
$\cc_{M_i}$ listed in Table \ref{tab:semsim}. 

\begin{align*}
\Gamma_1 &= W_0\\
\Gamma_2 &= \left\langle
  \begin{pmatrix} 1&0&0\\0&1&0\\0&0&-1\end{pmatrix},\,
  \tfrac{1}{2}\begin{pmatrix} 0 & 1+i & 1-i\\1+i & 1 & i \\1-i & i & 1
  \end{pmatrix} \right\rangle\\
\Gamma_3 &= \left\langle  \begin{pmatrix} 1&0\\0&i\end{pmatrix},\,
  \begin{pmatrix} 0&1\\1&0 \end{pmatrix}. \right\rangle\\
\Gamma_4 &= \left\langle   \begin{pmatrix}1&0\\1&-1\end{pmatrix},\,
  \begin{pmatrix}-\tfrac{1}{2}-\tfrac{1}{2}i & 1 \\ \tfrac{1}{2} &
    \tfrac{1}{2}-\tfrac{1}{2}i \end{pmatrix} \right\rangle\\
\Gamma_5 &= \left\langle \begin{pmatrix} 1&0\\0&i\end{pmatrix},
  (\tfrac{1}{2}+\tfrac{1}{2}i) \begin{pmatrix} 1&1\\1&-1\end{pmatrix}
    \right\rangle,  \\
\Gamma_{6} &=  \Gamma_7=\Gamma_8= \langle i \rangle. \\
\end{align*}

Now we give explicit polynomials defining the open sets $\cc_{M_1}^\circ$
inside $\cc_{M_i}$. If $\cc_{M_i}$ is 1-dimensional then this is obvious: let
$p$ be a basis element of $\cc_{M_i}$, then $\cc_{M_i}^\circ = \{ \mu p \mid
\mu \neq 0\}$. For the other cases we have the following statements which are
obtained by explicit computation in {\sf GAP}4. 

The element $x_1p_1+x_2p_2+x_3p_3+x_4p_4$ lies in $\cc_{M_1}^\circ$ if and only if
the following polynomials are nonzero:

\begin{align*}
& x_1x_2x_3x_4, \,\, x_2^4 - 2x_2^2x_3x_4 + \tfrac{1}{4}x_3^4 +
\tfrac{1}{2}x_3^2x_4^2 + \tfrac{1}{4}x_4^4,\\
& x_1^4-x_2^4,\,\, x_3^4-x_4^4, \,\, x_2^4 + 2x_2^2x_3x_4 + \tfrac{1}{4}x_3^4 +
\tfrac{1}{2}x_3^2x_4^2 +\tfrac{1}{4}x_4^4,\\
& x_1^2 - 2x_1x_2 + x_2^2 - x_3^2 - 2x_3x_4 - x_4^2, \,\, x_1^2 - 2x_1x_2 + x_2^2 -
x_3^2 + 2x_3x_4 - x_4^2,\\
& x_1^2 + 2x_1x_2 + x_2^2-x_3^2-2x_3x_4-x_4^2, \,\, x_1^2 + 2x_1x_2 + x_2^2-x_3^2
+ 2x_3x_4 - x_4^2,\\
& x_1^4 + 2x_1^2x_2^2 - 8x_1x_2x_4^2 + x_2^4 + 4x_4^4, \,\,
x_1^4 + 2x_1^2x_2^2 + 8x_1x_2x_4^2 + x_2^4 + 4x_4^4,\\
& x_1^4 + 2x_1^2x_2^2 - 8x_1x_2x_3^2 + x_2^4 + 4x_3^4, \,\,
x_1^4 + 2x_1^2x_2^2 + 8x_1x_2x_3^2 + x_2^4 + 4x_3^4, \\
& x_1^4 - 2x_1^2x_3x_4 +\tfrac{1}{4}x_3^4+\tfrac{1}{2}x_3^2x_4^2+\tfrac{1}{4}x_4^4, \,\,
x_1^4 + 2x_1^2x_3x_4 +\tfrac{1}{4}x_3^4+\tfrac{1}{2}x_3^2x_4^2+\tfrac{1}{4}x_4^4,\\
& x_1^2 - 2x_1x_2 + x_2^2 + x_3^2 - 2x_3x_4 + x_4^2, \,\,
x_1^2 - 2x_1x_2 + x_2^2 + x_3^2 + 2x_3x_4 + x_4^2,\\
& x_1^2 + 2x_1x_2 + x_2^2 + x_3^2 - 2x_3x_4 + x_4^2, \,\,
x_1^2 + 2x_1x_2 + x_2^2 + x_3^2 + 2x_3x_4 + x_4^2.\\
\end{align*}

The element $x_1p_2+x_2p_3+x_3p_4$ lies in $\cc_{M_2}^\circ$ if and only if the
following polynomials are nonzero:

\begin{align*}
& x_1x_2x_3, \,\, x_1^2 - 2x_1x_2 + x_2^2 - x_3^2,\,\,
x_2^4-x_3^4, \,\, x_1^2 + 2x_1x_2 + x_2^2 - x_3^2,\\
& x_1^2 + x_2^2 - 2x_2x_3 + x_3^2, \,\,  x_1^2 + x_2^2 + 2x_2x_3 + x_3^2,\,\,
x_1^4 + 4x_3^4,\,\,  x_1^4 + 4x_2^4, \\
& x_1^4 - 2x_1^2x_2x_3 + \tfrac{1}{4}x_2^4 + \tfrac{1}{2}x_2^2x_3^2 + \tfrac{1}{4}x_3^4,\,\,
x_1^4 + 2x_1^2x_2x_3 + \tfrac{1}{4}x_2^4 + \tfrac{1}{2}x_2^2x_3^2 + \tfrac{1}{4}x_3^4.\\
\end{align*}  

The element $x_1p_2+x_2p_3$ lies in $\cc_{M_3}^\circ$ if and only if the following
polynomials are nonzero:

$$x_1x_2,\, x_1^4-x_2^4,\, x_1^8 + \tfrac{17}{4}x_1^4x_2^4 + x_2^8.$$

The element $x_1(p_2-p_3)+x_2(p_3-p_4)$ lies in $\cc_{M_4}^\circ$ if and only if
the following polynomials are nonzero:

$$x_1x_2,\, x_1^2 - 3x_1x_2 + 2x_2^2,\, x_1^4 + 4x_2^4,\,
x_1^2 - \tfrac{6}{5}x_1x_2 + \tfrac{2}{5}x_2^2,\,
x_1^2 - \tfrac{2}{5}x_1x_2 + \tfrac{2}{5}x_2^2.$$

The element $x_1p_1+x_2p_2$ lies in $\cc_{M_5}^\circ$ if and only if the following
polynomials are nonzero:

$$x_1x_2,\, x_1^4-x_2^4.$$

\section{Presentation and invariants of the small Weyl group}

For an explicit identification of $W_0$ with the group numbered 31 in the Shephard-Todd classification,
\cite{broue} can be used (Table A.3, page 129): it must have a presentation with 5 involutions
$s$, $t$, $u$, $v$, $w$ that obey the relations $sw=ws$, $uv=vu$, $svs=vsv$, $vtv=tvt$, $wtw=twt$, $wuw=uwu$ and $stu=tus=ust$.
In $W_0$ those can be chosen, for example, as follows:
\begin{itemize}
\item $s=s_1s_5s_3s_4s_3s_5s_1$, reflection in $(1,i,1+i,0)$;
\item $t=s_4$, reflection in $(1,1,1,1)$;
\item $u=s_2s_1s_5s_1s_2$, reflection in $(0,1+i,1,i)$;
\item $v=s_1$, reflection in $(1,0,0,0)$;
\item $w=s_4s_2s_3s_5s_4s_5s_3s_2s_4$, reflection in $(0,0,0,1)$.
\end{itemize}

Invariants of the Shephard-Todd group 31 have been determined in \cite{shephard} (page 287; see also \cite{orlik}, page 285),
based on the work of Maschke \cite{maschke}. For more details see \cite{dimca}. The algebra of invariants is polynomial
with generators of degrees 8, 12, 20, 24. While invariants of degrees 8 and 12 are determined uniquely up to scalars,
we note that there are some alternatives for the choice of generators of degrees 20 and 24.

With respect to our basis of the Cartan subspace, the group $W_0$ acts on polynomials in coefficients $x_1$, ..., $x_4$
of $x_1p_1+...+x_4p_4$. There are ten quadrics, corresponding to the classical Klein quadrics from the works cited above,
on which $W_0$ acts by permuting them up to scalar multiples. Let
\begin{align*}
Q_1&=x_1x_3+ix_1x_4-x_2x_4-ix_2x_3\\
Q_2&=x_1x_3+ix_2x_3-x_2x_4-ix_1x_4\\
Q_3&=x_1^2+ix_3^2-x_2^2-ix_4^2\\
Q_4&=x_1^2+ix_4^2-x_2^2-ix_3^2\\
Q_5&=x_1x_3-ix_1x_4-ix_2x_3+x_2x_4\\
Q_6&=x_1x_3+ix_1x_4+ix_2x_3+x_2x_4\\
Q_7&=x_1^2+x_2^2+2x_3x_4\\
Q_8&=2x_1x_2+x_3^2+x_4^2\\
Q_9&=x_1^2+x_2^2-2x_3x_4\\
Q_{10}&=2x_1x_2-x_3^2-x_4^2,
\end{align*}
then the above generating reflections $s_1$, ..., $s_5$ act on these quadrics in the following way:

\begin{center}
\begin{tabular}{c|rrrrr}
&$s_1$&$s_2$&$s_3$&$s_4$&$s_5$\\
\hline
$Q_1$&$-Q_6$&$iQ_2$&$Q_5$&$-(1+i)Q_4/2$&$(1+i)Q_{10}/2$\\
$Q_2$&$-Q_5$&$-iQ_1$&$-Q_6$&$(-1+i)Q_3/2$&$Q_2$\\
$Q_3$&$Q_3$&$-Q_4$&$Q_3$&$-(1+i)Q_2$&$Q_3$\\
$Q_4$&$Q_4$&$-Q_3$&$Q_4$&$(-1+i)Q_1$&$-Q_9$\\
$Q_5$&$-Q_2$&$iQ_6$&$Q_1$&$Q_5$&$Q_5$\\
$Q_6$&$-Q_1$&$-iQ_5$&$-Q_2$&$Q_6$&$-(1+i)Q_8/2$\\
$Q_7$&$Q_7$&$Q_7$&$-Q_9$&$Q_7$&$Q_7$\\
$Q_8$&$Q_8$&$Q_8$&$Q_8$&$Q_8$&$(-1+i)Q_6$\\
$Q_9$&$Q_9$&$Q_9$&$-Q_7$&$-Q_{10}$&$-Q_4$\\
$Q_{10}$&$Q_{10}$&$Q_{10}$&$Q_{10}$&$-Q_9$&$(1-i)Q_1$
\end{tabular}
\end{center}

It follows that the product $\Pi_{20}:=Q_1\cdots Q_{10}$ is invariant under the action of $W_0$.

Further, there are six fundamental quartics from \cite{maschke} which in our basis are
\begin{align*}
A_1&=2x_1^4+2x_2^4-x_3^4-x_4^4+12x_1x_2(x_3^2+x_4^2)+6x_3^2x_4^2\\
A_2&=-x_1^4-x_2^4+2x_3^4+2x_4^4+6x_1^2x_2^2-12(x_1^2+x_2^2)x_3x_4\\
A_3&=2x_1^4+2x_2^4-x_3^4-x_4^4-12x_1x_2(x_3^2+x_4^2)+6x_3^2x_4^2\\
A_4&=-x_1^4-x_2^4+2x_3^4+2x_4^4+6x_1^2x_2^2+12(x_1^2+x_2^2)x_3x_4\\
A_5&=-x_1^4-x_2^4-x_3^4-x_4^4-6x_1^2x_2^2+6i(x_1^2-x_2^2)(x_3^2-x_4^2)-6x_3^2x_4^2\\
A_6&=-x_1^4-x_2^4-x_3^4-x_4^4-6x_1^2x_2^2-6i(x_1^2-x_2^2)(x_3^2-x_4^2)-6x_3^2x_4^2
\end{align*}
with $A_1+...+A_6=0$ that are permuted by the action of the generating reflections as follows:

\begin{center}
\begin{tabular}{c|ccccc}
&$s_1$&$s_2$&$s_3$&$s_4$&$s_5$\\
\hline
$A_1$&$A_3$&$A_1$&$A_1$&$A_4$&$A_1$\\
$A_2$&$A_2$&$A_2$&$A_4$&$A_2$&$A_6$\\
$A_3$&$A_1$&$A_3$&$A_3$&$A_3$&$A_3$\\
$A_4$&$A_4$&$A_4$&$A_2$&$A_1$&$A_4$\\
$A_5$&$A_5$&$A_6$&$A_5$&$A_5$&$A_5$\\
$A_6$&$A_6$&$A_5$&$A_6$&$A_6$&$A_2$
\end{tabular}
\end{center}

Thus any symmetric function of $A_1$, ..., $A_6$ is also invariant under $W_0$. Denoting by $\sigma_k$ the $k$th elementary symmetric function of $A_1$, ..., $A_6$,
the choices made in \cite{orlik} are as follows. For the degree 8 invariant, $F_8=-\sigma_2/6$; for the degree 12 invariant, $F_{12}=-\sigma_3/4$;
for the degree 20 invariant, $F_{20}=\sigma_5/12$; and for the degree 24 invariant $F_{24}$, $1/265531392$ times the Hessian determinant of $F_8$.

Note that
\[
F_{20}=F_8F_{12}+81\Pi_{20}
\]
and
\[
F_{24}=\Pi_{24}-4F_{12}^2,
\]
where $\Pi_{24}=\sigma_6=A_1\cdots A_6$. Thus $F_8$, $F_{12}$, $\Pi_{20}$, $\Pi_{24}$ also are polynomial generators for the invariant ring.

Note also that the invariants from \cite{maschke,orlik,dimca} are defined over rationals.
One can obtain the same expressions using the change of variables
\begin{align*}
z_1&=\frac{x_1+x_2}{\sqrt2}\\
z_2&=-i\frac{x_3-x_4}{\sqrt2}\\
z_3&=i\frac{x_1-x_2}{\sqrt2}\\
z_4&=-i\frac{x_3+x_4}{\sqrt2},
\end{align*}
i. e. passing to the basis
\[
\frac{p_1+p_2}{\sqrt2}, i\frac{p_3+p_4}{\sqrt2}, -i\frac{p_1-p_2}{\sqrt2}, i\frac{p_3-p_4}{\sqrt2}
\]
of the Cartan subspace. In this basis,
\begin{align*}
Q_1&=z_1z_2+z_3z_4\\
Q_2&=z_1z_2-z_3z_4\\
Q_3&=z_1z_3+z_2z_4\\
Q_4&=z_1z_3-z_2z_4\\
Q_5&=z_1z_4+z_2z_3\\
Q_6&=z_1z_4-z_2z_3\\
Q_7&=z_1^2+z_2^2-z_3^2-z_4^2\\
Q_8&=z_1^2-z_2^2+z_3^2-z_4^2\\
Q_9&=z_1^2-z_2^2-z_3^2+z_4^2\\
Q_{10}&=z_1^2+z_2^2+z_3^2+z_4^2
\intertext{and}
A_1&=z_1^4+z_2^4+z_3^4+z_4^4-6(z_1^2z_2^2+z_1^2z_3^2+z_1^2z_4^2+z_2^2z_3^2+z_2^2z_4^2+z_3^2z_4^2)\\
A_2&=z_1^4+z_2^4+z_3^4+z_4^4-6(z_1^2z_2^2-z_1^2z_3^2-z_1^2z_4^2-z_2^2z_3^2-z_2^2z_4^2+z_3^2z_4^2)\\
A_3&=z_1^4+z_2^4+z_3^4+z_4^4-6(-z_1^2z_2^2+z_1^2z_3^2-z_1^2z_4^2-z_2^2z_3^2+z_2^2z_4^2-z_3^2z_4^2)\\
A_4&=z_1^4+z_2^4+z_3^4+z_4^4-6(-z_1^2z_2^2-z_1^2z_3^2+z_1^2z_4^2+z_2^2z_3^2-z_2^2z_4^2-z_3^2z_4^2)\\
A_5&=-2z_1^4-2z_2^4-2z_3^4-2z_4^4-24z_1z_2z_3z_4\\
A_6&=-2z_1^4-2z_2^4-2z_3^4-2z_4^4+24z_1z_2z_3z_4
\end{align*}

\section{Determining stabilizers of semisimple elements}\label{sec:stab}

From \cite[Corollary 3.13]{grale} we recall the following fact.

\begin{theorem}\label{thm:stab}
For $p,p'\in \cc_{M_i}^\circ$ we have $Z_{G_0}(p) = Z_{G_0}(p')$.
\end{theorem}

In this section we show how we determined the stabilizers $Z_{G_0}(p) =
\{ g\in G_0 \mid g(p)=p\}$ of the semisimple
elements $p$ in the sets $\cc_{M_i}^\circ$ listed in Table \ref{tab:semsim}.
By the previous theorem these are independent of the chosen element of
$\cc_{M_i}^\circ$. 

The Lie algebra of $Z_{G_0}(p)$ is the centralizer $\z_{\g_0}(p)$.
The latter can be explicitly calculated and determines the identity component
$Z_{G_0}(p)^\circ$. It remains to determine the component groups. For this we
want to find one explicit element (i.e., automorphism of $\g$) in each
component of $Z_{G_0}(p)$. We do this for $1\leq i\leq 9$, where $M_i$ is as
in Table \ref{tab:semsim}.

First we need to recall a number of facts. Let $\a$ be a semisimple complex
Lie algebra. Consider the root system of $\a$ with respect to a fixed Cartan
subalgebra. There are root vectors $e_1,\ldots,e_\ell$ (corresponding to the
simple positive roots), $f_1,\ldots,f_\ell$ (corresponding to the negative
simple roots) and $h_1,\ldots,h_\ell$ in the Cartan subalgebra such that
$$[h_i,h_j]=0,\,[e_i,f_j] = \delta_{ij} h_i,\, [h_j,e_i] = C(i,j)e_i,\,
[h_j,f_i]=-C(i,j)f_i \text{ for } 1\leq i,j\leq\ell,$$
where $C$ is the Cartan matrix of the root system. These elements generate
$\a$ and are called a canonical set of generators of $\a$. 
Let $\pi$ be a permutation of $\{1,\ldots,\ell\}$ such that $C(i,j)=
C(\pi(i),\pi(j))$ for all $i,j$. Then mapping $e_i\mapsto e_{\pi(i)}$,
$f_i\mapsto f_{\pi(i)}$, $h_i\mapsto h_{\pi(i)}$ extends to a unique automorphism
$\sigma_\pi$ of $\a$ (cf. \cite[Theorem IV.3]{jac}). Here we call
$\sigma_\pi$ a {\em pure diagram automorphism} of $\a$ (with respect to the
fixed choice of a canonical generating set). Let $\Gamma$ be
the group of all pure diagram automorphisms. Let $\Int(\a)$ be the inner
automorphism group of $\a$ (this is the algebraic subgroup of $\GL(\a)$
generated by $\exp(\ad x)$ for all nilpotent $x\in \a$). Then $\Int(\a)$
is the identity component of the automorphism group of $\a$. Moreover, we
have $\Aut(\a) = \Gamma\ltimes \Int(\a)$ (cf. \cite[VIII.5 no 3, Cor 1]{bou3},
\cite[\S IX.4]{jac}).

Now let $A\subset \GL(V)$ be an algebraic group with Lie algebra $\a\subset
\gl(V)$. Suppose that $\a$ is semisimple. Let $g\in A$; then $\Ad(g) :
\a\to \a$ with $\Ad(g)(x) =gxg^{-1}$ is an automorphism of $\a$. So we get a 
homomorphism $\Ad : A\to \Aut(\a)$. We have that $\Ad : A^\circ \to \Int(\a)$ is
surjective. Let $g\in A$ then we can write $\Ad(g) = \sigma_\pi \tau$
with $\sigma_\pi\in \Gamma$ and $\tau\in \Int(\a)$. Let $g_\tau$ be a preimage
of $\tau$ in $A^\circ$.
Then $\Ad(gg_\tau^{-1})=\sigma_\pi$. It follows that every component of $A$
contains an element whose image under $\Ad$ is a pure diagram automorphism.
In our situation we always have that $A$ is a subgroup of $G=\Aut(\g)$.
The Lie algebra of $G$ is $\ad \g = \{ \ad x \mid x\in \g\}$. So
$\a = \{ \ad x \mid x\in \tilde\a\}$ where $\tilde\a$ is a subalgebra of $\g$.
For $g\in G$, $x\in \g$ we have $\Ad(g)(\ad x) = \ad g(x)$. It follows that
instead of the adjoint action of $A$ on $\a\subset \ad \g$ we can also work with
directly with the action of $A$ on $\tilde \a$. 

Let $A$, $\a$ be as above. We fix a Cartan subalgebra of $\a$ and corresponding
root system $\Psi$ with a fixed set of simple roots $\beta_1,\ldots,\beta_\ell$.
We fix a Chevalley basis of $\a$ consisting of
root vectors $x_\beta$ for $\beta\in \Psi$ and $h_1,\ldots,h_\ell$ in the Cartan
subalgebra (cf. \cite[Theorem 25.2]{hum}). For $\beta\in \Psi$ define the
elements
\begin{align*}
  x_\beta(t) &= \exp( tx_\beta),\, (t\in \C),\\
  w_\beta(t) &= x_\beta(t)x_{-\beta}(-t^{-1})x_\beta(t),\, (t\in \C^*),\\
  h_\beta(t) & = w_\beta(t)w_\beta(1)^{-1}, \,(t\in \C^*).
\end{align*}
Let $w$ be an element of the Weyl group of $\Psi$ and let $w=s_{\beta_{i_1}}\cdots
s_{\beta_{i_r}}$ be a reduced expression (so the $\beta_{i_k}$ are simple roots).
Then we define $\dot{w} = w_{\beta_{i_1}}(1)\cdots w_{\beta_{i_r}}(1)$. Also let
$\Psi_w$ be the set of positive roots $\beta\in \Psi$ such that $w(\beta)$ is
a negative root. Let $\beta_1,\ldots,\beta_m$ be the positive roots of
$\Psi$ and write $\Psi_w = \{\gamma_1,\ldots,\gamma_q\}$. Then every element of
$A^\circ$ can be written as
\begin{equation}\label{eq:bruhat}
x_{\beta_1}(u_1)\cdots x_{\beta_m}(u_m) h_{\beta_1}(t_1)\cdots h_{\beta_\ell}(t_\ell)
\dot{w} x_{\gamma_1}(s_1)\cdots x_{\gamma_q}(s_q)
\end{equation}
where $u_i,s_j\in \C$, $t_k\in \C^*$ and $w$ runs over the Weyl group of $\Psi$
(cf. \cite[Theorem 5.2.23]{gra16}, \cite[Corollary 8.3.9]{springer}).
For $w\in W$ let $C_w$ denote the set of all elements of the form
\eqref{eq:bruhat}. Then $A^\circ$ is the disjoint union of the sets $C_w$.
This is called the Bruhat decomposition of $A^\circ$ and the $C_w$ are called
the Bruhat cells of $A^\circ$. For our purposes this way of writing elements of
$A^\circ$ is useful because it gives a parametrization of $A^\circ$. Using it we
can use Gr\"obner basis techniques for finding elements with certain
properties of $A$ (for example those that stabilize a given $p\in \cc$).

We also remark here that there exists an algorithm that given a {\em connected}
algebraic
group $A\subset \GL(V)$ with Lie algebra $\a\subset \gl(V)$ decides whether a
given $g\in \GL(V)$ lies in $A$, taking as input the element $g$ and a basis of
$\a$ (see \cite[Remark 5.8]{borwdg}). When computing component groups this
algorithm is very useful, as it allows to decide whether a given element lies
in the identity component of an algebraic group, or whether two elements lie
in the same component of the group.

In the explanations below many statements come from an explicit computation
in {\sf GAP}4. If this is the case then we add a ({\sf GAP}) to the statement.
On some occasions we have also used the computer algebra system {\sc Magma}
(\cite{magma}) for Gr\"obner basis computations. 

Now we turn to the task of finding the component groups of the stabilizers of
$p\in \cc_{M_i}^\circ$. For each case we write a paragraph.

{\bf Let $p\in \cc_{M_1}^\circ$.} Then $W_p=M_1$ which is trivial.
The centralizer of $p$ in $\g$ is the Cartan subalgebra $\h = \z_{\g}(\cc)$ of
$\g$.  We have $\h = \h_1\oplus \h_3$, where $\h_i = \g_i
\cap \h$ and both intersections are of dimension 4. Let $g\in G_0$
stabilize $p$ then it also stabilizes $\h$ and hence it induces an element $w$
of the
Weyl group $W = N_G(\h)/Z_G(\h)$. As $g\in G_0$ it commutes with $\theta$.
Hence $w\in W^\theta$. As seen in the previous section this means that the
restriction of $w$ to $\cc$ lies in $W_0$. Since $p$ has
trivial stabilizer in $W_0$ and only the identity in $W^\theta$ restricts to the
identity in $W_0$, it follows that a $g\in Z_{G}(\h)$.
Hence the stabilizer of $p$ in $G_0$ is $G_0\cap Z_G(\h)$. Since $G$ is
simply connected a theorem of Steinberg (\cite[Theorem 8.1]{steinberg68})
states that $G_0=G^\theta$ (the centralizer of $\theta$). Another theorem of
Steinberg (\cite[Corollary 3.11]{steinberg75}) states that
the group $Z_G(\h)$ is connected. Its Lie algebra is $\z_\g(\h) = \h$.
As shown in \cite[\S 6.5]{borwdg} we can explicitly compute an isomorphism
of algebraic groups $\lambda : (\C^*)^8 \to Z_G(\h)$. Then the condition
$\theta \lambda(t_1,\ldots,t_8) = \lambda(t_1,\ldots,t_8)\theta$ is equivalent
to a set of polynomial equations in $t_1,\ldots,t_8$ and their inverses.
In these equations we write $s_i$ in place of $t_i^{-1}$ and add the equation
$s_it_i=1$. The zero locus ({\sf GAP}) of the resulting polynomials gives an
elementary abelian 2-group of order 16.

The latter group has been previously also determined by A.~M.~Popov, see entry 13 of Table 1 in \cite{popov}.

{\bf Let $p\in \cc_{M_2}^\circ$.} Then $Z_{G_0}(p)^\circ$ is a
1-dimensional torus ({\sf GAP}). Let $q\in\g_0$ span
$\z_{\g_0}(p)$. The centralizer in $\g$ of the subalgebra spanned by $p,q$ is
a Cartan subalgebra $\hat\h$ of $\g$ ({\sf GAP}). We remark that $\hat\h$ does not
contain $\cc$, and is therefore not equal to $\h$.
Let $g\in G_0$ satisfy $g(p)=p$. Then $g(\hat\h)=\hat\h$. Hence $g$ induces an
element of the Weyl group $\widehat{W}=N_G(\hat\h)/Z_G(\hat\h)$. We have that
$\hat\h$
is $\theta$-stable; hence $\theta$ also induces an element of the same Weyl
group. Let $\widehat{\Phi}$ denote the root system of
$\g$ with respect to $\hat\h$. By $g$, $\theta$
we also denote the elements of $\widehat{W}$ induced by $g\in G_0$ and
$\theta\in G$.
As $g\in G_0$ we have that $g$ commutes with $\theta$. By writing elements of
$\widehat{W}$ as permutations of the roots of $\widehat{\Phi}$ and using
permutation group
algorithms we can compute the centralizer $\widehat{W}^\theta$ of $\theta$ in
$\widehat{W}$.
It turns out to have order 768 ({\sf GAP}). For $w\in \widehat{W}$ we denote an element of
$G$
inducing it by $\dot{w}$. A computation ({\sf GAP}) shows that there are exactly two
elements in $\widehat{W}^\theta$ that stabilize $p$, and let $w$ be one of
these elements.
We want to find the elements of $G_0$ that induce $w$. 
Let $\widehat{H}$ be the connected subgroup of $G$ with Lie algebra $\hat\h$
(or, more precisely, $\ad\hat\h$). Then $\widehat{H}=Z_G(\hat\h)$ and
$\dot{w}\widehat{H}$ is precisely
the set of elements of $G$ that stabilize $\hat\h$ and induce $w$. Since
$\widehat{H}$ is a connected torus, in the same way as in the previous case,
we can parametrize its elements with eight nonzero parameters.
As $G_0=G^\theta$ we have that
$G_0\cap \dot{w}\widehat{H}$ is exactly the set of elements of
$\dot{w}\widehat{H}$ commuting
with $\theta$. The condition that an element commutes with $\theta$ translates
to polynomial equations in the eight parameters. As in the previous case
we introduce extra indeterminates for the inverses of the parameters. The
zero locus of the resulting polynomials is a variety
of dimension 1. To reduce this dimension we consider the identity component of
$Z_{G_0}(p)$. This is a 1-dimensional torus whose Lie algebra is spanned by
$\ad q$. Let $\alpha\in \widehat{\Phi}$ and let $y_\alpha$ be a corresponding
root vector. Let $y_\beta = \dot{w}(y_\alpha)$; here $\beta = w(\alpha)$. 
Then for $h\in \dot{w}\widehat{H}$ we have $h(y_\alpha) = c_h y_\beta$ with
$c_h\in \C$.
Let $t\in Z_{G_0}(p)^\circ$ then $ht(y_\alpha) = \alpha(t)c_h y_\beta$. So if we
select $\alpha$ such that $Z_{G_0}(p)^\circ$ acts nontrivially on $y_\alpha$ there
is a $t$ such that $\alpha(t)c_h=1$. It follows that every component of
$Z_{G_0}(p)$ contains an element $h$ with $c_h=1$. This requirement yields one
more polynomial equation. The resulting zero-locus is 0-dimensional. It turns
out that in total we get eight solutions ({\sf GAP}) and the component group
is elementary abelian of order 8. 

{\bf Let $p \in \cc_{M_3}^\circ$.} Then $Z_{G_0}(p)^\circ$ is a
2-dimensional torus ({\sf GAP}). Let  $\z_1 = \z_{\g_0}(p)$ (which is toral of 
dimension 2) and $\z_2=\z_{\g_0}(\z_1)$, which is a Cartan subalgebra of
$\g_0$ ({\sf GAP}). An element of $Z_{G_0}(p)$ stabilizes $\z_1$ and $\z_2$; hence it
induces an element of the Weyl group of $\g_0$ with respect to $\z_2$.
We run through this Weyl group and select all elements that stabilize $\z_1$;
there are 8 of them ({\sf GAP}). Then for each
such element $w$ we check whether there is a $z\in Z_{G_0}(\z_2)$ such that
$zw(p)=p$; note that $Z_{G_0}(\z_2)$ is a connected torus whose Lie algebra
is $\ad \z_2$, so we can compute an isomorphism $\lambda : (\C^*)^8 \to
Z_{G_0}(\z_2)$ and the check reduces to a Gr\"obner basis computation.
There are four elements $w$ such that there is a $z\in Z_{G_0}(\z_2)$ with
$zw(p)=p$ ({\sf GAP}). Each of them gives exactly one element of the component group ({\sf GAP})
which is elementary abelian of order 4. 

{\bf Let $p\in \cc_{M_4}^\circ$.} Then $Z_{G_0}(p)^\circ$ is simple of type $A_1$
({\sf GAP}). The component
group of $Z_{G_0}(p)$ is an elementary abelian 2-group of size 4. This is shown
in the following way. Since Lie algebras of type $A_1$ have no outer
automorphisms,
every component of $Z_{G_0}(p)$ has an element that is the identity on
$\z_1=\z_{\g_0}(p)$. Such elements lie in $Z_{G_0}(\z_1)$. We first determine the
component group of the latter group. Let $\z_2=
\z_{\g_0}(\z_1)$ which is the Lie algebra of $Z_{G_0}(\z_1)$.
Then $\z_2$ is semisimple of type $A_1+A_1$  ({\sf GAP}). Every component of
$Z_{G_0}(\z_1)$ has an element $g$ such that restricts to 
a pure diagram automorphism of $\z_2$. Consider the subalgebra
$\z = \z_1+\z_2$. The $\z$-module
$\g_1$ splits as a direct sum of four irreducible modules with highest weights
$(0;31)$, $(2;11)$, $(4;11)$, $(2;31)$  ({\sf GAP}).
(Here the weights of $\z_1$ and
$\z_2$ are separated by a semicolon.) We denote the corresponding highest
weight vectors by $v_1,\ldots,v_4$. All four weights are different, hence
this decomposition
is unique. Let $g\in Z_{G_0}(\z_1)$ restrict to a pure diagram
automorphism of
$\z_2$. If the diagram automorphism would be nontrivial then $g$ would map
$v_1$ to
a highest weight vector of weight $(0;31)$. But there is no such highest
weight. It follows that $g$ is the identity on $\z_2$
(and hence on $\z$). Furthermore $g$ must map the highest weight vectors of
the above modules to nonzero scalar multiples of themselves. The subalgebra
$\z$ along with the four highest weight vectors generate $\g$  ({\sf GAP}).
We now
define automorphisms of $\g$ by requiring that they are the identity on
$\z$ and map $v_i\mapsto z_i v_i$, $1\leq i\leq 4$. Imposing the condition that
this defines an automorphism yields polynomial equations on the $z_i$.
It turns out that there are exactly four solutions  ({\sf GAP});
which are representatives
of the elements of the component group of $Z_{G_0}(\z_1)$.
This group is cyclic of order 4.
Let $h_0\in G_0$ denote (a representative of) a generator of this component
group. Then $Z_{G_0}(\z_1)$ is the disjoint union of the sets
$h_0^i Z_{G_0}(\z_1)^\circ$. By using the Bruhat decomposition of the semisimple
group $Z_{G_0}(\z_1)^\circ$ we find polynomial equations for the set of elements
of each of these sets that stabilize $p$  ({\sf GAP}).
It turns out that for $i=1,3$ there
are no elements stabilizing $p$, whereas for $i=0,2$ there are. By Gr\"obner
basis techniques we compute all these elements and find the component group.

{\bf Let $p\in \cc_{M_5}^\circ$.} Then $Z_{G_0}(p)^\circ$ is of type $A_1+T_3$
({\sf GAP}).
The component group of $Z_{G_0}(p)$ is of order 4. We have computed it in the
following way. The Lie algebra of $Z_{G_0}(p)$ is $\z_1=\z_{\g_0}(p)$. Let
$\z_1'$ be its derived subalgebra, which is simple of type $A_1$  ({\sf GAP}).
As $\z_1'$ has no outer automorphisms each component of $Z_{G_0}(p)$ has
an element which restricts to the identity on $\z_1'$. These elements lie
in $Z_2=Z_{G_0}(\z_1')$. First we show that this group is connected.
The Lie algebra of $Z_2$ is $\z_2 = \z_{\g_0}(\z_1')$. This is a semisimple
subalgebra of $\g_0$ of type $A_1+2A_3$  ({\sf GAP}). Each component of $Z_2$
has an
element that restricts to a pure diagram automorphism of $\z_2$.
Set $\z = \z_1'\oplus \z_2$. Let $V$ be the orthogonal complement of
$\z$ in $\g_0$ with respect to the Killing form of $\g$. Then $V$ is an
irreducible $\z$-module of highest weight $(1; 1; 010; 000 )$  ({\sf GAP}).
(Here we
enumerate the Dynkin diagram of $\z$ as follows: first the $A_1$ corresponding
to $\z_1'$, then the $A_1$ in $\z_2$, then the two $A_3$'s.) Let $g\in Z_2$
restrict to a pure diagram automorphism of $\z_2$; then $g$ maps $V$ to itself.
Let $v_1$ be a fixed highest-weight vector of $V$. Since $g$ permutes the
elements of the canonical generating set of $\z$ we must have that $g(v_1)$
is a multiple of $v_1$. This implies that the diagram automorphism induced
by $g$ cannot interchange the two $A_3$'s. Now we let $U$ be the space $\g_1$
viewed as $\z$-module. It splits into a direct sum of two irreducible
modules $U_1$, $U_2$ with highest weights $(1;0;010;000)$ and $(0;1;001;100)$
 ({\sf GAP}).
Hence this decomposition is uniquely determined. We have that $g(U_i)$ is an
irreducible $\z$-submodule of $U$. But because of the weights of the first
$A_1$, $g$ cannot interchange the two modules. Let $u_1$, $u_2$ be highest
weight vectors; then $g(u_i)$ is a nonzero scalar multiple of $u_i$.
Since the $A_3$-parts of the highest weights are not invariant under
the diagram automorphisms of these $A_3$'s, we see that the restriction of
$g$ to $\z$ cannot be a nontrivial diagram automorphism. In other words, it
must be the identity. Now $\z$ along with $v_1$, $u_1$, $u_2$ generate $\g$
 ({\sf GAP}).
So we consider the set of all automorphisms of $\g$ that restrict to the
identity on $\z$ and map $v_1\mapsto z_1v_1$, $u_1\mapsto z_2u_1$, $u_2\mapsto
z_3u_2$. This set of automorphisms corresponds to the solution set of
a set of polynomial equations in $z_1,z_2,z_3$. It turns out that there are
8 solutions  ({\sf GAP}). So we get a group of 8 automorphisms that satisfy
the above
requirements. However, it turns out that they all lie in the identity
component of $Z_2$  ({\sf GAP}). We conclude that $Z_2$ is connected.

Now we consider the stabilizer of $p$ in $Z_2$. We parametrize the group $Z_2$
using the Bruhat decomposition  ({\sf GAP}).
The Weyl group of $\z_2$ has 1152 elements,
hence there are 1152 cells to consider. By Gr\"obner basis computations we
established that only 4 cells contain elements that stabilize $p$. Let
$T$ denote the connected subgroup of $G_0$ corresponding to the centre of
$\z_1$. It is a 3-dimensional torus. This group is contained in $Z_2$,
hence for each cell the set of elements that stabilize $p$ is a 3-dimensional
variety (if non-empty). We can compute an explicit isomorphism $\lambda :
(\C^*)^3\to T$ ({\sf GAP}). Using $\lambda$ we can divide the solutions
into cosets of $T$. This yields a finite number of elements to consider.
In the end, modulo the identity component, each of the four cells
with elements that stabilize $p$ gives rise to exactly one element of the
component group, which therefore has order 4.

{\bf Let $p\in \cc_{M_6}^\circ$.} Then $Z_{G_0}(p)^\circ$ is of type $A_1+T_1$
({\sf GAP}).
The component group is of order 2. This is establised in the following way.
We set $\z_1' = [\z_1,\z_1]$, where $\z_1= \z_{\g_0}(p)$. The procedure is highly
analogous to the case where $p\in \cc_{M_4}^\circ$, but we use $\z_1'$ instead of $\z_1$. 
We let $\z_2$ denote the centralizer of $\z_1'$ in $\g_0$.
It is of type $2A_1$  ({\sf GAP}). We set $\z=\z_1'\oplus \z_2$. We are
interested in the
component group of $Z_{G_0}(\z_1')$. Every component of the latter group has
an element that restricts to a pure diagram automorphism of $\z_2$ and such
that its restriction to $\z_1'$ is the identity. As $\z$-module $\g_1$ has
highest weights $(0;31)$, $(2;11)$, $(4;11)$, $(2;31)$ ({\sf GAP}).
So again, an element of $Z_{G_0}(\z_1')$ cannot act
on $\z_2$ as a non-trivial pure diagram automorphism. Again we denote the
highest weight vectors by $v_1,\ldots,v_4$.
The elements of $G$ that are the identity
on $\z$ and map $v_i\mapsto z_iv_i$ for $1\leq i\leq 4$ form a cyclic
group of order 4  ({\sf GAP}). However, the square of a generator lies in
$Z_{G_0}(\z_1')^\circ$  ({\sf GAP}). So the component group of
$Z_{G_0}(\z_1')$ is of order 2.
Let $h_0$ denote a nontrivial element.
Using the Bruhat decomposition of $Z_{G_0}(\z_1')^\circ$
we find the set $U$ of elements $g$ of that group with $g(p)=p$  ({\sf GAP}).
In this case
every solution set has dimension 1 because of the following reason.
Write $A=Z_{G_0}(p)\cap Z_{G_0}(\z_1')^\circ$. Then $\a = \z_{\g_0}(p)\cap
\z_{\g_0}(\z_1')$ is its Lie algebra. We have that $\a$ is the 1-dimensional
centre of $\z_{\g_0}(p)$. Hence the solution set of our equations is
1-dimensional. Actually, the Bruhat decomposition of $Z_{G_0}(\z_1')^\circ$ has
four cells, two of which have empty intersection with $U$ and the intersection
of the other two with $U$ is a variety of dimension 1. We can compute
an isomorphism $\lambda : \C^*\to A^\circ$  ({\sf GAP}).
We use this to partition $U$ into
$A^\circ$-cosets, of which there are two. One of the two representatives of
these cosets lies in $Z_{G_0}(p)^\circ$, whereas the other does not. 
Secondly, by Bruhat decomposition
again we parametrize the set of $g$ in $Z_{G_0}(\z_1)^\circ$ with $g(p)=h_0(p)$.
For this there turns out to be no solution  ({\sf GAP}).
So the component group of $Z_{G_0}(p)$ is of order 2.

{\bf Let $p\in \cc_{M_7}^\circ$.} Then $Z_{G_0}(p)^\circ$ is of type $2A_1$
({\sf GAP}), and the
component group of $Z_{G_0}(p)$ is of order 2. In order to see this
write $\z_1 = \z_{\g_0}(p)$ then $\z_1$ is semisimple of type $2A_1$  ({\sf GAP}).
As $\z_1$-module $\g_1$ splits as a direct sum of 14 irreducible submodules.
Among these there is a unique module with highest weight $(3,1)$ and a
unique module with highest weight $(1,3)$  ({\sf GAP}). Let $v_1,v_2$ be
corresponding
highest weight vectors. A computation shows that $\g$ is generated by
$\z_1$, $p$, $v_1$, $v_2$  ({\sf GAP}). Every component of $Z_{G_0}(p)$
contains an
element that restricts to either the identity or a pure diagram automorphism
of $\z_1$. Furthermore if $g\in Z_{G_0}(p)$ restricts to the identity on
$\z_1$ then $g(v_1) = z_1v_1$, $g(v_2)=z_2v_2$. All automorphisms of
$\g$ satisfying these conditions correspond to the solution set of a
set of polynomial equations in $z_1,z_2$.
It turns out that there are two solutions and both yield elements
that lie in $Z_{G_0}(p)^\circ$  ({\sf GAP}). Next we consider the elements that
restrict to a
pure diagram automorphism of $\z_1$. In this case we have $g(v_1) = z_1v_2$,
$g(v_2)=z_2v_1$. We do the same thing, and again get two solutions yielding
two automorphisms of order 2  ({\sf GAP}). They do not lie in the identity
component, but
are equal modulo the identity component  ({\sf GAP}).
Hence only one element of order 2 remains.

{\bf Let $p\in \cc_{M_8}^\circ$.} Then $Z_{G_0}(p)^\circ$ is of type $2A_1+A_2+T_1$
 ({\sf GAP}).
The stabilizer $Z_{G_0}(p)$ has two components. This is very
similar to the previous case. Here we work with the derived algebra $\z_1'$
of $\z_1=\z_{\g_0}(p)$. It is semisimple, we enumerate its Dynkin diagram
as follows: first the two $A_1's$, then the $A_2$. Now in $\g_1$ there are
unique $\z_1'$-submodules of highest weights $(2;0;01)$, $(0;2;10)$,
$(1;1;10)$, $(1;1;01)$  ({\sf GAP}). We see that the only possible diagram
automorphism is the one that
simultaneously interchanges the two $A_1's$ and reverses the Dynkin diagram of
$A_2$. We denote the four highest weight vectors by $v_1,\ldots,v_4$. 
Then $\g$ is generated by $\z_1'$, $p$, $v_1,\ldots,v_4$
({\sf GAP}).
Now we proceed exactly as in the previous case. However, due to the
presence of the 1-dimensional centre, the solution set of our polynomial
equations is also 1-dimensional. Let $T$ denote the connected algebraic subgroup
of $G$ whose Lie algebra is the centre of $\z_1$ (or, more precisely,
$ad \z_1$). We can compute an explicit isomorphism $\lambda : \C^*\to T$
({\sf GAP}).
Explicit computation shows that $\lambda(t)(v_1) = t^{-2}v_1$  ({\sf GAP}).
Let $g\in Z_{G_0}(p)$ restrict to the identity on $\z_1'$ and map $v_i$ to
$z_i v_i$ for $1\leq i\leq 4$, or restrict to the only possible nontrivial diagram automorphism
of $\z_1'$ and such that $g(v_1) = z_1v_2$, $g(v_2)=z_2v_1$, $g(v_3) = z_3v_4$,
$g(v_4)=z_4v_3$. Then by multiplying $g$ by a suitable element of the form
$\lambda(t)$ we find an element of the same component of $Z_{G_0}(p)$ and
such that $z_1=1$. This extra condition makes the solution set
finite. Similarly to the previous case, we find two elements that restrict
to the identity, and two elements that restrict to a pure diagram automorphism
 ({\sf GAP}).
The former elements lie in the identity component, whereas the latter are
equal modulo the identity component. So also here we have a component
group of order 2.

{\bf Let  $p\in \cc_{M_9}^\circ$.} We have $\cc_{M_9}^\circ=0$
so that $Z_{G_0}(p) = G_0$ which is connected.

\section{The mixed orbits}\label{sec:mixed}

In this section we determine the orbits of mixed type. Such an orbit has a
representative of the form $p+e$ where $p$ is semisimple, $e$ is nilpotent
$[p,e]=0$ and $p,e$ are both nonzero. A first remark is that $e$ lies in the
graded subalgebra $\z_{\g}(p)$. Secondly, we may assume that $p$ lies in one
of the classes $\cc_{M_i}^\circ$ listed in Table \ref{tab:semsim}. We have that
$p+e$ and $p+e'$ are $G_0$-conjugate if and only if $e,e'$ are
$Z_{G_0}(p)$-conjugate. By Theorem \ref{thm:stab} the stabilizer $Z_{G_0}(p)$
does not depend on the choice of the point $p$ in $\cc_{M_i}^\circ$. 
Furthermore, in the previous section the component groups of these stabilizers
have been determined. Using the algorithms of \cite{gra15} we can determine
the nilpotent $Z_{G_0}(p)^\circ$-orbits in $\z_{\g}(p)_1$. This yields a finite set
of representatives. Some of them are conjugate under representatives of the
component group of $Z_{G_0}(p)$. Since we have determined those representatives
explictly, we can decide this and obtain an irredundant list of nilpotent parts
of the mixed elements with semisimple part from $\cc_{M_i}(p)$.

Below we list the nilpotent parts of the representatives of the mixed orbits
with semisimple part in $\cc_{M_i}^\circ$ for $1\leq i\leq 9$. In each table the
second column has the dimension of the given nilpotent orbit in $\z_\g(p)_1$;
this is the same as the dimension of the space $[\z_{g}(p)_0,e]$. The third
column has the isomorphism type of the centralizer $\z_\g(p+e)$. Here we use
the following notation: $\ttt_k$ means a toral (i.e., commutative and
consisting of semisimple
elements) subalgebra of dimension $k$; $\u_k$ indicates an ideal consisting of
nilpotent elements of dimension $k$. A semisimple subalgebra is indicated by
the type of its root system.

If $p\in \cc_{M_1}^\circ$ then $\z_\g(p)_1=\cc$, so there are no nilpotent
elements in $\g_1$ that are centralized by $p$. Hence in this case we
do not obtain any mixed elements.

For $p\in \cc_{M_2}^\circ$ there are two nilpotent
$Z_{G_0}(p)^\circ$-orbits in $\z_{\g}(p)_1$. They are conjugate under the
component group of $Z_{G_0}(p)$, so only one orbit remains. It is shown in
Table \ref{tab:cm2}.

\begin{table}[htb]\caption{Nilpotent parts of mixed elements with semisimple
    part in $\cc_{M_2}^\circ$.}\label{tab:cm2}
\begin{tabular}{|l|r|r|}
  \hline
  element $e$ & $\dim$ & $\z_{\g_0}(p+e)$ \\
  \hline
  $(3,5)\otimes 1 +(1,3)\otimes 4$ & 1 & 0\\
  \hline
\end{tabular}
\end{table}

For $p\in \cc_{M_3}^\circ$ there are eight nilpotent
$Z_{G_0}(p)^\circ$-orbits in $\z_{\g}(p)_1$. The action of the component group
reduces this number to three. Table \ref{tab:cm3} lists their representatives.

\begin{table}[htb]\caption{Nilpotent parts of mixed elements with semisimple
    part in $\cc_{M_3}^\circ$.}\label{tab:cm3}
\begin{tabular}{|l|r|r|}
  \hline
  element $e$ & $\dim$ & $\z_{\g_0}(p+e)$ \\
  \hline
$(1,4)\otimes 1-(1,5)\otimes 3$ & 1 & $\ttt_1$ \\
$(3,5)\otimes 1+(1,3)\otimes 4$ & 1 &  $\ttt_1$ \\
$(1,4)\otimes 1-(3,5)\otimes 1-(1,5)\otimes 3-(1,3)\otimes 4$ & 2 & 0\\
  \hline
\end{tabular}
\end{table}

For $p\in \cc_{M_4}^\circ$ there are two nilpotent
$Z_{G_0}(p)^\circ$-orbits in $\z_{\g}(p)_1$. The action of the component group
is trivial on these orbits. Table \ref{tab:cm4} lists their representatives.

\begin{table}[htb]\caption{Nilpotent parts of mixed elements with semisimple
    part in $\cc_{M_4}^\circ$.}\label{tab:cm4}
\begin{tabular}{|l|r|r|}
  \hline
  element $e$ & $\dim$ & $\z_{\g_0}(p+e)$ \\
  \hline
$(3,5)\otimes 1 + (1,3)\otimes 4$ & 2 & $\ttt_1$\\
  $()\otimes 1+(2,3)\otimes 1+(1,3,4,5)\otimes 1-(3,5)\otimes 2$ & 3 & 0 \\
  $+(1,3)\otimes 3+(1,5)\otimes 4-(3,4)\otimes 4-(1,2,3,5)\otimes 4$ & & \\ 
  \hline
\end{tabular}
\end{table}

For $p\in \cc_{M_5}^\circ$ there are fourtyone nilpotent
$Z_{G_0}(p)^\circ$-orbits in $\z_{\g}(p)_1$. Up to the action of the component group
thirteen orbits remain. Table \ref{tab:cm6} lists their representatives.

\begin{table}[htb]\caption{Nilpotent parts of mixed elements with semisimple
    part in $\cc_{M_5}^\circ$.}\label{tab:cm6}
\begin{tabular}{|l|r|r|}
  \hline
  element $e$ & $\dim$ & $\z_{\g_0}(p+e)$ \\
  \hline
  $(1,4)\otimes 1$ & 2 & $\ttt_3+\u_1$ \\
  $(1,4)\otimes 1-(4,5)\otimes 4$ & 3 &  $\ttt_2+\u_1$ \\
  $(1,5)\otimes 3+(4,5)\otimes 4$ & 3 &  $\ttt_2+\u_1$ \\
  $()\otimes 2-(4,5)\otimes 4$ & 3 &  $\ttt_2+\u_1$ \\
  $(2,3)\otimes 2-(4,5)\otimes 4$ & 4 & $\ttt_2$ \\
  $(2,3)\otimes 2-(1,5)\otimes 3$ & 4 & $\ttt_2$\\
  $(1,4)\otimes 1+(2,3)\otimes 2$ & 4 & $\ttt_2$\\
  $(1,4)\otimes 1+()\otimes 2-(4,5)\otimes 4$ & 4 & $\ttt_1+\u_1$\\
  $(2,3)\otimes 2-(1,5)\otimes 3-(4,5)\otimes 4$ & 5 & $\ttt_1$\\
  $(1,4)\otimes 1+(2,3)\otimes 2-(4,5)\otimes 4$ & 5 & $\ttt_1$\\
  $(1,4)\otimes 1+(2,3)\otimes 2-(1,5)\otimes 3$ & 5 & $\ttt_1$\\
  $()\otimes 2 -(4,5)\otimes 4 +(1,4)\otimes 1 +(1,2,3,5)\otimes 3$ & 6 &0\\
  $(1,4)\otimes 1+()\otimes 2+(2,3)\otimes 2+(2,3,4,5)\otimes 4$
  & 6 &0\\
  \hline
\end{tabular}
\end{table}

For $p\in \cc_{M_6}^\circ$ there are eight nilpotent
$Z_{G_0}(p)^\circ$-orbits in $\z_{\g}(p)_1$. Up to the action of the component group
five orbits remain. Table \ref{tab:cm5} lists their representatives.

\begin{table}[htb]\caption{Nilpotent parts of mixed elements with semisimple
    part in $\cc_{M_6}^\circ$.}\label{tab:cm5}
\begin{tabular}{|l|r|r|}
  \hline
  element $e$ & $\dim$ & $\z_{\g_0}(p+e)$ \\
  \hline
$(2,3)\otimes 2 +(2,3,4,5)\otimes 4$ & 1 & $A_1$\\

$-(2,5)\otimes 1+(3,5)\otimes 1-(1,2,4,5)\otimes 2+(1,3,4,5)\otimes 2$ 
  & 2 & $\ttt_1+\u_1$\\

$-(2,4)\otimes 1+(3,4)\otimes 1-2(1,2,3,5)\otimes 1-(1,2)\otimes 2+(1,3)\otimes 2-2(4,5)\otimes 2$ &  3 & $\ttt_1$ \\
$-(2,5)\otimes 3+(3,5)\otimes 3+(1,2,4,5)\otimes 4-(1,3,4,5)\otimes 4$ & & \\
  
$(1,4)\otimes 1-(2,5)\otimes 1+(3,5)\otimes 1-(1,2,4,5)\otimes 2+(1,3,4,5)\otimes 2-(1,5)\otimes 3$ & 3 & $\u_1$ \\
  
$(1,4)\otimes 1-(2,4)\otimes 1+(3,4)\otimes 1-2(1,2,3,5)\otimes 1-(1,2)\otimes 2+(1,3)\otimes 2$ & 4 & 0 \\
$-2(4,5)\otimes 2-(1,5)\otimes 3-(2,5)\otimes 3+(3,5)\otimes 3+(1,2,4,5)\otimes 4-(1,3,4,5)\otimes 4$ & & \\

\hline
\end{tabular}
\end{table}

For $p\in \cc_{M_7}^\circ$ there are six nilpotent
$Z_{G_0}(p)^\circ$-orbits in $\z_{\g}(p)_1$. Up to the action of the component group
four orbits remain. Table \ref{tab:cm7} lists their representatives.

\begin{table}[htb]\caption{Nilpotent parts of mixed elements with semisimple
    part in $\cc_{M_7}^\circ$.}\label{tab:cm7}
\begin{tabular}{|l|r|r|}
  \hline
  element $e$ & $\dim$ & $\z_{\g_0}(p+e)$ \\
  \hline
$(3,5)\otimes 1 +(1,3)\otimes 4$ & 3 & $\ttt_1+\u_2$\\ 

$-()\otimes 1+(1,3,4,5)\otimes 1+(1,3)\otimes 3-(1,5)\otimes 4$ & 4 & $\ttt_1+\u_1$\\

$-()\otimes 1-(2,3)\otimes 1+(1,3,4,5)\otimes 1-(3,5)\otimes 2+(1,3)\otimes 3$ & 5 & $\u_1$ \\
$-(1,5)\otimes 4-(3,4)\otimes 4+(1,2,3,5)\otimes 4$ & & \\

$(1,4)\otimes 1-(2,3)\otimes 1-(3,5)\otimes 2-(1,5)\otimes 3-(3,4)\otimes 4+(1,2,3,5)\otimes 4$ & 6 & 0 \\
  \hline
  \end{tabular}
\end{table}

For $p\in \cc_{M_8}^\circ$ there are sixtyfour nilpotent
$Z_{G_0}(p)^\circ$-orbits in $\z_{\g}(p)_1$. The action of the component group reduces
this number to thirtyfive. Table \ref{tab:cm8} lists their representatives.

Also in the column ``char.'' we provide the \emph{characteristic} of the nilpotent part of the corresponding representative,
and the column ``Dynkin scheme'' shows the Dynkin scheme of weights for weight vectors having nonzero coefficients in the representative.
These are defined similarly to the table for nilpotent orbits below, except that characteristic is taken with respect to the
centralizer of $p_1$. Accordingly, the characteristic is an element $h$ of the Cartan subalgebra of the centralizer $\z_{\g_0}(p)$ of $p=p_1$,
i.~e. of the reductive Lie algebra of type $A_2+2A_1+T_1$, and is represented by a quadruple of nonnegative integers and a rational number.
The quadruple gives values of simple roots of $A_2+2A_1$ on $h$, while the rational number is the coordinate of $h$ on the 1-dimensional center
of $\z_{\g_0}(p_1)$.

{
\renewcommand{\arraystretch}{1.2}
\small\tabcolsep=3pt
\newcommand*\scaleMe{0.44}

\begin{tiny}

\begin{longtable}[l]{|P{6cm}|P{1.2cm}|P{.4cm}|P{1.5cm}|P{4.2cm}|}

\caption{Nilpotent parts of mixed elements with semisimple
    part in $\cc_{M_8}^\circ$.}\label{tab:cm8}\\
\hline
element $e$ & char. & $\dim$ & $\z_{\g_0}(p_1+e)$ & Dynkin scheme \\
\hline
$\text{(1,4)$\otimes $1}$&$(0110,\frac{1}{3})$&4&$2A_1 + \ttt_2 + \u_3$&\includegraphics[scale=\scaleMe]{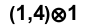}\\\hline
$\text{(1,4)$\otimes $1}\alb-\text{(4,5)$\otimes $4}$&$(0200,\frac{2}{3})$&5&$2A_1 + \ttt_2 + \u_2$&\includegraphics[scale=\scaleMe]{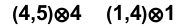}\\\hline
$\text{(1,2)$\otimes $1}\alb+\text{(1,4)$\otimes $4}\alb-\text{(4,5)$\otimes $1}$&$(1110,1)$&7&$A_1 + \ttt_2 + \u_3$&\includegraphics[scale=\scaleMe]{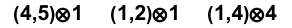}\\\hline
$\text{(1,4)$\otimes $1}\alb+\text{()$\otimes $2}$&$(1111,0)$&7&$\ttt_3 + \u_5$&\includegraphics[scale=\scaleMe]{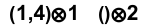}\\\hline
$\text{(1,2)$\otimes $1}\alb-\text{(4,5)$\otimes $4}$&$(2000,\frac{4}{3})$&8&$A_1 + \ttt_2 + \u_2$&\includegraphics[scale=\scaleMe]{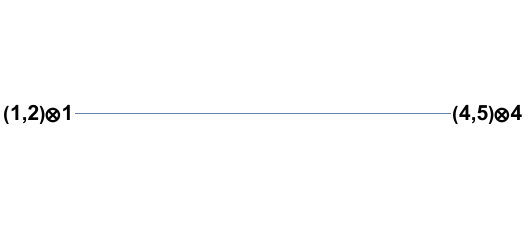}\\\hline
$\text{(1,4)$\otimes $1}\alb+\text{(2,3)$\otimes $2}$&$(0022,0)$&8&$A_1 + \ttt_2 + \u_2$&\includegraphics[scale=\scaleMe]{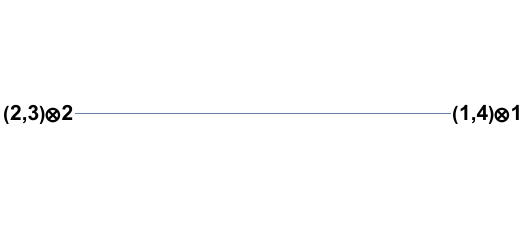}\\\hline
$\text{(1,4)$\otimes $1}\alb-\text{(4,5)$\otimes $4}\alb+\text{()$\otimes $2}$&$(1201,\frac{1}{3})$&8&$\ttt_3 + \u_4$&\includegraphics[scale=\scaleMe]{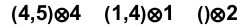}\\\hline
$\text{(1,2)$\otimes $1}\alb+\text{(1,4)$\otimes $4}\alb+\text{(2,3,4,5)$\otimes $4}\alb-\text{(4,5)$\otimes $1}$&$(0000,2)$&9&$2A_1$&\includegraphics[scale=\scaleMe]{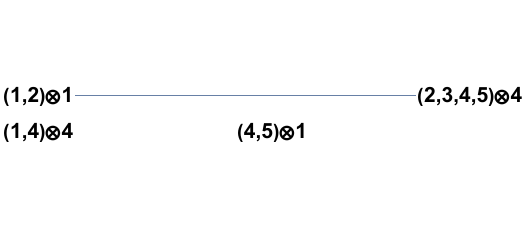}\\\hline
$\text{(1,4)$\otimes $1}\alb+\text{(1,5)$\otimes $2}\alb-\text{(4,5)$\otimes $4}\alb+\text{()$\otimes $3}$&$(2200,0)$&9&$\ttt_3 + \u_3$&\includegraphics[scale=\scaleMe]{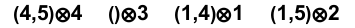}\\\hline
$\text{(1,4)$\otimes $1}\alb+\text{(2,3)$\otimes $2}\alb-\text{(4,5)$\otimes $4}$&$(0204,\frac{2}{3})$&9&$A_1 + \ttt_2 + \u_1$&\includegraphics[scale=\scaleMe]{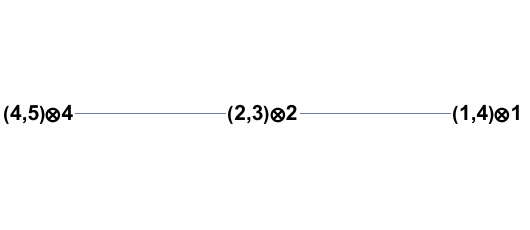}\\\hline
$\text{(1,2)$\otimes $1}\alb+\text{(1,4)$\otimes $4}\alb-\text{(4,5)$\otimes $1}\alb+\text{()$\otimes $2}$&$(2111,\frac{2}{3})$&9&$\ttt_2 + \u_4$&\includegraphics[scale=\scaleMe]{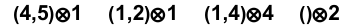}\\\hline
$\text{(1,2)$\otimes $1}\alb-\text{(4,5)$\otimes $4}\alb+\text{()$\otimes $2}$&$(3001,1)$&10&$\ttt_2 + \u_3$&\includegraphics[scale=\scaleMe]{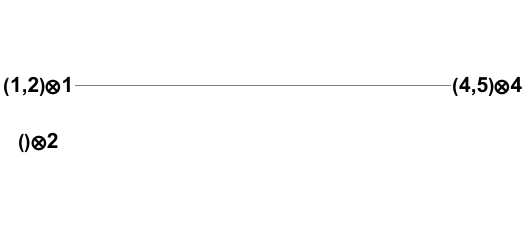}\\\hline
$\text{(1,2,3,4)$\otimes $1}\alb+\text{(1,4)$\otimes $4}\alb-\text{(4,5)$\otimes $1}\alb+\text{()$\otimes $2}$&$(0222,\frac{2}{3})$&10&$\ttt_2 + \u_3$&\includegraphics[scale=\scaleMe]{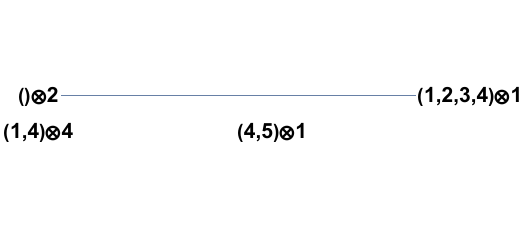}\\\hline
$\text{(1,2)$\otimes $1}\alb+\text{(1,4)$\otimes $4}\alb+\text{(1,5)$\otimes $2}\alb-\text{(4,5)$\otimes $1}\alb+\text{()$\otimes $3}$&$(3110,\frac{1}{3})$&10&$\ttt_2 + \u_3$&\includegraphics[scale=\scaleMe]{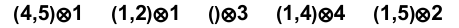}\\\hline
$\text{(1,2)$\otimes $1}\alb+\text{(1,5)$\otimes $2}\alb-\text{(4,5)$\otimes $4}\alb+\text{()$\otimes $3}$&$(4000,\frac{2}{3})$&11&$\ttt_2 + \u_2$&\includegraphics[scale=\scaleMe]{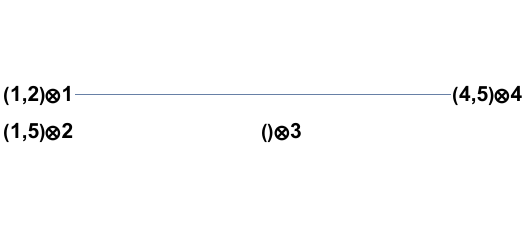}\\\hline
$\text{(1,2,3,4)$\otimes $1}\alb+\text{(1,4)$\otimes $4}\alb+\text{(1,5)$\otimes $2}\alb-\text{(4,5)$\otimes $1}\alb+\text{()$\otimes $3}$&$(2240,0)$&11&$\ttt_2 + \u_2$&\includegraphics[scale=\scaleMe]{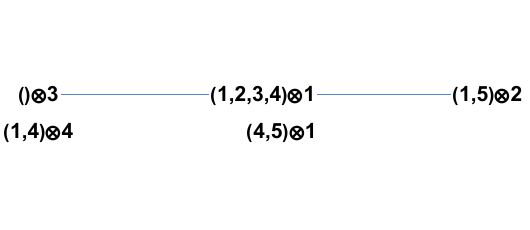}\\\hline
$\text{(1,2)$\otimes $1}\alb+\text{(1,4)$\otimes $4}\alb+\text{(1,5)$\otimes $2}\alb+\text{(3,4)$\otimes $2}\alb-\text{(4,5)$\otimes $1}\alb+\text{()$\otimes $3}$&$(2222,0)$&11&$\ttt_1 + \u_3$&\includegraphics[scale=\scaleMe]{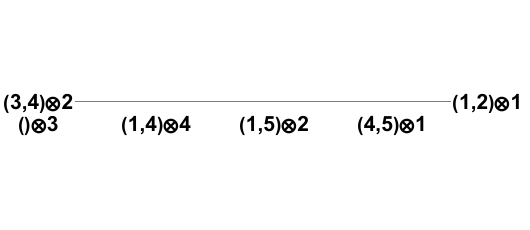}\\\hline
$\text{(1,2)$\otimes $1}\alb+\text{(3,4)$\otimes $2}\alb-\text{(4,5)$\otimes $4}$&$(1213,\frac{4}{3})$&11&$\ttt_2 + \u_2$&\includegraphics[scale=\scaleMe]{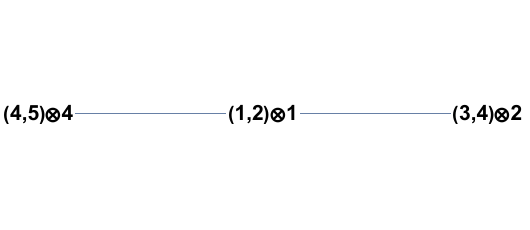}\\\hline
$\text{(1,2)$\otimes $1}\alb+\text{(1,4)$\otimes $4}\alb+\text{(2,3)$\otimes $2}\alb-\text{(4,5)$\otimes $1}$&$(1114,1)$&11&$\ttt_2 + \u_2$&\includegraphics[scale=\scaleMe]{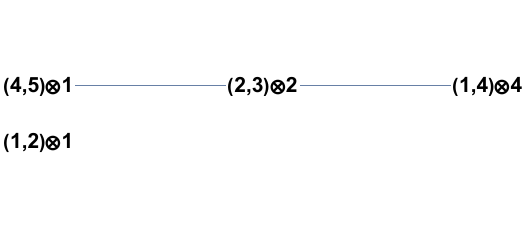}\\\hline
$\text{(1,2)$\otimes $1}\alb+\text{(2,3)$\otimes $2}\alb+\text{(3,4)$\otimes $2}\alb-\text{(4,5)$\otimes $4}$&$(2004,\frac{4}{3})$&12&$\ttt_1 + \u_2$&\includegraphics[scale=\scaleMe]{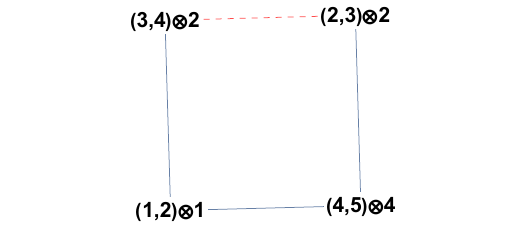}\\\hline
$\text{(1,2,3,4)$\otimes $1}\alb+\text{(1,5)$\otimes $2}\alb-\text{(4,5)$\otimes $4}\alb+\text{()$\otimes $3}$&$(2640,\frac{4}{3})$&12&$\ttt_2 + \u_1$&\includegraphics[scale=\scaleMe]{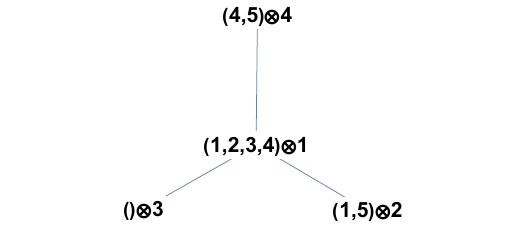}\\\hline
$\text{(1,2)$\otimes $1}\alb+\text{(1,4)$\otimes $4}\alb+\text{(1,5)$\otimes $2}\alb+\text{(2,3)$\otimes $2}\alb-\text{(4,5)$\otimes $1}\alb+\text{()$\otimes $3}$&$(3114,\frac{1}{3})$&12&$\ttt_1 + \u_2$&\includegraphics[scale=\scaleMe]{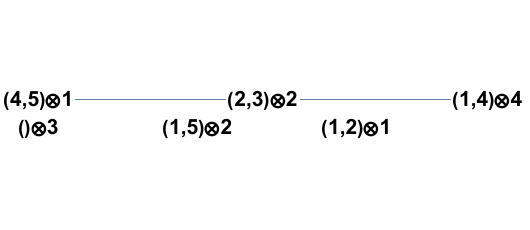}\\\hline
$\text{(1,2,3,4)$\otimes $1}\alb+\text{(1,2)$\otimes $4}\alb-\text{(2,5)$\otimes $1}\alb-\text{(4,5)$\otimes $4}\alb+\text{()$\otimes $2}$&$(1113,2)$&12&$\ttt_1 + \u_2$&\includegraphics[scale=\scaleMe]{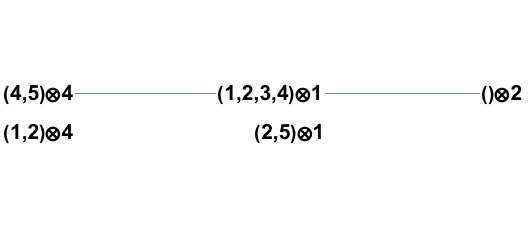}\\\hline
$\text{(1,2)$\otimes $1}\alb+\text{(1,5)$\otimes $2}\alb+\text{(3,4)$\otimes $2}\alb-\text{(4,5)$\otimes $4}\alb+\text{()$\otimes $3}$&$(3213,\frac{2}{3})$&12&$\ttt_1 + \u_2$&\includegraphics[scale=\scaleMe]{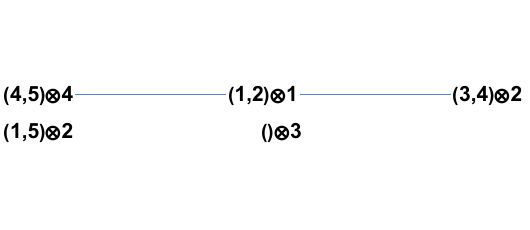}\\\hline
$-\text{(1,2,3,4)$\otimes $4}\alb+\text{(1,2)$\otimes $1}\alb+\text{(2,3,4,5)$\otimes $1}\alb-\text{(4,5)$\otimes $4}\alb+\text{()$\otimes $2}$&$(0004,2)$&13&$\ttt_1 + \u_1$&\includegraphics[scale=\scaleMe]{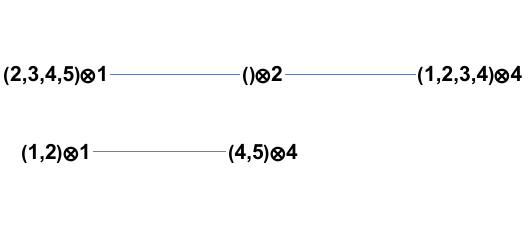}\\\hline
$\text{(1,2)$\otimes $1}\alb+\text{(1,5)$\otimes $2}\alb+\text{(2,3)$\otimes $2}\alb+\text{(3,4)$\otimes $2}\alb-\text{(4,5)$\otimes $4}\alb+\text{()$\otimes $3}$&$(4004,\frac{2}{3})$&13&$\u_2$&\includegraphics[scale=\scaleMe]{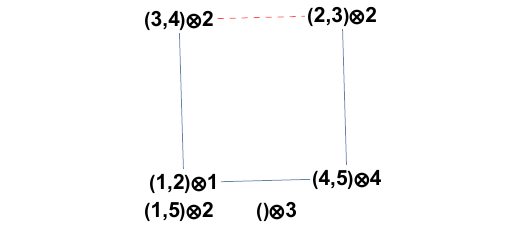}\\\hline
$\text{(1,2,3,4)$\otimes $1}\alb+\text{(1,2)$\otimes $4}\alb+\text{(1,5)$\otimes $2}\alb-\text{(2,5)$\otimes $1}\alb-\text{(4,5)$\otimes $4}\alb+\text{()$\otimes $3}$&$(4440,2)$&13&$\ttt_1 + \u_1$&\includegraphics[scale=\scaleMe]{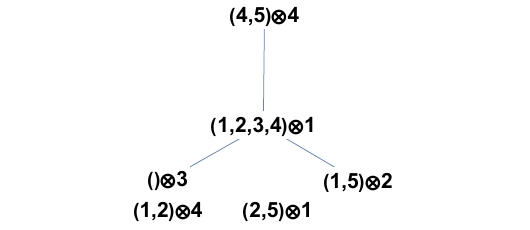}\\\hline
$\text{(1,2,3,4)$\otimes $1}\alb+\text{(1,4)$\otimes $4}\alb+\text{(1,5)$\otimes $2}\alb+\text{(2,3)$\otimes $2}\alb-\text{(4,5)$\otimes $1}\alb+\text{()$\otimes $3}$&$(2244,0)$&13&$\ttt_1 + \u_1$&\includegraphics[scale=\scaleMe]{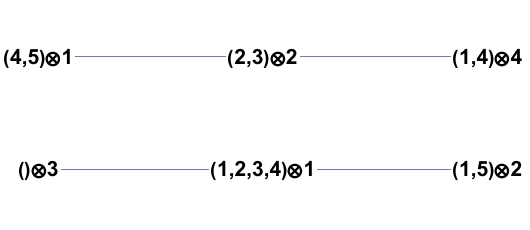}\\\hline
$\text{(1,2,3,4)$\otimes $1}\alb+\text{(1,5)$\otimes $2}\alb+\text{(3,4)$\otimes $2}\alb-\text{(4,5)$\otimes $4}\alb+\text{()$\otimes $3}$&$(1741,1)$&13&$\ttt_1 + \u_1$&\includegraphics[scale=\scaleMe]{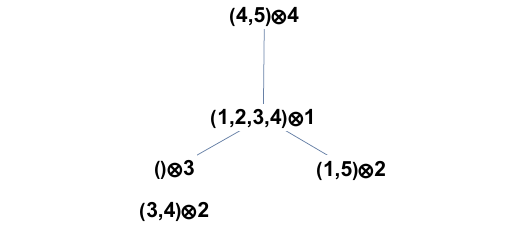}\\\hline
$\text{(1,2)$\otimes $1}\alb-\text{(1,5)$\otimes $3}\alb+\text{(3,4)$\otimes $2}\alb-\text{(4,5)$\otimes $4}$&$(4422,0)$&13&$\ttt_1 + \u_1$&\includegraphics[scale=\scaleMe]{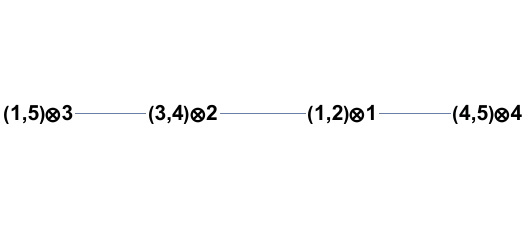}\\\hline
$\text{(1,2)$\otimes $1}\alb+\text{(1,3,4,5)$\otimes $3}\alb-\text{(1,5)$\otimes $3}\alb+\text{(3,4)$\otimes $2}\alb-\text{(4,5)$\otimes $4}$&$(0840,\frac{2}{3})$&14&$\u_1$&\includegraphics[scale=\scaleMe]{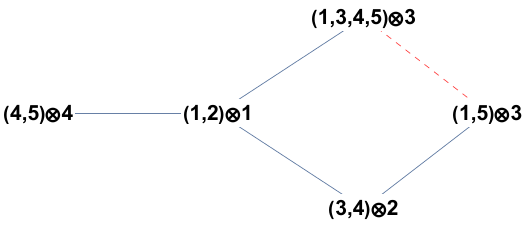}\\\hline
$\text{(1,2,3,4)$\otimes $1}\alb+\text{(1,2)$\otimes $4}\alb+\text{(1,5)$\otimes $2}\alb-\text{(2,5)$\otimes $1}\alb+\text{(3,4)$\otimes $2}\alb-\text{(4,5)$\otimes $4}\alb+\text{()$\otimes $3}$&$(4444,2)$&14&$\u_1$&\includegraphics[scale=\scaleMe]{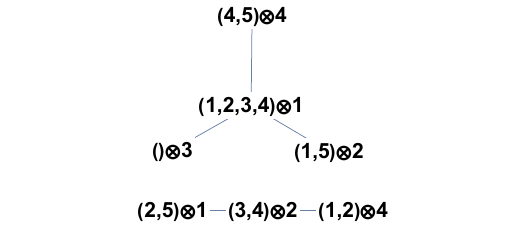}\\\hline
$\text{(1,2,3,4)$\otimes $1}\alb+\text{(1,5)$\otimes $2}\alb+\text{(2,3)$\otimes $2}\alb-\text{(4,5)$\otimes $4}\alb+\text{()$\otimes $3}$&$(2648,\frac{4}{3})$&14&$\ttt_1$&\includegraphics[scale=\scaleMe]{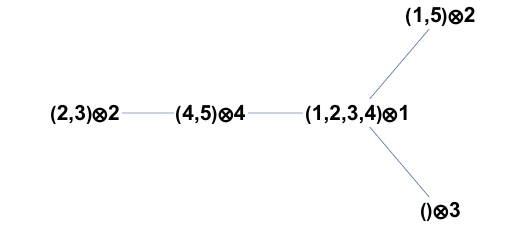}\\\hline
$\text{(1,2,3,4)$\otimes $1}\alb+\text{(1,2)$\otimes $4}\alb-\text{(1,5)$\otimes $3}\alb-\text{(2,5)$\otimes $1}\alb+\text{(3,4)$\otimes $2}\alb-\text{(4,5)$\otimes $4}$&$(8884,2)$&15&$0$&\includegraphics[scale=\scaleMe]{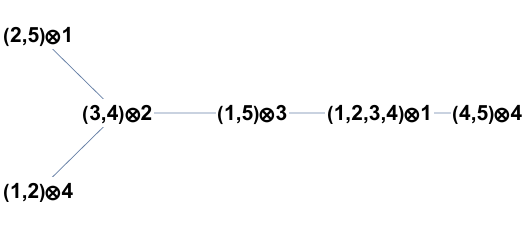}\\\hline
$\text{(1,2,3,4)$\otimes $1}\alb+\text{(1,2)$\otimes $4}\alb+\text{(1,5)$\otimes $2}\alb+\text{(2,3)$\otimes $2}\alb-\text{(2,5)$\otimes $1}\alb-\text{(4,5)$\otimes $4}\alb+\text{()$\otimes $3}$&$(4448,2)$&15&$0$&\includegraphics[scale=\scaleMe]{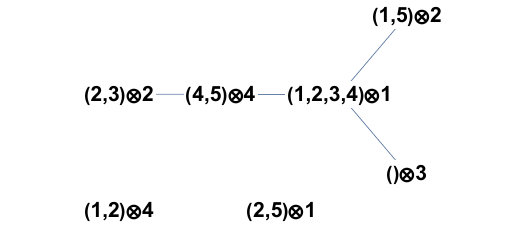}\\
\hline
\end{longtable}

\end{tiny}

}

\section{The nilpotent orbits}\label{sec:nilp}

This section contains the list of nilpotent orbits (Table \ref{tab:nilorb}) and
the Hasse diagram displaying their closure ordering (Figure \ref{fig}).
The representatives of the nilpotent orbits have been computed
with the algorithms of \cite{gra15}. The closure ordering has been computed
with the algorithm of \cite{gravinya}. For accounts of both algorithms we also
refer to \cite[Chapter 8]{gra16}. 

For the notation used to indicate the centralizer $\z_{\g_0}(e)$ we refer to
the previous section. In the fourth column we give the {\em characteristic}
of the nilpotent $e$ in the second column. This is defined as follows. The
nilpotent $e\in \g_1$ lies in a homogeneous $\ssl_2$-triple $(h,e,f)$ with
$h\in \g_0$, $f\in \g_{-1}$ and
$$[h,e]=2e,\, [h,f]=-2f,\, [e,f]=h.$$
Furthermore, it can be shown that $e,e'$ lying in homogeneous $\ssl_2$-triples
$(h,e,f)$, $(h',e',f')$ are $G_0$-conjugate if and only if the triples
are $G_0$-conjugate if and only if $h,h'$ are $G_0$-conjugate (cf.
\cite[Theorem 8.3.6]{gra16}). The element $h\in \g_0$ lies in a Cartan
subalgebra of $\g_0$. We consider the corresponding root system of $\g_0$
and its set of simple roots $\gamma_1,\ldots,\gamma_8$. We have that $h$ is
conjugate under the action of the Weyl group to a unique $\tilde h$ with
$\gamma_i(\tilde h)\geq 0$ for all $i$. In the last column of the table we
list the numbers $\gamma_i(\tilde h)$; they uniquely determine the orbit.
For this we use the following enumeration of the Dynkin diagram of $\g_0$, shown embedded into the extended Dynkin diagram $E_8$ for $\g$:
\begin{center}
\begin{tikzpicture}
\node[dnode,label=above:{\small $1$}] (1) at (0,1) {};
\node[dnode,label=above:{\small $2$}] (2) at (1,1) {};
\node[dnode,label=above:{\small $3$}] (3) at (2,1) {};
\node[dnode,label=below:{\small $4$}] (4) at (2,0) {};
\node[dnode,label=above:{\small $5$}] (5) at (3,1) {};
\node at (4,1) {$\oplus$};
\node[dnode,label=above:{\small $6$}] (6) at (5,1) {};
\node[dnode,label=above:{\small $7$}] (7) at (6,1) {};
\node[dnode,label=above:{\small $8$}] (8) at (7,1) {};
\path (1) edge[sedge] (2)
(2) edge[sedge] (3)
(3) edge[sedge] (4)
(3) edge[sedge] (5)
(6) edge[sedge] (7)
(7) edge[sedge] (8);
\end{tikzpicture}
\end{center}
In the last column the Dynkin scheme for the representative is given, as explained at the end of Section \ref{sec:vinb}.
In some cases the representatives have been chosen in such a way that their Dynkin schemes would have certain standard form.
In these cases, the criterion to decide that the chosen element $e$ indeed lies on the required orbit was as follows.
For the corresponding characteristic $h$ the $\z_{\g_0}(h)$-module $M_h:=\{x\in\g_1\mid[h,x]=2x\}$ has been computed,
and it was checked that the vector space $\{[a,e]\mid a\in\z_{\g_0}(h)\}$ coincides with $M_h$. 
The described criterion suffices as it ensures that $e$ lies on an open $Z_{G_0}(h)$-orbit in $M_h$,
as can be concluded, for example, from the more general fact given in \cite[Lemma 1]{elashvili}.
In fact, the direct sum decomposition of this module has been used to choose $e$ in the required form.
Under the required form here are meant Dynkin schemes with shapes that allow to identify the corresponding
nilpotent orbits in $\g$. In particular, inspecting the table one can conclude that every nilpotent $E_8$-orbit
occurs at least once.

{
\newcommand*\imgscale{0.304}
\setlength{\arraycolsep}{1pt}
\renewcommand{\arraystretch}{1.2}
\small\tabcolsep=3pt

\begin{tiny}
\begin{longtable}[l]{|P{.4cm}|P{4cm}|P{2cm}|P{.4cm}|P{1.8cm}|P{4.5cm}|}
\caption{Nilpotent orbits}
\label{tab:nilorb}\\\hline

N &   element $e$ &  $\z_{\g_0}(e)$ & $\dim$  & characteristic & Dynkin scheme \\
  \hline

1&$\text{(1,2,4,5)$\otimes$4}\alb+\text{(1,2)$\otimes$1}\alb+\text{(1,3,4,5)$\otimes$1}\alb+\text{(1,3)$\otimes$2}\alb+\text{(1,5)$\otimes$3}\alb+\text{(2,3,4,5)$\otimes$3}\alb+\text{(2,4)$\otimes$2}\alb+\text{()$\otimes$4}$& 0&60&$\begin{array}{cccccccc}
 8 & 8 & 8 & 8 & \oplus  & 16 & 8 & 8 \\
  &  & 8\\
\end{array}$&\includegraphics[scale=\imgscale]{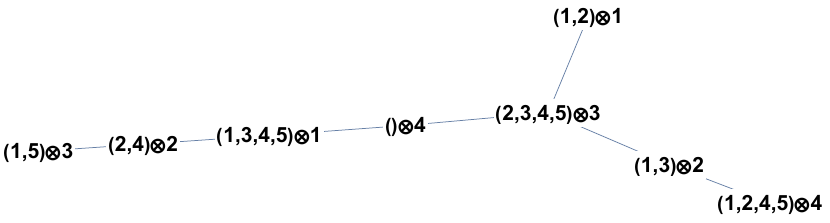}\\\hline
2&$\text{(1,2,3,4)$\otimes$2}\alb+\text{(1,2,4,5)$\otimes$4}\alb+\text{(1,3,4,5)$\otimes$1}\alb+\text{(1,5)$\otimes$3}\alb+\text{(2,3,4,5)$\otimes$3}\alb+\text{(2,3)$\otimes$4}\alb+\text{(4,5)$\otimes$4}\alb+\text{()$\otimes$1}$&  0 &60&$\begin{array}{cccccccc}
 8 & 8 & 0 & 8 & \oplus  & 8 & 8 & 16 \\
  &  & 8\\
\end{array}$&\includegraphics[scale=\imgscale]{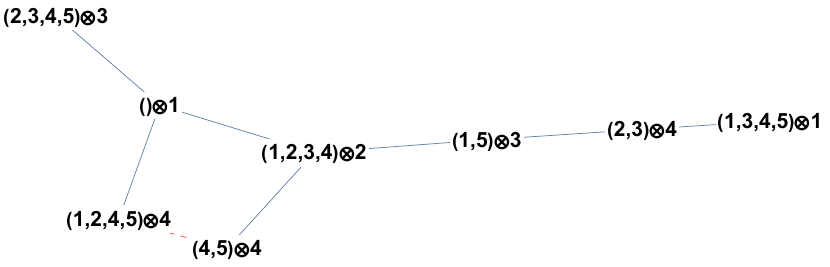}\\\hline
3&$\text{(1,2,3,5)$\otimes$3}\alb+\text{(1,2,4,5)$\otimes$4}\alb+\text{(1,2)$\otimes$1}\alb+\text{(1,3,4,5)$\otimes$1}\alb+\text{(1,3)$\otimes$2}\alb+\text{(2,3)$\otimes$4}\alb+\text{(2,4)$\otimes$1}\alb+\text{(4,5)$\otimes$3}$& $\u_1$ &59&$\begin{array}{cccccccc}
 8 & 0 & 8 & 0 & \oplus  & 8 & 8 & 8 \\
  &  & 8\\
\end{array}$&\includegraphics[scale=\imgscale]{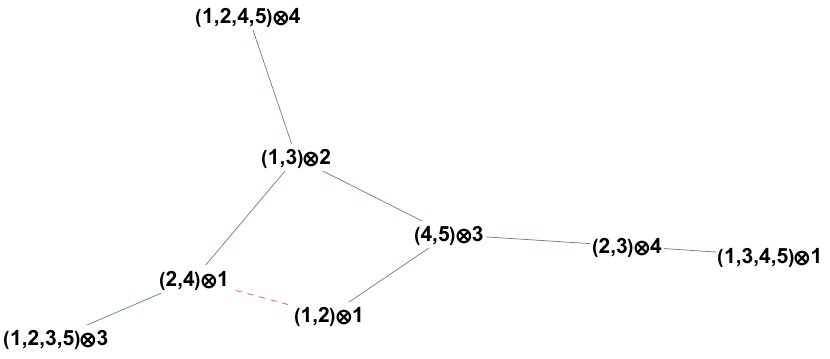}\\\hline
4&$\text{(1,2,4,5)$\otimes$4}\alb+\text{(1,2)$\otimes$1}\alb+\text{(1,3,4,5)$\otimes$1}\alb+\text{(1,3)$\otimes$2}\alb+\text{(1,5)$\otimes$3}\alb+\text{(2,3)$\otimes$3}\alb+\text{(2,4)$\otimes$2}\alb+\text{()$\otimes$4}$& $\u_1$ &59&$\begin{array}{cccccccc}
 4 & 4 & 4 & 8 & \oplus  & 8 & 4 & 4 \\
  &  & 4\\
\end{array}$&\includegraphics[scale=\imgscale]{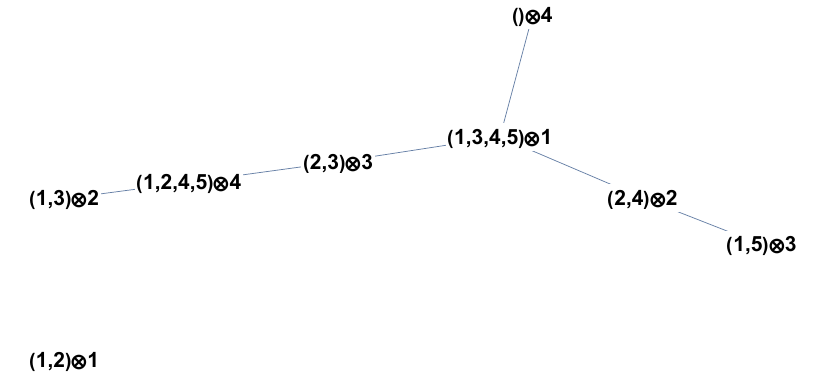}\\\hline
5&$\text{(1,2,3,5)$\otimes$3}\alb+\text{(1,2)$\otimes$2}\alb+\text{(1,3,4,5)$\otimes$4}\alb+\text{(1,3)$\otimes$1}\alb+\text{(2,4)$\otimes$1}\alb+\text{(3,4)$\otimes$2}\alb+\text{(4,5)$\otimes$3}\alb+\text{()$\otimes$4}$& $\u_2$&58&$\begin{array}{cccccccc}
 8 & 0 & 8 & 0 & \oplus  & 8 & 8 & 0 \\
  &  & 0\\
\end{array}$&\includegraphics[scale=\imgscale]{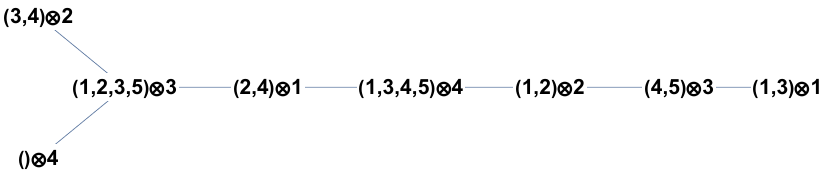}\\\hline
6&$\text{(1,2,4,5)$\otimes$4}\alb+\text{(1,2)$\otimes$1}\alb+\text{(1,3,4,5)$\otimes$1}\alb+\text{(1,3)$\otimes$2}\alb+\text{(2,3)$\otimes$3}\alb+\text{(2,4)$\otimes$2}\alb+\text{(4,5)$\otimes$3}\alb+\text{()$\otimes$4}$& $\u_2$&58&$\begin{array}{cccccccc}
 4 & 4 & 4 & 4 & \oplus  & 4 & 4 & 4 \\
  &  & 4\\
\end{array}$&\includegraphics[scale=\imgscale]{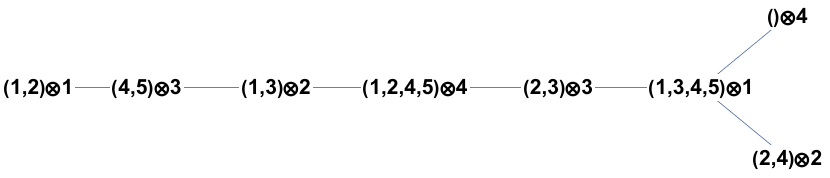}\\\hline
7&$\text{(2,4)$\otimes$1}\alb+\text{(1,3,4,5)$\otimes$1}\alb+\text{(1,4)$\otimes$2}\alb+\text{(1,2,3,4)$\otimes$2}\alb+\text{(1,5)$\otimes$3}\alb+\text{(2,3,4,5)$\otimes$3}\alb+\text{(2,3)$\otimes$4}\alb+\text{(4,5)$\otimes$4}$& $\u_2$ &58&$\begin{array}{cccccccc}
 0 & 0 & 8 & 0 & \oplus  & 0 & 8 & 8 \\
  &  & 0\\
\end{array}$&\includegraphics[scale=\imgscale]{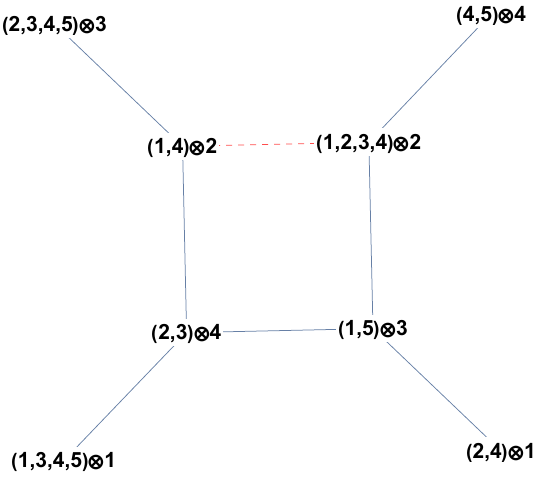}\\\hline
8&$\text{(1,2,4,5)$\otimes$4}\alb+\text{(1,2)$\otimes$1}\alb+\text{(1,3,4,5)$\otimes$4}\alb+\text{(1,5)$\otimes$3}\alb+\text{(2,3,4,5)$\otimes$3}\alb+\text{(3,4)$\otimes$2}\alb+\text{(4,5)$\otimes$1}\alb+\text{()$\otimes$4}$& $\u_2$ &58&$\begin{array}{cccccccc}
 0 & 8 & 0 & 8 & \oplus  & 8 & 0 & 8 \\
  &  & 0\\
\end{array}$&\includegraphics[scale=\imgscale]{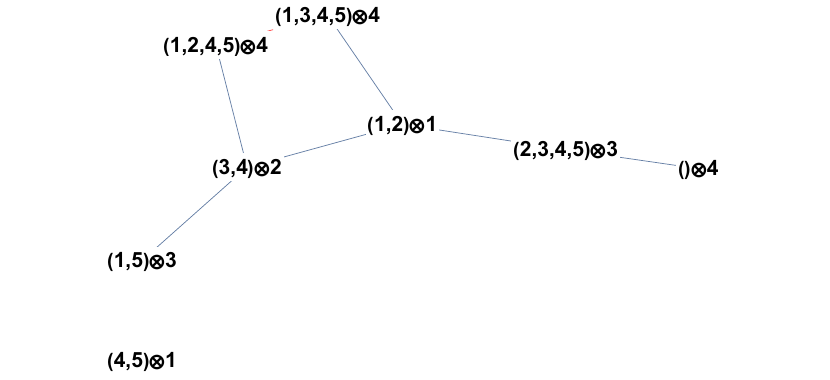}\\\hline
9&$\text{(1,2,4,5)$\otimes$4}\alb+\text{(1,3,4,5)$\otimes$1}\alb+\text{(1,3)$\otimes$2}\alb+\text{(1,5)$\otimes$3}\alb+\text{(2,3)$\otimes$3}\alb+\text{(2,4)$\otimes$2}\alb+\text{()$\otimes$4}$& $\ttt_1\alb+\u_1$&58&$\begin{array}{cccccccc}
 3 & 5 & 3 & 8 & \oplus  & 8 & 5 & 3 \\
  &  & 5\\
\end{array}$&\includegraphics[scale=\imgscale]{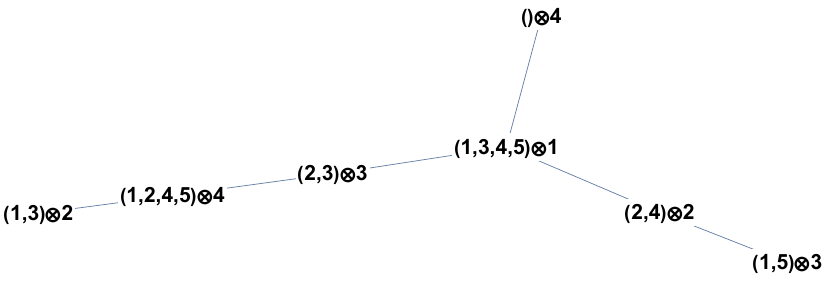}\\\hline
10&$\text{(1,2,4,5)$\otimes$4}\alb+\text{(1,2)$\otimes$1}\alb+\text{(1,3,4,5)$\otimes$1}\alb+\text{(1,5)$\otimes$3}\alb+\text{(2,3)$\otimes$3}\alb+\text{(3,4)$\otimes$2}\alb+\text{(4,5)$\otimes$3}\alb+\text{()$\otimes$4}$& $\u_3$ &57&$\begin{array}{cccccccc}
 4 & 0 & 4 & 4 & \oplus  & 4 & 4 & 4 \\
  &  & 4\\
\end{array}$&\includegraphics[scale=\imgscale]{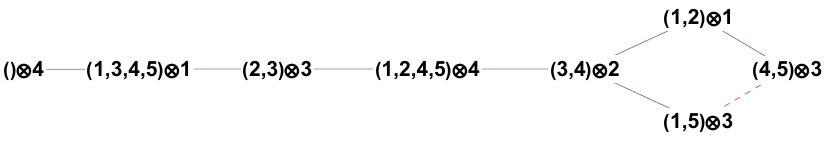}\\\hline
11&$\text{(1,2,3,4)$\otimes$2}\alb+\text{(1,2,3,5)$\otimes$3}\alb+\text{(1,2,4,5)$\otimes$4}\alb+\text{(1,2)$\otimes$1}\alb+\text{(1,3,4,5)$\otimes$1}\alb+\text{(1,4)$\otimes$2}\alb+\text{(2,3)$\otimes$4}\alb+\text{(4,5)$\otimes$3}$& $\u_3$&57&$\begin{array}{cccccccc}
 4 & 0 & 4 & 0 & \oplus  & 4 & 4 & 8 \\
  &  & 4\\
\end{array}$&\includegraphics[scale=\imgscale]{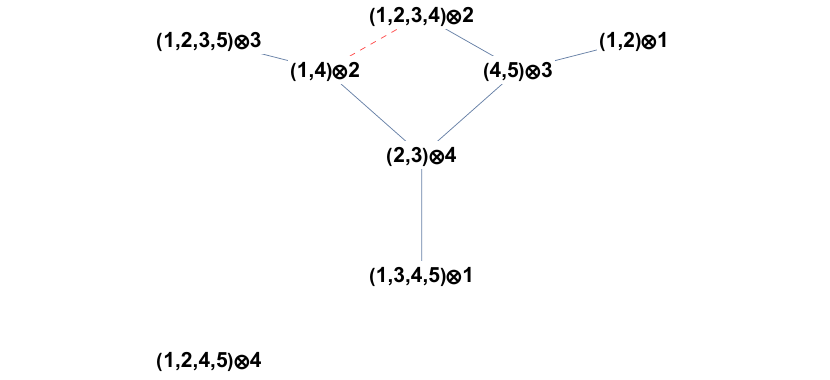}\\\hline
12&$\text{(1,2,3,4)$\otimes$2}\alb+\text{(1,3)$\otimes$1}\alb+\text{(1,5)$\otimes$3}\alb+\text{(2,3,4,5)$\otimes$3}\alb+\text{(2,3)$\otimes$4}\alb+\text{(2,4)$\otimes$1}\alb+\text{(4,5)$\otimes$4}$& $\ttt_1\alb+\u_2$&57&$\begin{array}{cccccccc}
 1 & 1 & 6 & 1 & \oplus  & 1 & 8 & 7 \\
  &  & 1\\
\end{array}$&\includegraphics[scale=\imgscale]{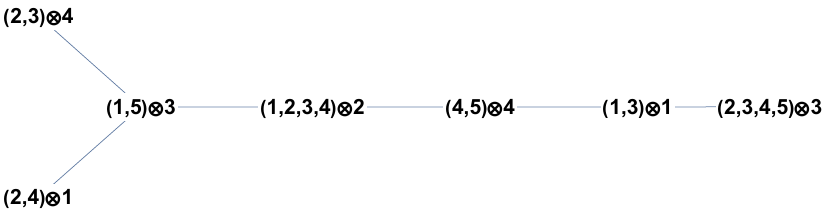}\\\hline
13&$\text{(1,2,4,5)$\otimes$4}\alb+\text{(1,2)$\otimes$1}\alb+\text{(1,3,4,5)$\otimes$4}\alb+\text{(1,5)$\otimes$3}\alb+\text{(2,3,4,5)$\otimes$3}\alb+\text{(3,4)$\otimes$2}\alb+\text{()$\otimes$4}$& $\ttt_1\alb+\u_2$ &57&$\begin{array}{cccccccc}
 1 & 7 & 1 & 7 & \oplus  & 8 & 1 & 7 \\
  &  & 0\\
\end{array}$&\includegraphics[scale=\imgscale]{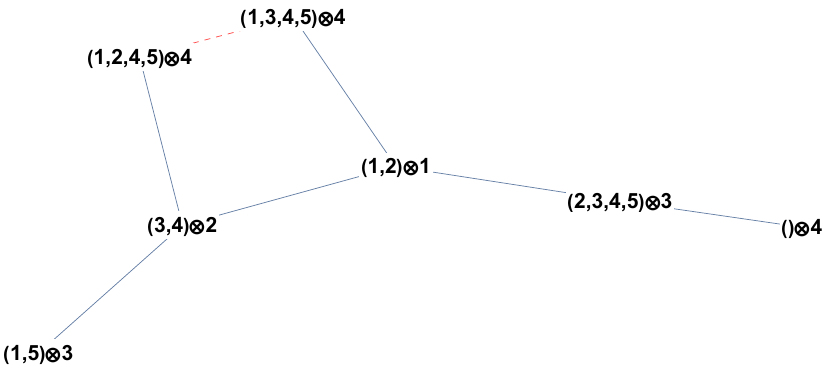}\\\hline
14&$\text{(1,2,3,5)$\otimes$3}\alb+\text{(1,2,4,5)$\otimes$2}\alb+\text{(1,2)$\otimes$4}\alb+\text{(1,3,4,5)$\otimes$4}\alb+\text{(1,3)$\otimes$1}\alb+\text{(2,4)$\otimes$1}\alb+\text{(3,4)$\otimes$2}\alb+\text{()$\otimes$3}$& $\u_4$ &56&$\begin{array}{cccccccc}
 0 & 0 & 8 & 0 & \oplus  & 8 & 0 & 0 \\
  &  & 0\\
\end{array}$&\includegraphics[scale=\imgscale]{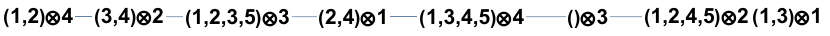}\\\hline
15&$\text{(1,2,4,5)$\otimes$4}\alb+\text{(1,3)$\otimes$1}\alb+\text{(1,5)$\otimes$3}\alb+\text{(2,3,4,5)$\otimes$3}\alb+\text{(2,3)$\otimes$3}\alb+\text{(2,4)$\otimes$1}\alb+\text{(3,4)$\otimes$2}\alb+\text{()$\otimes$4}$& $\u_4$ &56&$\begin{array}{cccccccc}
 0 & 4 & 0 & 4 & \oplus  & 4 & 4 & 4 \\
  &  & 4\\
\end{array}$&\includegraphics[scale=\imgscale]{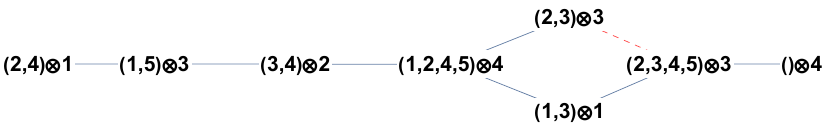}\\\hline
16&$\text{(1,2)$\otimes$2}\alb+\text{(1,3,4,5)$\otimes$4}\alb+\text{(1,3)$\otimes$1}\alb+\text{(2,3)$\otimes$3}\alb+\text{(2,4)$\otimes$1}\alb+\text{(3,4)$\otimes$2}\alb+\text{(4,5)$\otimes$3}\alb+\text{()$\otimes$4}$& $\u_4$ &56&$\begin{array}{cccccccc}
 4 & 0 & 4 & 4 & \oplus  & 4 & 4 & 0 \\
  &  & 0\\
\end{array}$&\includegraphics[scale=\imgscale]{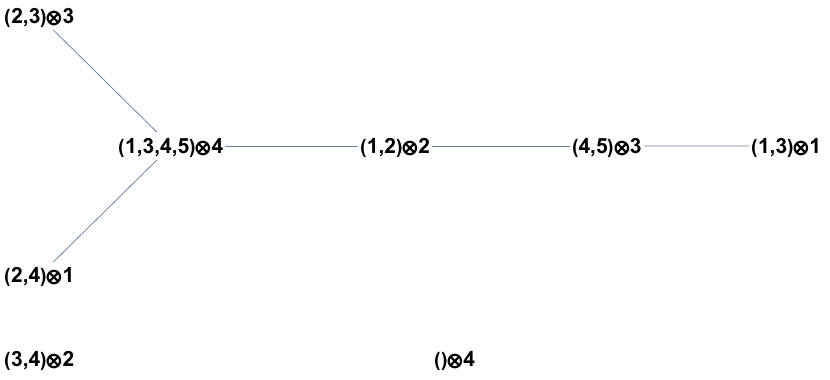}\\\hline
17&$\text{(1,2,4,5)$\otimes$4}\alb+\text{(1,2)$\otimes$1}\alb+\text{(1,3,4,5)$\otimes$1}\alb+\text{(1,4)$\otimes$2}\alb+\text{(2,3)$\otimes$3}\alb+\text{(2,4)$\otimes$2}\alb+\text{(4,5)$\otimes$3}\alb+\text{()$\otimes$4}$& $\u_4$ &56&$\begin{array}{cccccccc}
 0 & 8 & 0 & 0 & \oplus  & 0 & 0 & 8 \\
  &  & 0\\
\end{array}$&\includegraphics[scale=\imgscale]{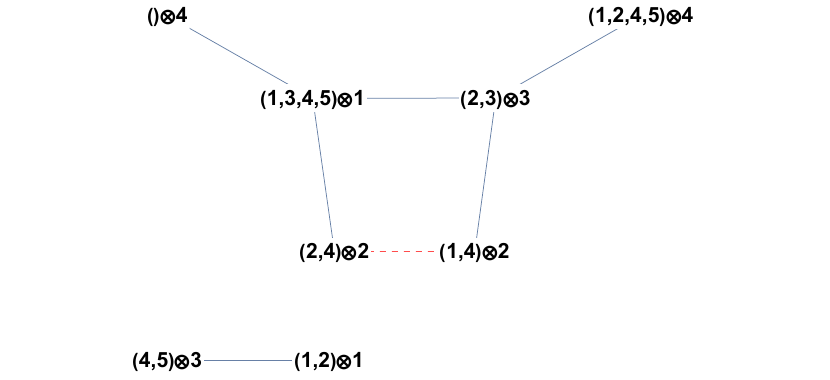}\\\hline
18&$\text{(1,2,4,5)$\otimes$4}\alb+\text{(1,3,4,5)$\otimes$1}\alb+\text{(1,3)$\otimes$2}\alb+\text{(2,3)$\otimes$3}\alb+\text{(2,4)$\otimes$1}\alb+\text{(4,5)$\otimes$3}\alb+\text{()$\otimes$4}$& $\ttt_1\alb+\u_3$ &56&$\begin{array}{cccccccc}
 3 & 2 & 3 & 3 & \oplus  & 3 & 5 & 3 \\
  &  & 5\\
\end{array}$&\includegraphics[scale=\imgscale]{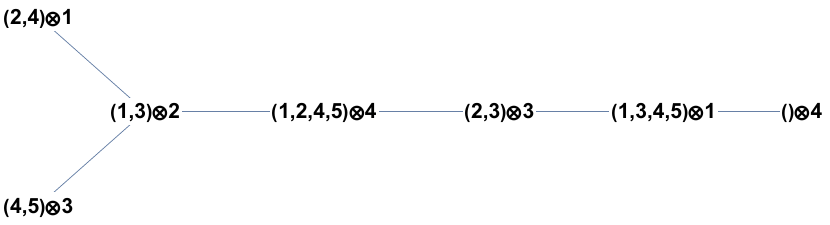}\\\hline
19&$\text{(1,2,3,4)$\otimes$2}\alb+\text{(1,2,3,5)$\otimes$3}\alb+\text{(1,2)$\otimes$1}\alb+\text{(1,3,4,5)$\otimes$1}\alb+\text{(1,4)$\otimes$2}\alb+\text{(2,3)$\otimes$4}\alb+\text{(4,5)$\otimes$3}$& $\ttt_1\alb+\u_3$&56&$\begin{array}{cccccccc}
 5 & 0 & 3 & 0 & \oplus  & 5 & 3 & 8 \\
  &  & 5\\
\end{array}$&\includegraphics[scale=\imgscale]{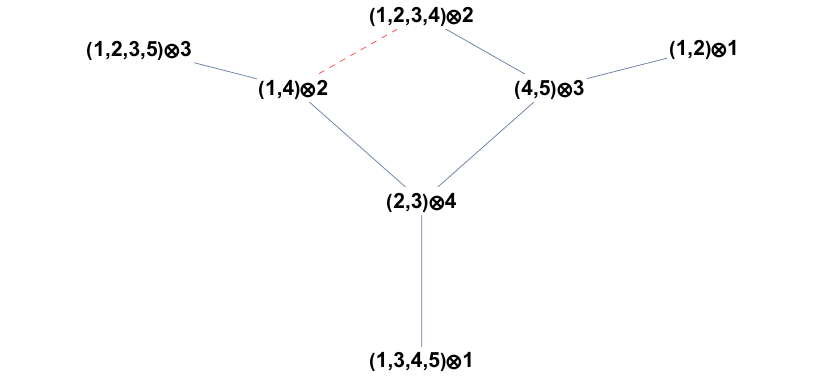}\\\hline
20&$\text{(1,2,4,5)$\otimes$4}\alb+\text{(1,3,4,5)$\otimes$1}\alb+\text{(1,4)$\otimes$2}\alb+\text{(1,5)$\otimes$3}\alb+\text{(2,3,4,5)$\otimes$3}\alb+\text{(2,3)$\otimes$4}\alb+\text{(2,4)$\otimes$1}$& $\ttt_1\alb+\u_3$ &56&$\begin{array}{cccccccc}
 3 & 2 & 1 & 2 & \oplus  & 3 & 5 & 8 \\
  &  & 5\\
\end{array}$&\includegraphics[scale=\imgscale]{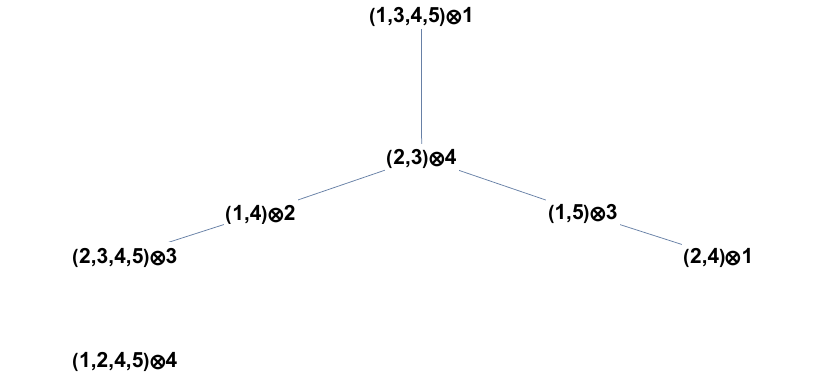}\\\hline
21&$\text{(1,2,3,4)$\otimes$2}\alb+\text{(1,3,4,5)$\otimes$4}\alb+\text{(1,3)$\otimes$1}\alb+\text{(1,5)$\otimes$3}\alb+\text{(2,3,4,5)$\otimes$3}\alb+\text{(2,3)$\otimes$4}\alb+\text{(2,4)$\otimes$1}\alb+\text{()$\otimes$4}$& $\u_5$ &55&$\begin{array}{cccccccc}
 4 & 0 & 4 & 0 & \oplus  & 4 & 4 & 4 \\
  &  & 0\\
\end{array}$&\includegraphics[scale=\imgscale]{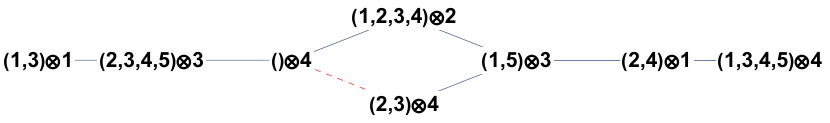}\\\hline
22&$\text{(1,2,3,4)$\otimes$2}\alb+\text{(1,2,4,5)$\otimes$4}\alb+\text{(1,2)$\otimes$1}\alb+\text{(1,3,4,5)$\otimes$1}\alb+\text{(2,3)$\otimes$3}\alb+\text{(4,5)$\otimes$3}\alb+\text{()$\otimes$4}$& $\ttt_1\alb+\u_4$ &55&$\begin{array}{cccccccc}
 1 & 5 & 1 & 1 & \oplus  & 1 & 1 & 7 \\
  &  & 1\\
\end{array}$&\includegraphics[scale=\imgscale]{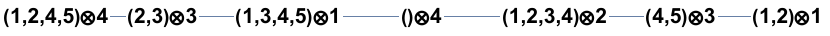}\\\hline
23&$\text{(1,2,4,5)$\otimes$4}\alb+\text{(1,3)$\otimes$1}\alb+\text{(1,5)$\otimes$3}\alb+\text{(2,3,4,5)$\otimes$3}\alb+\text{(2,3)$\otimes$3}\alb+\text{(2,4)$\otimes$1}\alb+\text{(3,4)$\otimes$2}$& $\ttt_1\alb+\u_4$ &55&$\begin{array}{cccccccc}
 0 & 2 & 4 & 2 & \oplus  & 6 & 2 & 2 \\
  &  & 2\\
\end{array}$&\includegraphics[scale=\imgscale]{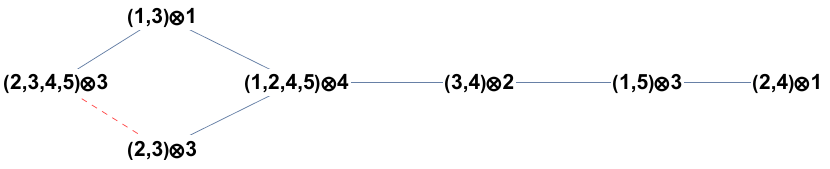}\\\hline
24&$\text{(1,2)$\otimes$2}\alb+\text{(1,3,4,5)$\otimes$4}\alb+\text{(1,3)$\otimes$1}\alb+\text{(1,5)$\otimes$3}\alb+\text{(2,3)$\otimes$3}\alb+\text{(2,4)$\otimes$1}\alb+\text{(3,4)$\otimes$2}$& $\ttt_1\alb+\u_4$ &55&$\begin{array}{cccccccc}
 3 & 0 & 5 & 3 & \oplus  & 5 & 3 & 0 \\
  &  & 0\\
\end{array}$&\includegraphics[scale=\imgscale]{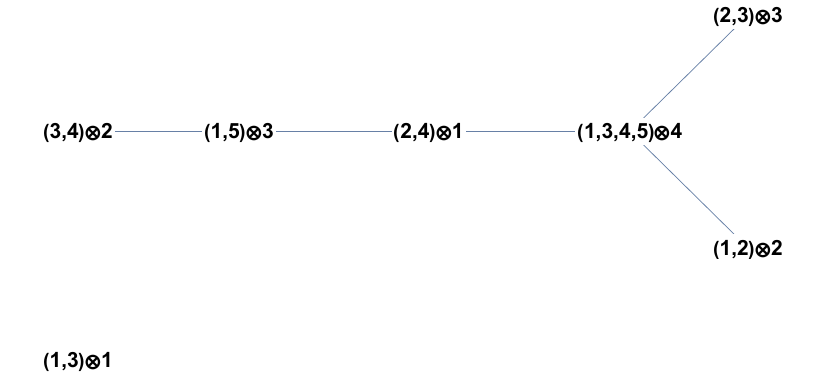}\\\hline
25&$\text{(1,2,4,5)$\otimes$4}\alb+\text{(1,3)$\otimes$1}\alb+\text{(1,5)$\otimes$3}\alb+\text{(2,3)$\otimes$3}\alb+\text{(2,4)$\otimes$1}\alb+\text{(3,4)$\otimes$2}\alb+\text{()$\otimes$4}$& $\ttt_1\alb+\u_4$ &55&$\begin{array}{cccccccc}
 3 & 1 & 3 & 4 & \oplus  & 4 & 4 & 1 \\
  &  & 1\\
\end{array}$&\includegraphics[scale=\imgscale]{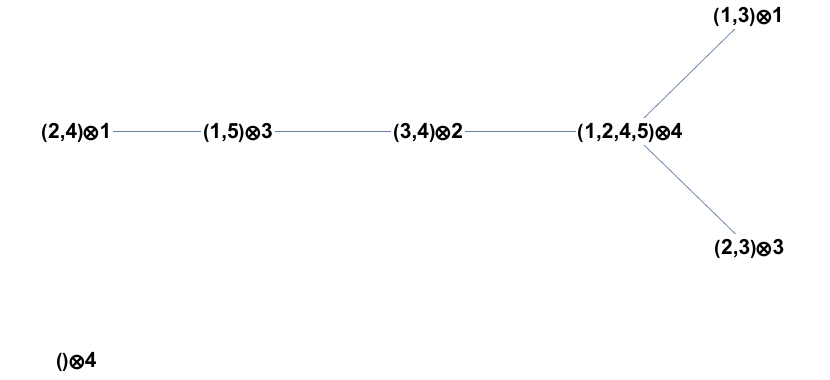}\\\hline
26&$\text{(1,2)$\otimes$1}\alb+\text{(1,3,4,5)$\otimes$4}\alb+\text{(2,3)$\otimes$3}\alb+\text{(2,4)$\otimes$2}\alb+\text{(3,4)$\otimes$2}\alb+\text{(4,5)$\otimes$3}\alb+\text{()$\otimes$4}$& $\ttt_1\alb+\u_4$ &55&$\begin{array}{cccccccc}
 1 & 6 & 1 & 1 & \oplus  & 1 & 1 & 6 \\
  &  & 0\\
\end{array}$&\includegraphics[scale=\imgscale]{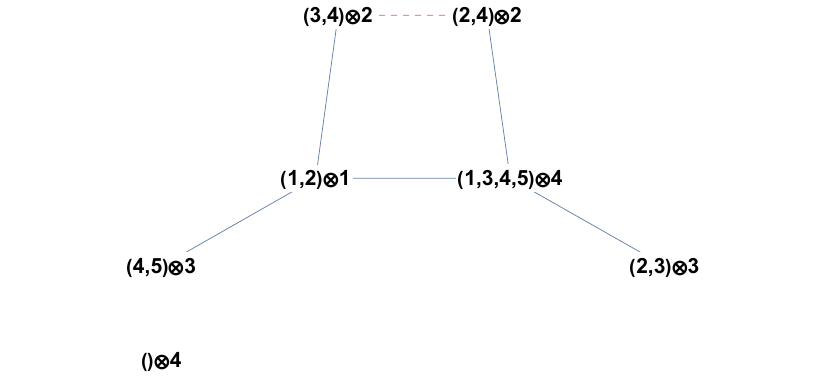}\\\hline
27&$\text{(1,2,4,5)$\otimes$4}\alb+\text{(1,3,4,5)$\otimes$1}\alb+\text{(1,4)$\otimes$2}\alb+\text{(2,3)$\otimes$3}\alb+\text{(2,4)$\otimes$1}\alb+\text{(3,4)$\otimes$2}\alb+\text{(4,5)$\otimes$3}\alb+\text{()$\otimes$4}$& $\u_6$ &54&$\begin{array}{cccccccc}
 0 & 4 & 0 & 0 & \oplus  & 0 & 4 & 4 \\
  &  & 4\\
\end{array}$&\includegraphics[scale=\imgscale]{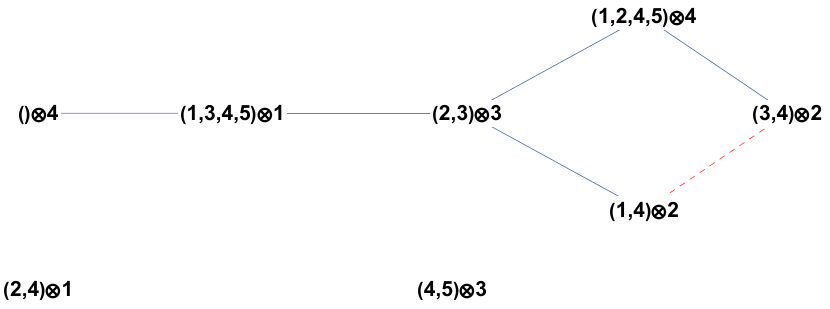}\\\hline
28&$\text{(1,2,3,4)$\otimes$2}\alb+\text{(1,2,4,5)$\otimes$1}\alb+\text{(1,2)$\otimes$4}\alb+\text{(1,3,4,5)$\otimes$4}\alb+\text{(1,3)$\otimes$3}\alb+\text{(2,4)$\otimes$3}\alb+\text{(4,5)$\otimes$2}\alb+\text{()$\otimes$1}$& $\u_6$ &54&$\begin{array}{cccccccc}
 0 & 0 & 0 & 8 & \oplus  & 0 & 0 & 0 \\
  &  & 0\\
\end{array}$&\includegraphics[scale=\imgscale]{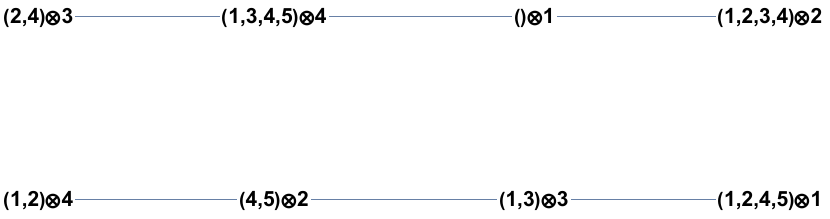}\\\hline
29&$\text{(1,2,3,4)$\otimes$2}\alb+\text{(1,3,4,5)$\otimes$4}\alb+\text{(1,3)$\otimes$1}\alb+\text{(1,5)$\otimes$3}\alb+\text{(2,3,4,5)$\otimes$3}\alb+\text{(2,4)$\otimes$1}\alb+\text{()$\otimes$4}$& $\ttt_1\alb+\u_5$ &54&$\begin{array}{cccccccc}
 3 & 2 & 1 & 2 & \oplus  & 5 & 3 & 3 \\
  &  & 2\\
\end{array}$&\includegraphics[scale=\imgscale]{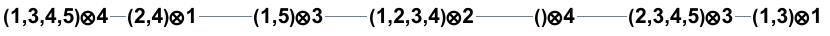}\\\hline
30&$\text{(1,2,4,5)$\otimes$1}\alb+\text{(1,2)$\otimes$4}\alb+\text{(1,3,4,5)$\otimes$4}\alb+\text{(1,3)$\otimes$1}\alb+\text{(2,3,4,5)$\otimes$3}\alb+\text{(3,4)$\otimes$2}\alb+\text{()$\otimes$3}$& $\ttt_2\alb+\u_4$ &54&$\begin{array}{cccccccc}
 0 & 6 & 2 & 0 & \oplus  & 2 & 0 & 6 \\
  &  & 0\\
\end{array}$&\includegraphics[scale=\imgscale]{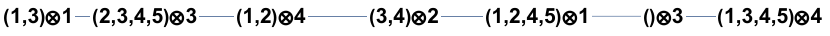}\\\hline
31&$\text{(1,2,4,5)$\otimes$3}\alb+\text{(1,2)$\otimes$1}\alb+\text{(1,3,4,5)$\otimes$4}\alb+\text{(2,3)$\otimes$3}\alb+\text{(3,4)$\otimes$2}\alb+\text{(4,5)$\otimes$3}\alb+\text{()$\otimes$4}$& $\ttt_1\alb+\u_5$ &54&$\begin{array}{cccccccc}
 3 & 2 & 0 & 3 & \oplus  & 3 & 0 & 5 \\
  &  & 3\\
\end{array}$&\includegraphics[scale=\imgscale]{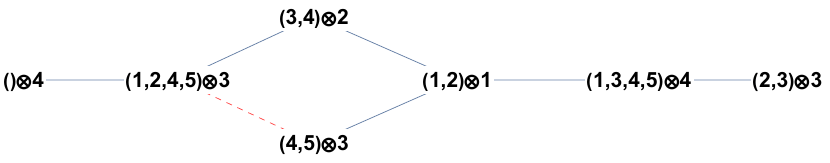}\\\hline
32&$\text{(1,2,3,4)$\otimes$2}\alb+\text{(1,2,4,5)$\otimes$4}\alb+\text{(1,3)$\otimes$1}\alb+\text{(1,4)$\otimes$2}\alb+\text{(1,5)$\otimes$3}\alb+\text{(2,3,4,5)$\otimes$3}\alb+\text{()$\otimes$4}$& $\ttt_1\alb+\u_5$ &54&$\begin{array}{cccccccc}
 1 & 0 & 1 & 0 & \oplus  & 1 & 7 & 1 \\
  &  & 6\\
\end{array}$&\includegraphics[scale=\imgscale]{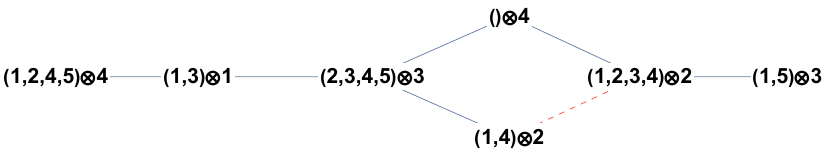}\\\hline
33&$\text{(1,2,4,5)$\otimes$4}\alb+\text{(1,3,4,5)$\otimes$1}\alb+\text{(1,4)$\otimes$2}\alb+\text{(1,5)$\otimes$3}\alb+\text{(2,3)$\otimes$3}\alb+\text{(2,4)$\otimes$1}\alb+\text{()$\otimes$4}$& $\ttt_1\alb+\u_5$ &54&$\begin{array}{cccccccc}
 2 & 0 & 2 & 2 & \oplus  & 2 & 2 & 6 \\
  &  & 2\\
\end{array}$&\includegraphics[scale=\imgscale]{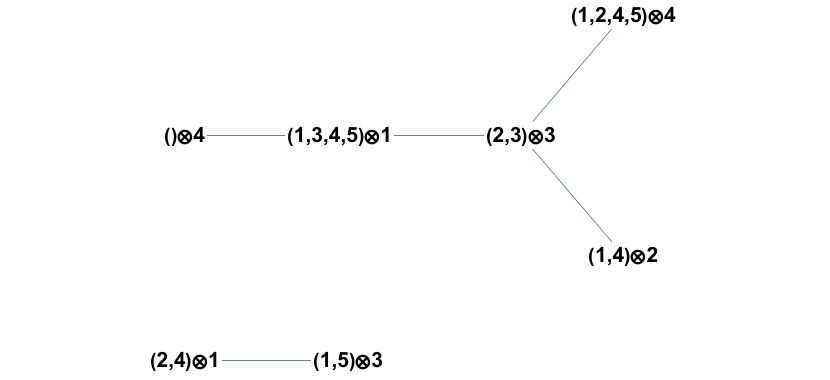}\\\hline
34&$\text{(1,3,4,5)$\otimes$1}\alb+\text{(1,4)$\otimes$2}\alb+\text{(1,5)$\otimes$3}\alb+\text{(2,3,4,5)$\otimes$3}\alb+\text{(2,3)$\otimes$4}\alb+\text{(2,4)$\otimes$1}$& $A_1\alb+\ttt_1\alb+\u_2$&54&$\begin{array}{cccccccc}
 4 & 2 & 0 & 2 & \oplus  & 4 & 4 & 8 \\
  &  & 6\\
\end{array}$&\includegraphics[scale=\imgscale]{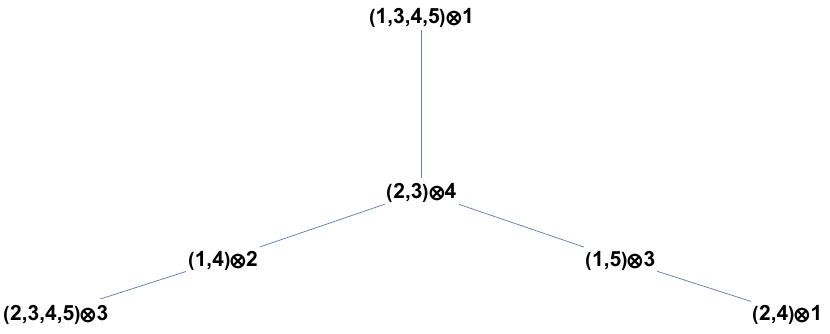}\\\hline
35&$\text{(1,2,4,5)$\otimes$4}\alb+\text{(1,3)$\otimes$1}\alb+\text{(1,5)$\otimes$3}\alb+\text{(2,3)$\otimes$3}\alb+\text{(2,4)$\otimes$1}\alb+\text{(3,4)$\otimes$2}$& $\ttt_2\alb+\u_4$ &54&$\begin{array}{cccccccc}
 2 & 1 & 4 & 3 & \oplus  & 5 & 3 & 1 \\
  &  & 1\\
\end{array}$&\includegraphics[scale=\imgscale]{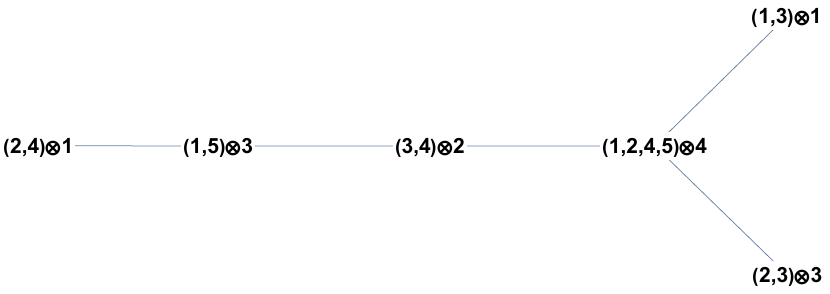}\\\hline
36&$\text{(1,3,4,5)$\otimes$4}\alb+\text{(1,3)$\otimes$1}\alb+\text{(1,4)$\otimes$2}\alb+\text{(1,5)$\otimes$3}\alb+\text{(2,3,4,5)$\otimes$3}\alb+\text{(2,3)$\otimes$4}\alb+\text{(2,4)$\otimes$1}$& $A_1\alb+\u_4$&53&$\begin{array}{cccccccc}
 4 & 2 & 0 & 2 & \oplus  & 4 & 4 & 4 \\
  &  & 2\\
\end{array}$&\includegraphics[scale=\imgscale]{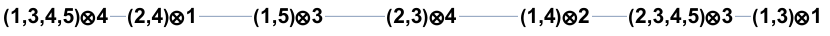}\\\hline
37&$\text{(1,2,4,5)$\otimes$4}\alb+\text{(1,3,4,5)$\otimes$1}\alb+\text{(1,4)$\otimes$2}\alb+\text{(2,3)$\otimes$3}\alb+\text{(2,4)$\otimes$1}\alb+\text{(4,5)$\otimes$3}\alb+\text{()$\otimes$4}$& $\ttt_1\alb+\u_6$ &53&$\begin{array}{cccccccc}
 1 & 2 & 1 & 1 & \oplus  & 1 & 3 & 5 \\
  &  & 3\\
\end{array}$&\includegraphics[scale=\imgscale]{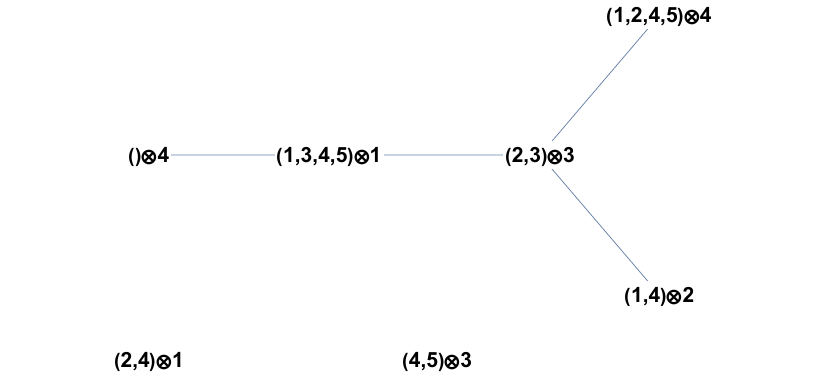}\\\hline
38&$\text{(1,2,3,4)$\otimes$2}\alb+\text{(1,2,4,5)$\otimes$3}\alb+\text{(1,2)$\otimes$4}\alb+\text{(1,3,4,5)$\otimes$4}\alb+\text{(1,3)$\otimes$1}\alb+\text{(4,5)$\otimes$1}\alb+\text{()$\otimes$3}$& $\ttt_1\alb+\u_6$ &53&$\begin{array}{cccccccc}
 0 & 1 & 0 & 7 & \oplus  & 0 & 0 & 1 \\
  &  & 0\\
\end{array}$&\includegraphics[scale=\imgscale]{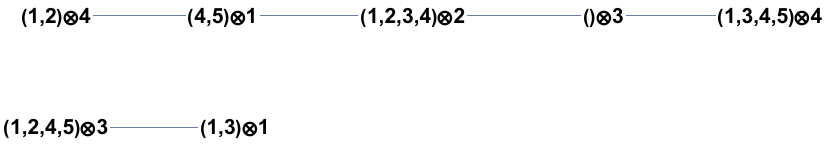}\\\hline
39&$\text{(1,2,3,4)$\otimes$2}\alb+\text{(1,2,4,5)$\otimes$4}\alb+\text{(1,3)$\otimes$1}\alb+\text{(2,3)$\otimes$3}\alb+\text{(2,4)$\otimes$1}\alb+\text{(4,5)$\otimes$3}\alb+\text{()$\otimes$4}$& $\ttt_1\alb+\u_6$ &53&$\begin{array}{cccccccc}
 2 & 1 & 2 & 2 & \oplus  & 2 & 3 & 3 \\
  &  & 1\\
\end{array}$&\includegraphics[scale=\imgscale]{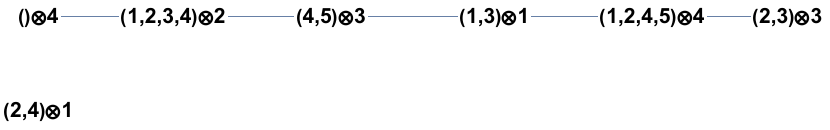}\\\hline
40&$\text{(1,2,4,5)$\otimes$3}\alb+\text{(1,2)$\otimes$4}\alb+\text{(1,3,4,5)$\otimes$1}\alb+\text{(2,3)$\otimes$3}\alb+\text{(2,4)$\otimes$1}\alb+\text{(3,4)$\otimes$2}\alb+\text{(4,5)$\otimes$3}$& $\ttt_1\alb+\u_6$ &53&$\begin{array}{cccccccc}
 1 & 3 & 0 & 1 & \oplus  & 1 & 3 & 4 \\
  &  & 4\\
\end{array}$&\includegraphics[scale=\imgscale]{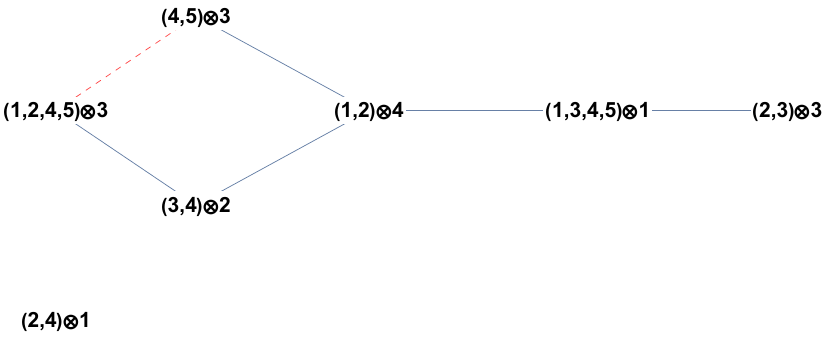}\\\hline
41&$\text{(1,2,4,5)$\otimes$4}\alb+\text{(1,2)$\otimes$3}\alb+\text{(1,3,4,5)$\otimes$1}\alb+\text{(2,3,4,5)$\otimes$4}\alb+\text{(3,4)$\otimes$2}\alb+\text{(4,5)$\otimes$3}\alb+\text{()$\otimes$4}$& $\ttt_1\alb+\u_6$ &53&$\begin{array}{cccccccc}
 0 & 3 & 0 & 0 & \oplus  & 0 & 5 & 3 \\
  &  & 5\\
\end{array}$&\includegraphics[scale=\imgscale]{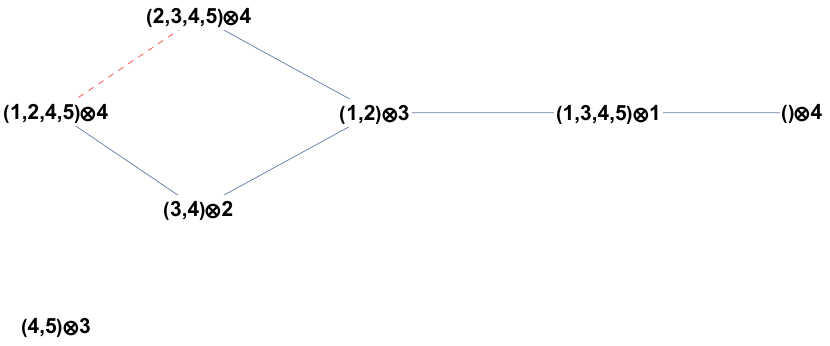}\\\hline
42&$\text{(1,2,4,5)$\otimes$3}\alb+\text{(1,2)$\otimes$1}\alb+\text{(1,3,4,5)$\otimes$4}\alb+\text{(2,3)$\otimes$3}\alb+\text{(3,4)$\otimes$2}\alb+\text{()$\otimes$4}$& $\ttt_2\alb+\u_5$ &53&$\begin{array}{cccccccc}
 4 & 0 & 4 & 0 & \oplus  & 0 & 4 & 0 \\
  &  & 0\\
\end{array}$&\includegraphics[scale=\imgscale]{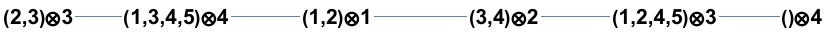}\\\hline
43&$\text{(1,2,3,4)$\otimes$2}\alb+\text{(1,2)$\otimes$4}\alb+\text{(1,3,4,5)$\otimes$4}\alb+\text{(1,3)$\otimes$1}\alb+\text{(2,3)$\otimes$3}\alb+\text{(2,4)$\otimes$1}\alb+\text{(4,5)$\otimes$3}\alb+\text{()$\otimes$4}$& $\u_8$ &52&$\begin{array}{cccccccc}
 2 & 2 & 0 & 2 & \oplus  & 2 & 2 & 2 \\
  &  & 2\\
\end{array}$&\includegraphics[scale=\imgscale]{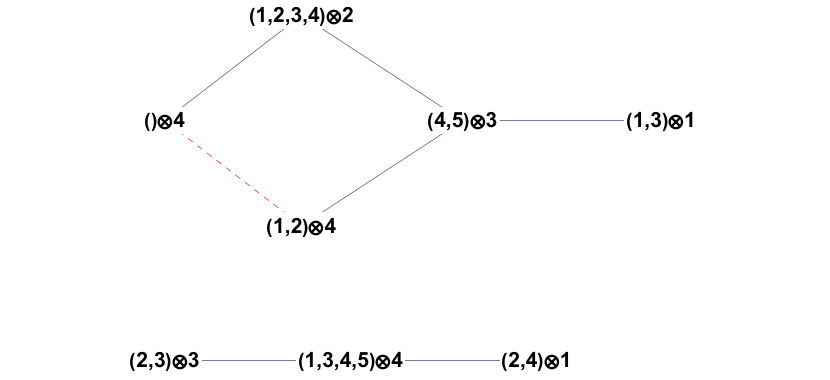}\\\hline
44&$\text{(1,2,3,4)$\otimes$2}\alb+\text{(1,2,4,5)$\otimes$1}\alb+\text{(1,2)$\otimes$4}\alb+\text{(1,3,4,5)$\otimes$4}\alb+\text{(1,3)$\otimes$1}\alb+\text{(2,3,4,5)$\otimes$3}\alb+\text{(2,4)$\otimes$4}\alb+\text{()$\otimes$3}$& $\u_8$ &52&$\begin{array}{cccccccc}
 0 & 0 & 4 & 0 & \oplus  & 4 & 0 & 4 \\
  &  & 0\\
\end{array}$&\includegraphics[scale=\imgscale]{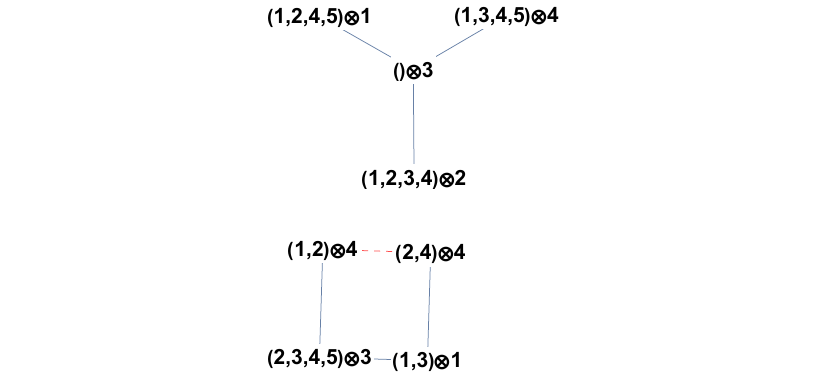}\\\hline
45&$\text{(1,2,4,5)$\otimes$3}\alb+\text{(1,2)$\otimes$4}\alb+\text{(1,3,4,5)$\otimes$1}\alb+\text{(2,4)$\otimes$1}\alb+\text{(3,4)$\otimes$2}\alb+\text{(4,5)$\otimes$4}\alb+\text{()$\otimes$3}$& $\ttt_1\alb+\u_7$ &52&$\begin{array}{cccccccc}
 0 & 1 & 0 & 6 & \oplus  & 0 & 1 & 1 \\
  &  & 1\\
\end{array}$&\includegraphics[scale=\imgscale]{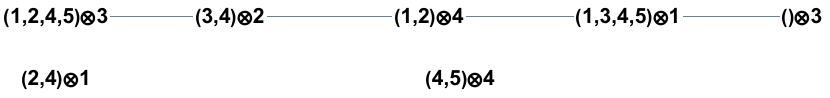}\\\hline
46&$\text{(1,2,4,5)$\otimes$4}\alb+\text{(1,2)$\otimes$1}\alb+\text{(1,3,4,5)$\otimes$3}\alb+\text{(1,3)$\otimes$2}\alb+\text{(2,4)$\otimes$2}\alb+\text{(3,4)$\otimes$1}\alb+\text{()$\otimes$3}$& $\ttt_2\alb+\u_6$ &52&$\begin{array}{cccccccc}
 1 & 0 & 1 & 6 & \oplus  & 0 & 1 & 0 \\
  &  & 0\\
\end{array}$&\includegraphics[scale=\imgscale]{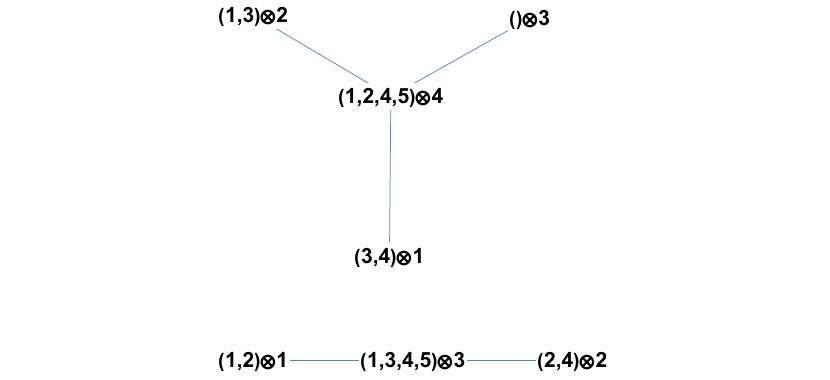}\\\hline
47&$\text{(1,2,3,4)$\otimes$2}\alb+\text{(1,2,4,5)$\otimes$4}\alb+\text{(1,5)$\otimes$3}\alb+\text{(2,3)$\otimes$3}\alb+\text{(3,4)$\otimes$1}\alb+\text{()$\otimes$4}$& $\ttt_2\alb+\u_6$ &52&$\begin{array}{cccccccc}
 2 & 0 & 2 & 2 & \oplus  & 2 & 4 & 2 \\
  &  & 2\\
\end{array}$&\includegraphics[scale=\imgscale]{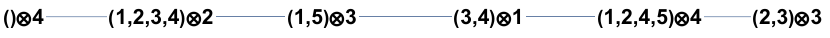}\\\hline
48&$\text{(1,2,4,5)$\otimes$3}\alb+\text{(1,2)$\otimes$4}\alb+\text{(1,3,4,5)$\otimes$1}\alb+\text{(2,3)$\otimes$3}\alb+\text{(3,4)$\otimes$2}\alb+\text{(4,5)$\otimes$3}$& $\ttt_2\alb+\u_6$&52&$\begin{array}{cccccccc}
 1 & 2 & 0 & 1 & \oplus  & 1 & 4 & 3 \\
  &  & 5\\
\end{array}$&\includegraphics[scale=\imgscale]{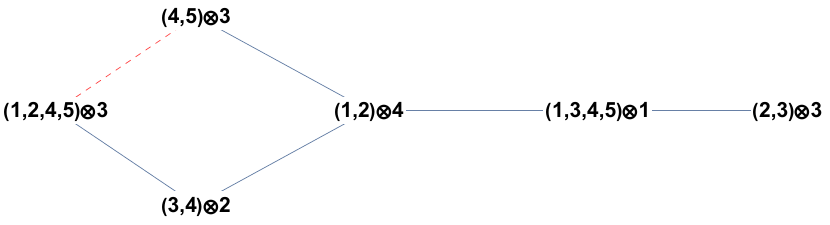}\\\hline
49&$\text{(1,2,4,5)$\otimes$4}\alb+\text{(1,3,4,5)$\otimes$1}\alb+\text{(1,4)$\otimes$2}\alb+\text{(2,3)$\otimes$3}\alb+\text{(4,5)$\otimes$3}\alb+\text{()$\otimes$4}$& $\ttt_2\alb+\u_6$ &52&$\begin{array}{cccccccc}
 1 & 1 & 1 & 1 & \oplus  & 1 & 4 & 4 \\
  &  & 4\\
\end{array}$&\includegraphics[scale=\imgscale]{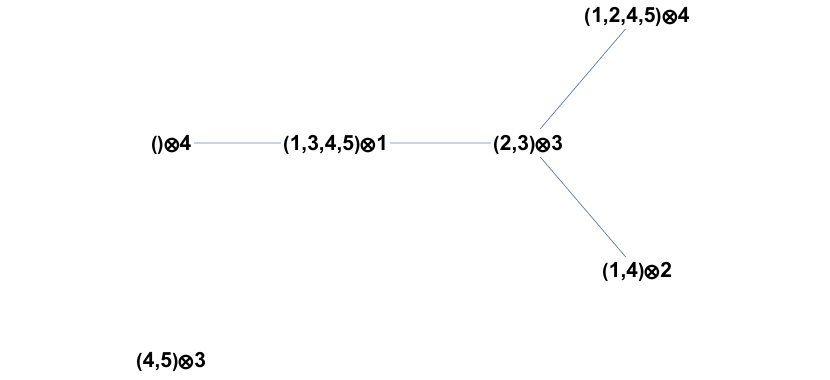}\\\hline
50&$\text{(1,2,3,4)$\otimes$2}\alb+\text{(1,2,4,5)$\otimes$3}\alb+\text{(1,3,4,5)$\otimes$1}\alb+\text{(2,3)$\otimes$3}\alb+\text{(2,4)$\otimes$1}\alb+\text{()$\otimes$4}$& $\ttt_2\alb+\u_6$ &52&$\begin{array}{cccccccc}
 2 & 1 & 2 & 0 & \oplus  & 0 & 3 & 5 \\
  &  & 3\\
\end{array}$&\includegraphics[scale=\imgscale]{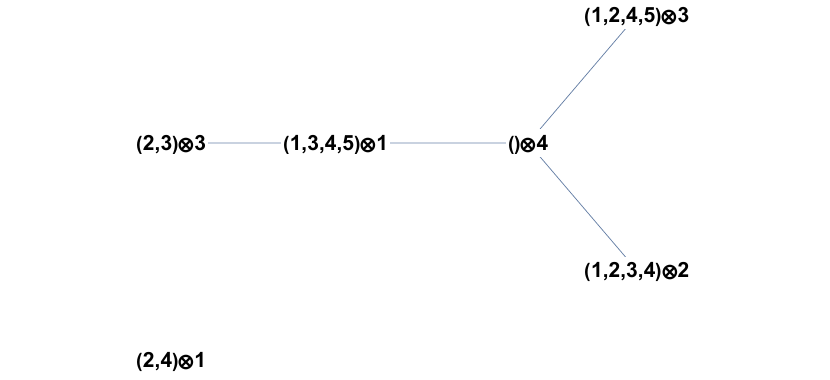}\\\hline
51&$\text{(1,2)$\otimes$4}\alb+\text{(1,3,4,5)$\otimes$4}\alb+\text{(1,3)$\otimes$1}\alb+\text{(1,4)$\otimes$2}\alb+\text{(1,5)$\otimes$3}\alb+\text{(2,3)$\otimes$3}\alb+\text{(2,4)$\otimes$1}$& $\ttt_1\alb+\u_8$ &51&$\begin{array}{cccccccc}
 1 & 1 & 1 & 2 & \oplus  & 3 & 1 & 4 \\
  &  & 1\\
\end{array}$&\includegraphics[scale=\imgscale]{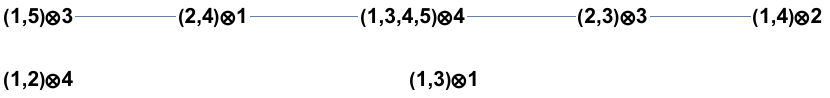}\\\hline
52&$\text{(1,2,3,4)$\otimes$2}\alb+\text{(1,2,3,5)$\otimes$3}\alb+\text{(1,2)$\otimes$4}\alb+\text{(1,3,4,5)$\otimes$4}\alb+\text{(1,3)$\otimes$1}\alb+\text{(2,4)$\otimes$4}\alb+\text{(4,5)$\otimes$3}$&  $\ttt_1\alb+\u_8$ &51&$\begin{array}{cccccccc}
 1 & 0 & 3 & 0 & \oplus  & 4 & 1 & 3 \\
  &  & 1\\
\end{array}$&\includegraphics[scale=\imgscale]{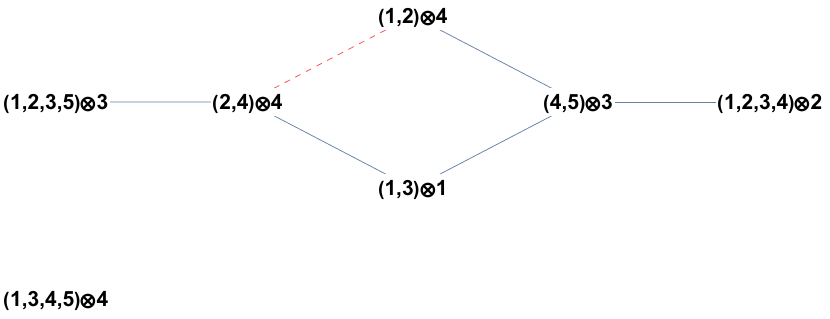}\\\hline
53&$\text{(1,2,3,4)$\otimes$2}\alb+\text{(1,2)$\otimes$4}\alb+\text{(1,3,4,5)$\otimes$4}\alb+\text{(1,3)$\otimes$1}\alb+\text{(2,4)$\otimes$1}\alb+\text{(3,4)$\otimes$2}\alb+\text{(4,5)$\otimes$3}$& $\ttt_1\alb+\u_8$ &51&$\begin{array}{cccccccc}
 1 & 1 & 0 & 5 & \oplus  & 1 & 1 & 1 \\
  &  & 1\\
\end{array}$&\includegraphics[scale=\imgscale]{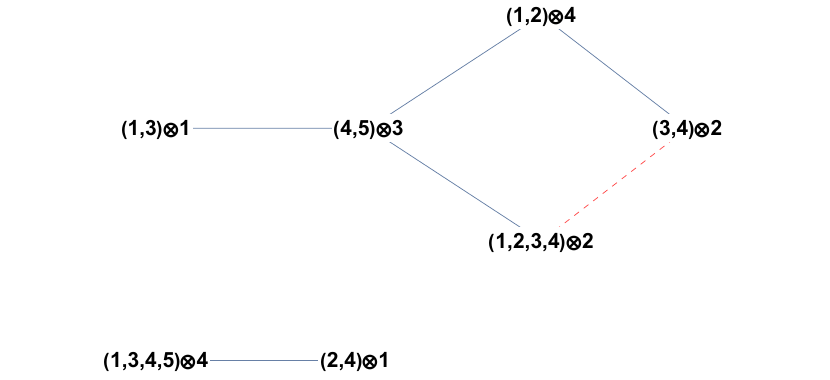}\\\hline
54&$\text{(1,2,3,4)$\otimes$2}\alb+\text{(1,2)$\otimes$4}\alb+\text{(1,3,4,5)$\otimes$4}\alb+\text{(1,3)$\otimes$1}\alb+\text{(2,3)$\otimes$3}\alb+\text{(2,4)$\otimes$1}\alb+\text{(4,5)$\otimes$3}$&  $\ttt_1\alb+\u_8$&51&$\begin{array}{cccccccc}
 1 & 1 & 2 & 1 & \oplus  & 3 & 1 & 3 \\
  &  & 1\\
\end{array}$&\includegraphics[scale=\imgscale]{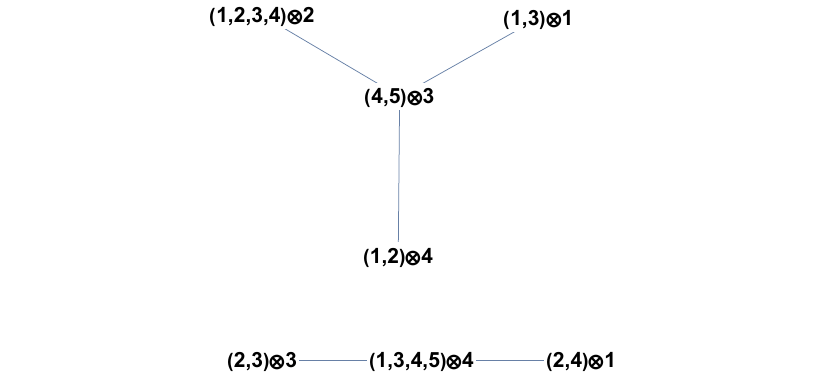}\\\hline
55&$\text{(1,2,3,4)$\otimes$1}\alb+\text{(1,4)$\otimes$2}\alb+\text{(1,5)$\otimes$3}\alb+\text{(2,3,4,5)$\otimes$3}\alb+\text{(2,3)$\otimes$4}\alb+\text{(4,5)$\otimes$4}$& $2A_1\alb+\u_3$ &51&$\begin{array}{cccccccc}
 0 & 0 & 4 & 0 & \oplus  & 0 & 8 & 0 \\
  &  & 0\\
\end{array}$&\includegraphics[scale=\imgscale]{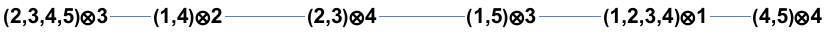}\\\hline
56&$\text{(1,2)$\otimes$4}\alb+\text{(1,3)$\otimes$1}\alb+\text{(1,5)$\otimes$3}\alb+\text{(2,3,4,5)$\otimes$3}\alb+\text{(2,4)$\otimes$1}\alb+\text{(3,4)$\otimes$2}$& $2A_1\alb+\u_3$ &51&$\begin{array}{cccccccc}
 8 & 0 & 0 & 0 & \oplus  & 8 & 0 & 0 \\
  &  & 0\\
\end{array}$&\includegraphics[scale=\imgscale]{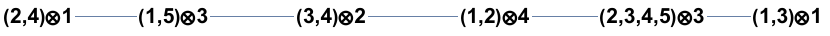}\\\hline
57&$\text{(1,2,4,5)$\otimes$3}\alb+\text{(1,2)$\otimes$4}\alb+\text{(1,3,4,5)$\otimes$1}\alb+\text{(2,4)$\otimes$1}\alb+\text{(3,4)$\otimes$2}\alb+\text{()$\otimes$3}$& $\ttt_2\alb+\u_7$ &51&$\begin{array}{cccccccc}
 1 & 0 & 1 & 5 & \oplus  & 1 & 0 & 1 \\
  &  & 1\\
\end{array}$&\includegraphics[scale=\imgscale]{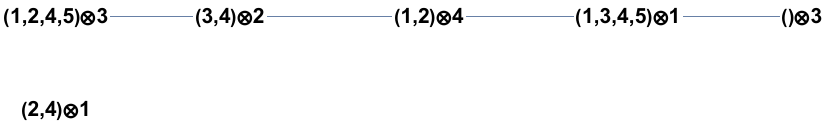}\\\hline
58&$\text{(1,2,4,5)$\otimes$4}\alb+\text{(1,2)$\otimes$1}\alb+\text{(1,3,4,5)$\otimes$1}\alb+\text{(2,3)$\otimes$3}\alb+\text{(4,5)$\otimes$3}\alb+\text{()$\otimes$4}$& $2A_1\alb+\u_3$ &51&$\begin{array}{cccccccc}
 0 & 0 & 0 & 0 & \oplus  & 0 & 0 & 8 \\
  &  & 0\\
\end{array}$&\includegraphics[scale=\imgscale]{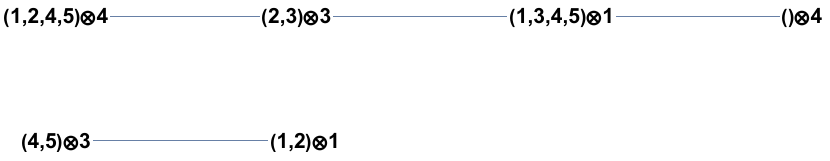}\\\hline
59&$\text{(1,2,3,4)$\otimes$2}\alb+\text{(1,2)$\otimes$4}\alb+\text{(1,3,4,5)$\otimes$4}\alb+\text{(1,3)$\otimes$1}\alb+\text{(2,3)$\otimes$3}\alb+\text{(3,4)$\otimes$1}\alb+\text{(4,5)$\otimes$3}$& $\ttt_1\alb+\u_9$ &50&$\begin{array}{cccccccc}
 2 & 0 & 1 & 2 & \oplus  & 3 & 2 & 1 \\
  &  & 2\\
\end{array}$&\includegraphics[scale=\imgscale]{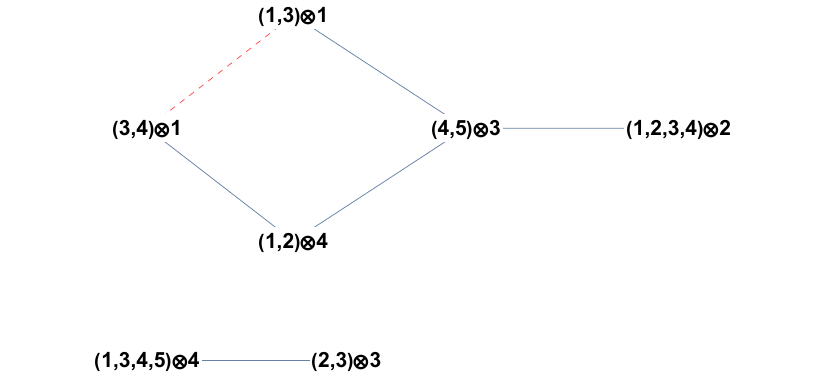}\\\hline
60&$\text{(1,2)$\otimes$1}\alb+\text{(1,3,4,5)$\otimes$3}\alb+\text{(1,3)$\otimes$2}\alb+\text{(2,3)$\otimes$3}\alb+\text{(2,4)$\otimes$4}\alb+\text{(3,4)$\otimes$1}\alb+\text{()$\otimes$2}$& $A_1\alb+\u_7$ &50&$\begin{array}{cccccccc}
 4 & 0 & 0 & 4 & \oplus  & 4 & 0 & 0 \\
  &  & 0\\
\end{array}$&\includegraphics[scale=\imgscale]{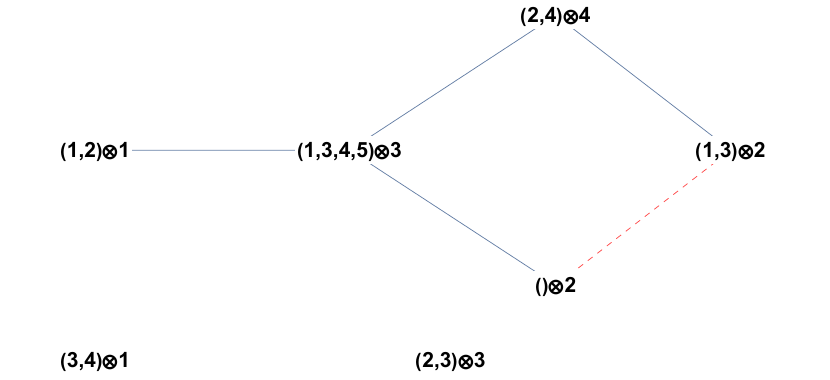}\\\hline
61&$\text{(1,2,3,4)$\otimes$2}\alb+\text{(1,2,3,5)$\otimes$3}\alb+\text{(1,2)$\otimes$1}\alb+\text{(1,3)$\otimes$4}\alb+\text{(2,4)$\otimes$4}\alb+\text{(3,4)$\otimes$1}\alb+\text{(4,5)$\otimes$3}$& $\ttt_2\alb+\u_8$ &50&$\begin{array}{cccccccc}
 2 & 0 & 3 & 0 & \oplus  & 5 & 0 & 3 \\
  &  & 0\\
\end{array}$&\includegraphics[scale=\imgscale]{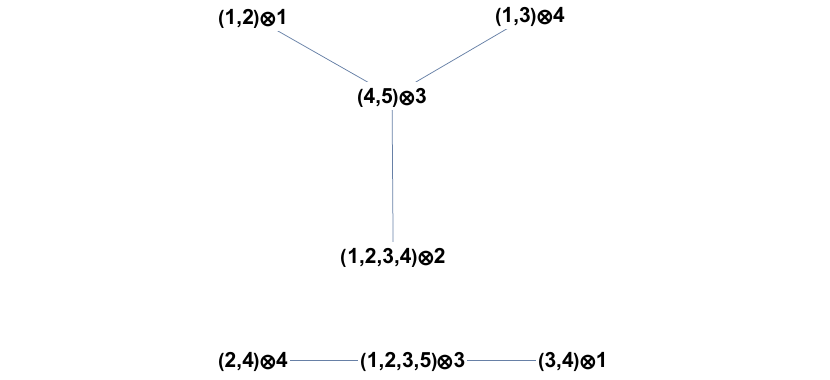}\\\hline
62&$\text{(1,2,4,5)$\otimes$3}\alb+\text{(1,3,4,5)$\otimes$4}\alb+\text{(1,3)$\otimes$1}\alb+\text{(1,4)$\otimes$2}\alb+\text{(2,3)$\otimes$3}\alb+\text{(2,4)$\otimes$1}\alb+\text{()$\otimes$4}$& $\ttt_1\alb+\u_9$ &50&$\begin{array}{cccccccc}
 2 & 1 & 1 & 1 & \oplus  & 1 & 2 & 3 \\
  &  & 1\\
\end{array}$&\includegraphics[scale=\imgscale]{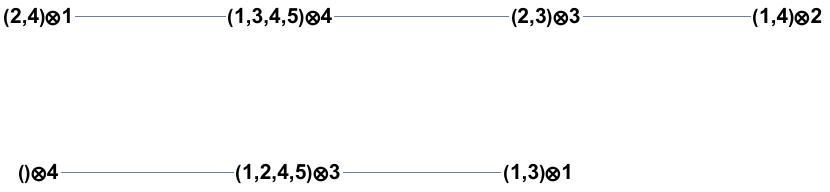}\\\hline
63&$\text{(1,2,4,5)$\otimes$3}\alb+\text{(1,2)$\otimes$4}\alb+\text{(1,3,4,5)$\otimes$4}\alb+\text{(1,3)$\otimes$1}\alb+\text{(2,4)$\otimes$1}\alb+\text{(3,4)$\otimes$2}\alb+\text{()$\otimes$3}$&  $\ttt_1\alb+\u_9$&50&$\begin{array}{cccccccc}
 1 & 1 & 1 & 3 & \oplus  & 1 & 1 & 1 \\
  &  & 1\\
\end{array}$&\includegraphics[scale=\imgscale]{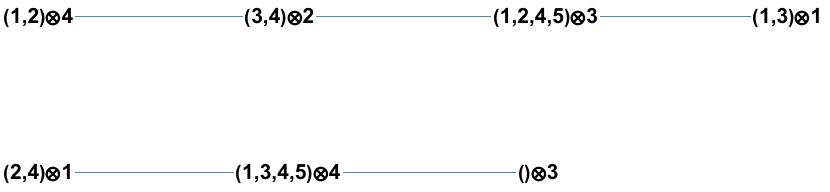}\\\hline
64&$\text{(1,2,3,4)$\otimes$1}\alb+\text{(1,2,4,5)$\otimes$4}\alb+\text{(1,3,4,5)$\otimes$2}\alb+\text{(1,5)$\otimes$3}\alb+\text{(2,4)$\otimes$3}\alb+\text{()$\otimes$4}$& $A_1\alb+\ttt_1\alb+\u_6$ &50&$\begin{array}{cccccccc}
 0 & 0 & 0 & 6 & \oplus  & 0 & 2 & 0 \\
  &  & 2\\
\end{array}$&\includegraphics[scale=\imgscale]{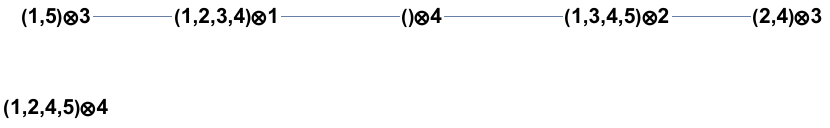}\\\hline
65&$\text{(1,2)$\otimes$4}\alb+\text{(1,3,4,5)$\otimes$4}\alb+\text{(1,4)$\otimes$2}\alb+\text{(1,5)$\otimes$3}\alb+\text{(2,3,4,5)$\otimes$3}\alb+\text{(3,4)$\otimes$1}$&  $\ttt_2\alb+\u_8$&50&$\begin{array}{cccccccc}
 0 & 2 & 1 & 2 & \oplus  & 3 & 2 & 3 \\
  &  & 0\\
\end{array}$&\includegraphics[scale=\imgscale]{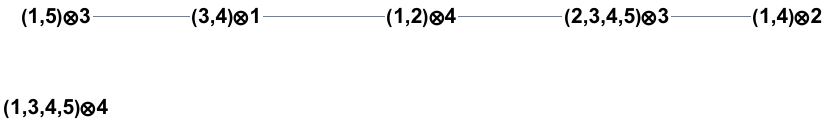}\\\hline
66&$\text{(1,2,3,4)$\otimes$2}\alb+\text{(1,2)$\otimes$4}\alb+\text{(1,3,4,5)$\otimes$4}\alb+\text{(1,3)$\otimes$1}\alb+\text{(2,4)$\otimes$1}\alb+\text{(4,5)$\otimes$3}$& $\ttt_2\alb+\u_8$ &50&$\begin{array}{cccccccc}
 0 & 0 & 2 & 4 & \oplus  & 2 & 0 & 2 \\
  &  & 0\\
\end{array}$&\includegraphics[scale=\imgscale]{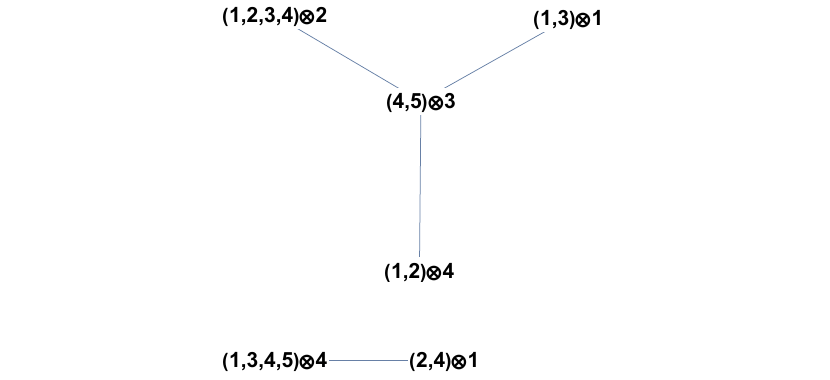}\\\hline
67&$\text{(1,2)$\otimes$4}\alb+\text{(1,3,4,5)$\otimes$4}\alb+\text{(1,3)$\otimes$1}\alb+\text{(1,5)$\otimes$3}\alb+\text{(2,3)$\otimes$3}\alb+\text{(2,4)$\otimes$1}$& $\ttt_2\alb+\u_8$ &50&$\begin{array}{cccccccc}
 0 & 0 & 1 & 0 & \oplus  & 1 & 0 & 7 \\
  &  & 0\\
\end{array}$&\includegraphics[scale=\imgscale]{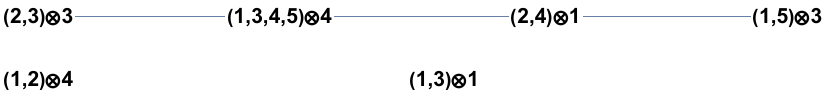}\\\hline
68&$\text{(1,2,3,4)$\otimes$2}\alb+\text{(1,2,4,5)$\otimes$3}\alb+\text{(1,3,4,5)$\otimes$1}\alb+\text{(2,3)$\otimes$3}\alb+\text{()$\otimes$4}$& $A_1\alb+\ttt_2\alb+\u_5$ &50&$\begin{array}{cccccccc}
 2 & 0 & 2 & 0 & \oplus  & 0 & 4 & 4 \\
  &  & 4\\
\end{array}$&\includegraphics[scale=\imgscale]{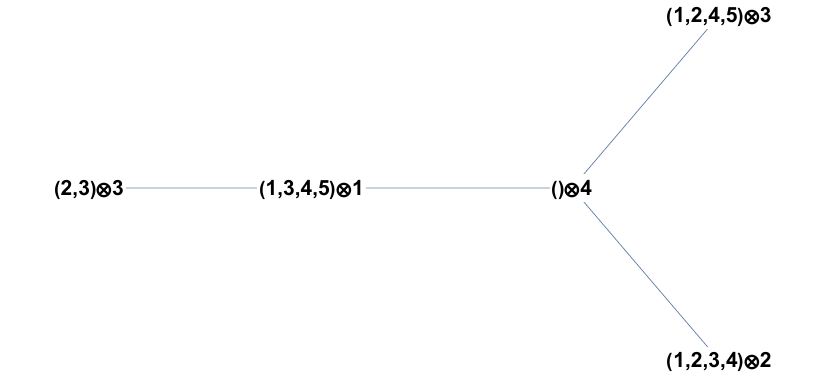}\\\hline
69&$\text{(1,2,3,4)$\otimes$2}\alb+\text{(1,2,4,5)$\otimes$3}\alb+\text{(1,2)$\otimes$4}\alb+\text{(1,3,4,5)$\otimes$4}\alb+\text{(1,3)$\otimes$1}\alb+\text{(2,4)$\otimes$1}\alb+\text{()$\otimes$3}$& $\ttt_1\alb+\u_{10}$ &49&$\begin{array}{cccccccc}
 1 & 1 & 1 & 2 & \oplus  & 1 & 1 & 2 \\
  &  & 1\\
\end{array}$&\includegraphics[scale=\imgscale]{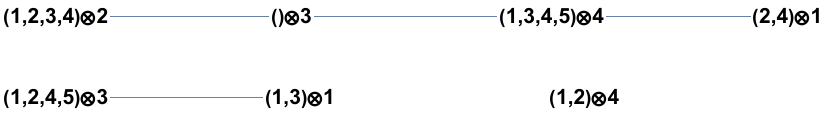}\\\hline
70&$\text{(1,2,3,4)$\otimes$2}\alb+\text{(1,2)$\otimes$4}\alb+\text{(1,3,4,5)$\otimes$3}\alb+\text{(1,3)$\otimes$1}\alb+\text{(2,4)$\otimes$1}\alb+\text{(3,4)$\otimes$4}\alb+\text{()$\otimes$2}$& $A_1\alb+\ttt_1\alb+\u_7$ &49&$\begin{array}{cccccccc}
 3 & 0 & 0 & 5 & \oplus  & 3 & 0 & 0 \\
  &  & 0\\
\end{array}$&\includegraphics[scale=\imgscale]{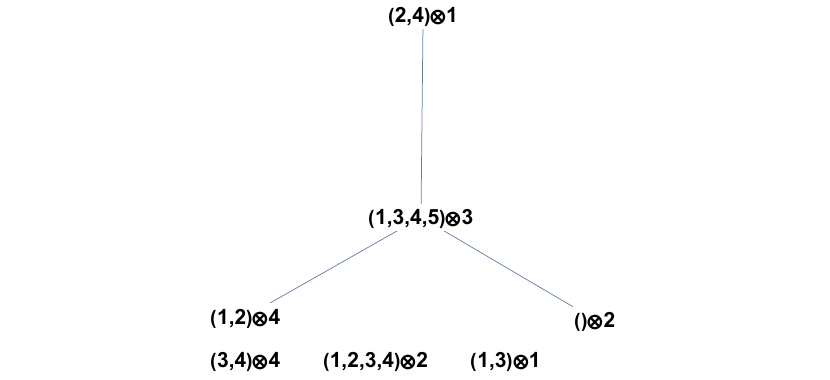}\\\hline
71&$\text{(1,3,4,5)$\otimes$4}\alb+\text{(1,3)$\otimes$1}\alb+\text{(1,4)$\otimes$2}\alb+\text{(1,5)$\otimes$3}\alb+\text{(2,3,4,5)$\otimes$3}\alb+\text{(2,4)$\otimes$4}$& $A_1\alb+\ttt_1\alb+\u_7$ &49&$\begin{array}{cccccccc}
 0 & 2 & 0 & 2 & \oplus  & 4 & 0 & 4 \\
  &  & 2\\
\end{array}$&\includegraphics[scale=\imgscale]{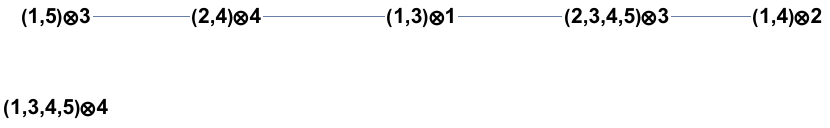}\\\hline
72&$\text{(1,2,3,4)$\otimes$2}\alb+\text{(1,2)$\otimes$1}\alb+\text{(1,3)$\otimes$4}\alb+\text{(2,3)$\otimes$3}\alb+\text{(2,4)$\otimes$1}\alb+\text{(4,5)$\otimes$3}$& $\ttt_2\alb+\u_9$ &49&$\begin{array}{cccccccc}
 3 & 0 & 1 & 3 & \oplus  & 4 & 0 & 1 \\
  &  & 0\\
\end{array}$&\includegraphics[scale=\imgscale]{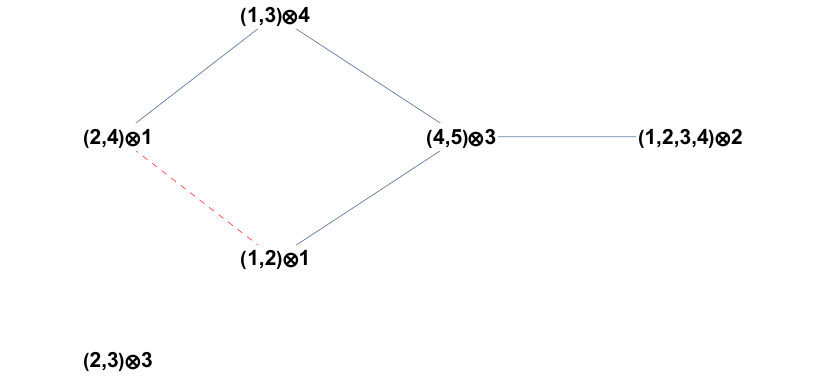}\\\hline
73&$\text{(1,2,3,4)$\otimes$2}\alb+\text{(1,2)$\otimes$4}\alb+\text{(1,3,4,5)$\otimes$4}\alb+\text{(1,3)$\otimes$1}\alb+\text{(3,4)$\otimes$1}\alb+\text{(4,5)$\otimes$3}$& $\ttt_2\alb+\u_9$ &49&$\begin{array}{cccccccc}
 1 & 0 & 1 & 4 & \oplus  & 2 & 1 & 1 \\
  &  & 1\\
\end{array}$&\includegraphics[scale=\imgscale]{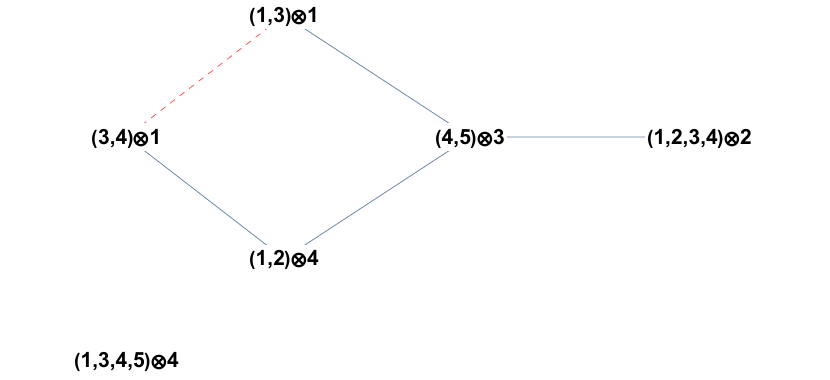}\\\hline
74&$\text{(1,2,4,5)$\otimes$3}\alb+\text{(1,3)$\otimes$3}\alb+\text{(1,4)$\otimes$1}\alb+\text{(2,3,4,5)$\otimes$4}\alb+\text{(2,4)$\otimes$2}\alb+\text{()$\otimes$4}$& $\ttt_2\alb+\u_9$&49&$\begin{array}{cccccccc}
 0 & 4 & 0 & 0 & \oplus  & 0 & 4 & 0 \\
  &  & 0\\
\end{array}$&\includegraphics[scale=\imgscale]{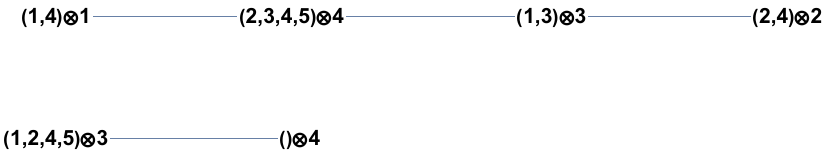}\\\hline
75&$\text{(1,2,4,5)$\otimes$3}\alb+\text{(1,3,4,5)$\otimes$1}\alb+\text{(1,3)$\otimes$2}\alb+\text{(2,4)$\otimes$4}\alb+\text{()$\otimes$3}$& $A_1\alb+\ttt_2\alb+\u_6$&49&$\begin{array}{cccccccc}
 0 & 1 & 0 & 5 & \oplus  & 1 & 1 & 0 \\
  &  & 2\\
\end{array}$&\includegraphics[scale=\imgscale]{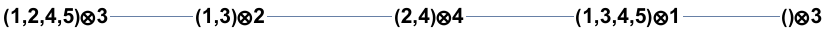}\\\hline
76&$\text{(1,2,4,5)$\otimes$3}\alb+\text{(1,2)$\otimes$4}\alb+\text{(1,3,4,5)$\otimes$4}\alb+\text{(1,4)$\otimes$2}\alb+\text{(2,3)$\otimes$3}\alb+\text{(3,4)$\otimes$1}\alb+\text{(4,5)$\otimes$3}$& $\ttt_1\alb+\u_{11}$ &48&$\begin{array}{cccccccc}
 1 & 2 & 0 & 1 & \oplus  & 1 & 3 & 1 \\
  &  & 1\\
\end{array}$&\includegraphics[scale=\imgscale]{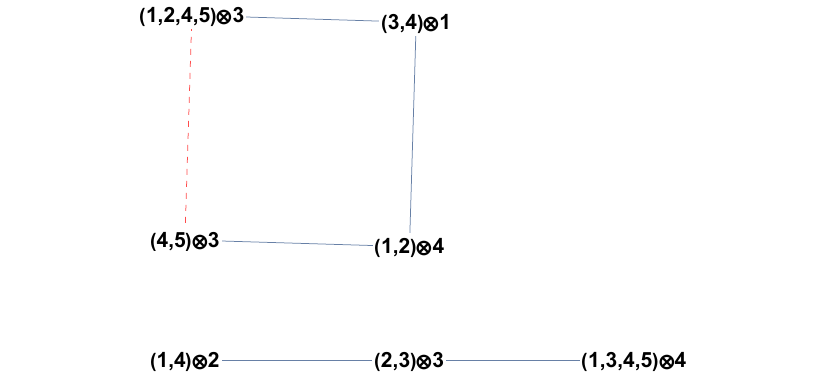}\\\hline
77&$\text{(1,2,3,4)$\otimes$2}\alb+\text{(1,3)$\otimes$4}\alb+\text{(1,5)$\otimes$3}\alb+\text{(2,3,4,5)$\otimes$3}\alb+\text{(2,4)$\otimes$4}\alb+\text{(3,4)$\otimes$1}$& $A_1\alb+\ttt_1\alb+\u_8$ &48&$\begin{array}{cccccccc}
 0 & 2 & 0 & 2 & \oplus  & 4 & 2 & 0 \\
  &  & 2\\
\end{array}$&\includegraphics[scale=\imgscale]{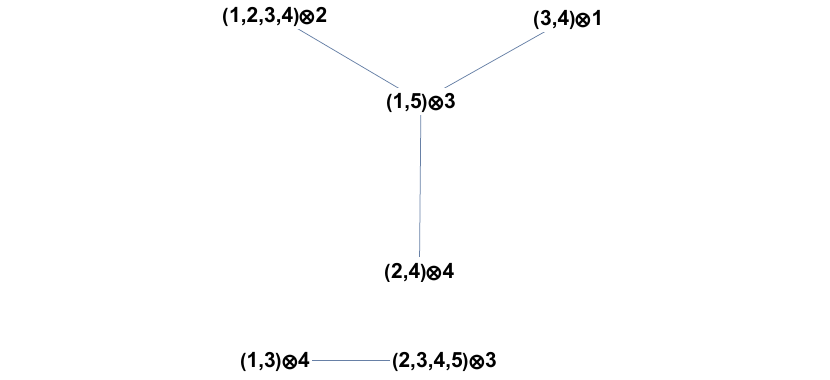}\\\hline
78&$\text{(1,2,4,5)$\otimes$3}\alb+\text{(1,3,4,5)$\otimes$4}\alb+\text{(1,3)$\otimes$1}\alb+\text{(2,4)$\otimes$4}\alb+\text{(3,4)$\otimes$2}\alb+\text{()$\otimes$3}$& $\ttt_2\alb+\u_{10}$ &48&$\begin{array}{cccccccc}
 0 & 2 & 0 & 2 & \oplus  & 2 & 0 & 2 \\
  &  & 2\\
\end{array}$&\includegraphics[scale=\imgscale]{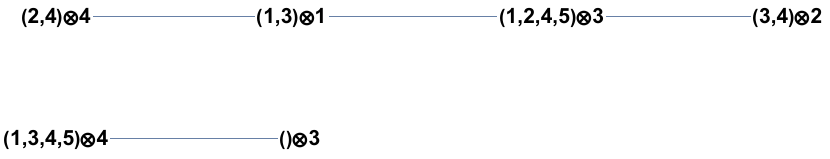}\\\hline
79&$\text{(1,2,3,4)$\otimes$1}\alb+\text{(1,2,4,5)$\otimes$3}\alb+\text{(1,2)$\otimes$4}\alb+\text{(1,3)$\otimes$3}\alb+\text{(2,4)$\otimes$2}\alb+\text{(4,5)$\otimes$4}$& $\ttt_2\alb+\u_{10}$ &48&$\begin{array}{cccccccc}
 0 & 3 & 0 & 2 & \oplus  & 0 & 3 & 0 \\
  &  & 0\\
\end{array}$&\includegraphics[scale=\imgscale]{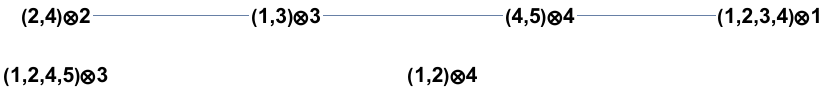}\\\hline
80&$\text{(1,2)$\otimes$4}\alb+\text{(1,3,4,5)$\otimes$3}\alb+\text{(1,3)$\otimes$1}\alb+\text{(2,4)$\otimes$1}\alb+\text{(3,4)$\otimes$2}\alb+\text{()$\otimes$3}$&  $\ttt_2\alb+\u_{10}$ &48&$\begin{array}{cccccccc}
 2 & 1 & 1 & 1 & \oplus  & 1 & 1 & 1 \\
  &  & 1\\
\end{array}$&\includegraphics[scale=\imgscale]{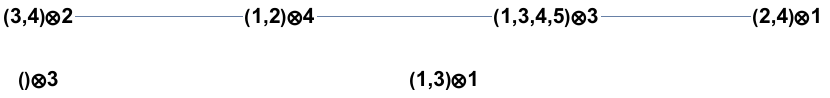}\\\hline
81&$\text{(1,2,4,5)$\otimes$3}\alb+\text{(1,3,4,5)$\otimes$4}\alb+\text{(1,3)$\otimes$1}\alb+\text{(2,3)$\otimes$3}\alb+\text{(2,4)$\otimes$1}\alb+\text{()$\otimes$4}$& $A_1\alb+\ttt_1\alb+\u_8$ &48&$\begin{array}{cccccccc}
 1 & 0 & 1 & 0 & \oplus  & 0 & 1 & 6 \\
  &  & 0\\
\end{array}$&\includegraphics[scale=\imgscale]{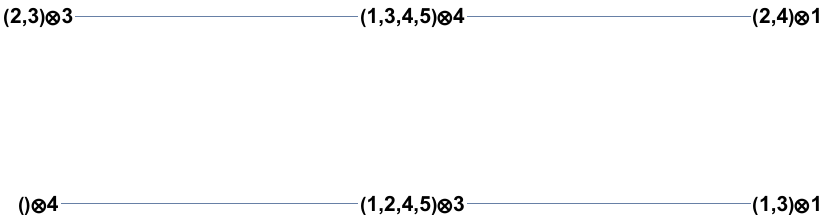}\\\hline
82&$\text{(1,3)$\otimes$4}\alb+\text{(1,4)$\otimes$2}\alb+\text{(1,5)$\otimes$3}\alb+\text{(2,3,4,5)$\otimes$3}\alb+\text{(2,4)$\otimes$1}$& $A_1\alb+\ttt_2\alb+\u_7$ &48&$\begin{array}{cccccccc}
 1 & 2 & 0 & 2 & \oplus  & 4 & 1 & 3 \\
  &  & 1\\
\end{array}$&\includegraphics[scale=\imgscale]{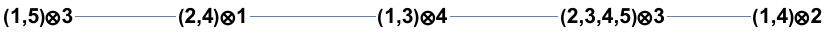}\\\hline
83&$\text{(1,2,3,4)$\otimes$2}\alb+\text{(1,2)$\otimes$1}\alb+\text{(1,3)$\otimes$4}\alb+\text{(2,4)$\otimes$1}\alb+\text{(4,5)$\otimes$3}$& $\ttt_3\alb+\u_9$ &48&$\begin{array}{cccccccc}
 2 & 0 & 1 & 4 & \oplus  & 3 & 0 & 1 \\
  &  & 0\\
\end{array}$&\includegraphics[scale=\imgscale]{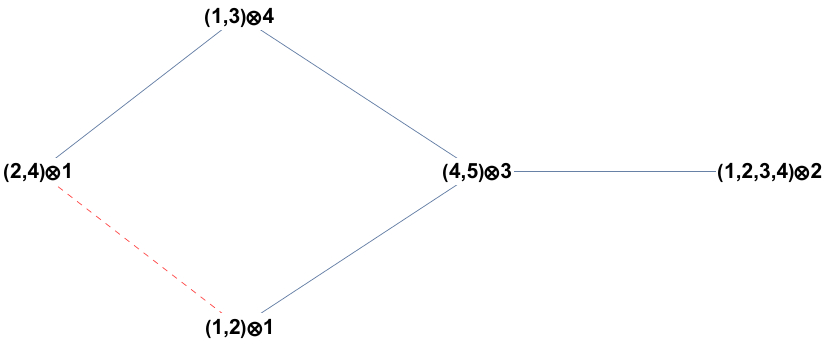}\\\hline
84&$\text{(1,2)$\otimes$4}\alb+\text{(1,3,4,5)$\otimes$4}\alb+\text{(1,5)$\otimes$3}\alb+\text{(2,3,4,5)$\otimes$3}\alb+\text{(3,4)$\otimes$1}$& $A_1\alb+\ttt_2\alb+\u_7$ &48&$\begin{array}{cccccccc}
 1 & 0 & 0 & 0 & \oplus  & 1 & 1 & 6 \\
  &  & 1\\
\end{array}$&\includegraphics[scale=\imgscale]{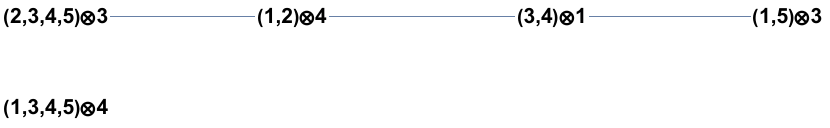}\\\hline
85&$\text{(1,2,4,5)$\otimes$3}\alb+\text{(1,2)$\otimes$4}\alb+\text{(1,3,4,5)$\otimes$4}\alb+\text{(1,3)$\otimes$1}\alb+\text{(1,4)$\otimes$2}\alb+\text{(2,4)$\otimes$1}\alb+\text{()$\otimes$3}$& $A_1\alb+\u_{10}$&47&$\begin{array}{cccccccc}
 0 & 0 & 0 & 4 & \oplus  & 0 & 0 & 4 \\
  &  & 0\\
\end{array}$&\includegraphics[scale=\imgscale]{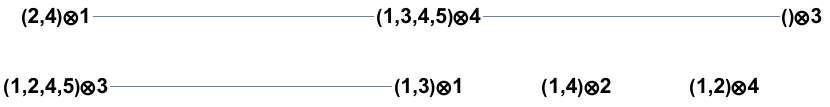}\\\hline
86&$\text{(1,2,4,5)$\otimes$3}\alb+\text{(1,2)$\otimes$4}\alb+\text{(1,3,4,5)$\otimes$4}\alb+\text{(1,3)$\otimes$1}\alb+\text{(2,4)$\otimes$1}\alb+\text{(3,4)$\otimes$2}$& $2A_1\alb+\u_7$ &47&$\begin{array}{cccccccc}
 0 & 0 & 4 & 0 & \oplus  & 0 & 0 & 0 \\
  &  & 0\\
\end{array}$&\includegraphics[scale=\imgscale]{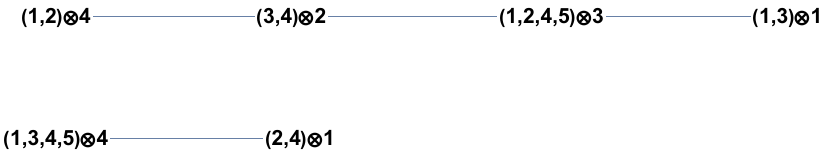}\\\hline
87&$\text{(1,2,3,4)$\otimes$2}\alb+\text{(1,3)$\otimes$4}\alb+\text{(1,5)$\otimes$3}\alb+\text{(2,3)$\otimes$3}\alb+\text{(2,4)$\otimes$4}\alb+\text{(3,4)$\otimes$1}$& $A_1\alb+\ttt_1\alb+\u_9$ &47&$\begin{array}{cccccccc}
 2 & 1 & 0 & 3 & \oplus  & 4 & 1 & 0 \\
  &  & 1\\
\end{array}$&\includegraphics[scale=\imgscale]{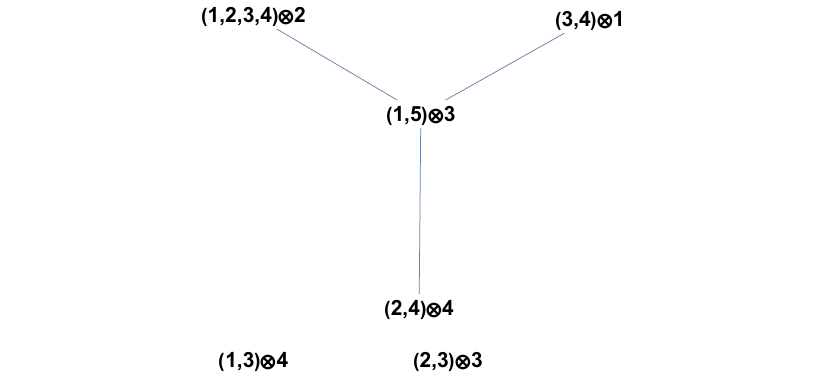}\\\hline
88&$\text{(1,2,3,4)$\otimes$2}\alb+\text{(1,2,4,5)$\otimes$3}\alb+\text{(1,2)$\otimes$4}\alb+\text{(1,3,4,5)$\otimes$4}\alb+\text{(3,4)$\otimes$1}\alb+\text{()$\otimes$3}$& $\ttt_2\alb+\u_{11}$ &47&$\begin{array}{cccccccc}
 1 & 1 & 1 & 1 & \oplus  & 1 & 2 & 1 \\
  &  & 1\\
\end{array}$&\includegraphics[scale=\imgscale]{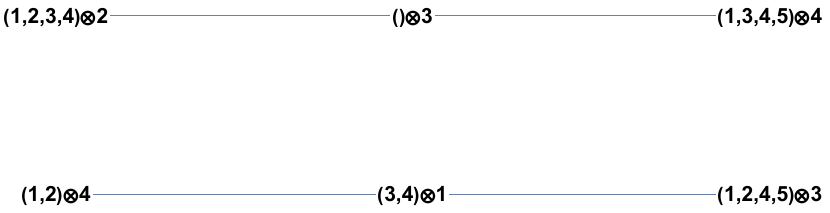}\\\hline
89&$\text{(1,2,4,5)$\otimes$3}\alb+\text{(1,3)$\otimes$4}\alb+\text{(2,4)$\otimes$1}\alb+\text{(3,4)$\otimes$2}\alb+\text{()$\otimes$3}$& $\ttt_3\alb+\u_{10}$ &47&$\begin{array}{cccccccc}
 1 & 2 & 1 & 1 & \oplus  & 1 & 2 & 0 \\
  &  & 0\\
\end{array}$&\includegraphics[scale=\imgscale]{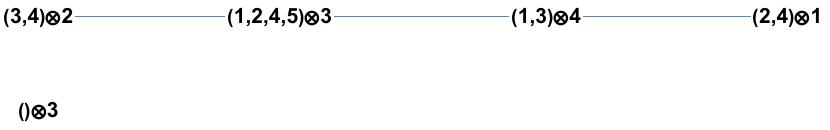}\\\hline
90&$\text{(1,2,4,5)$\otimes$3}\alb+\text{(1,2)$\otimes$4}\alb+\text{(1,3,4,5)$\otimes$3}\alb+\text{(1,4)$\otimes$2}\alb+\text{(2,3)$\otimes$3}\alb+\text{(3,4)$\otimes$1}$& $\ttt_2\alb+\u_{12}$ &46&$\begin{array}{cccccccc}
 2 & 0 & 2 & 0 & \oplus  & 2 & 0 & 2 \\
  &  & 0\\
\end{array}$&\includegraphics[scale=\imgscale]{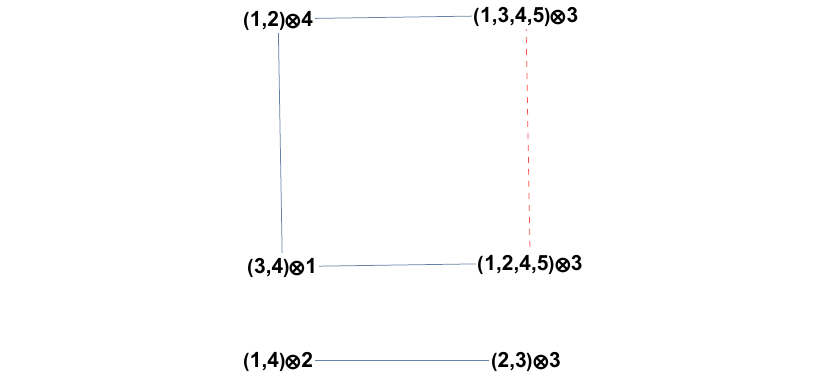}\\\hline
91&$\text{(1,2,4,5)$\otimes$3}\alb+\text{(1,2)$\otimes$4}\alb+\text{(1,3,4,5)$\otimes$4}\alb+\text{(2,3)$\otimes$3}\alb+\text{(3,4)$\otimes$1}\alb+\text{(4,5)$\otimes$3}$&  $A_1\alb+\ttt_1\alb+\u_{10}$ &46&$\begin{array}{cccccccc}
 0 & 2 & 0 & 0 & \oplus  & 0 & 2 & 4 \\
  &  & 0\\
\end{array}$&\includegraphics[scale=\imgscale]{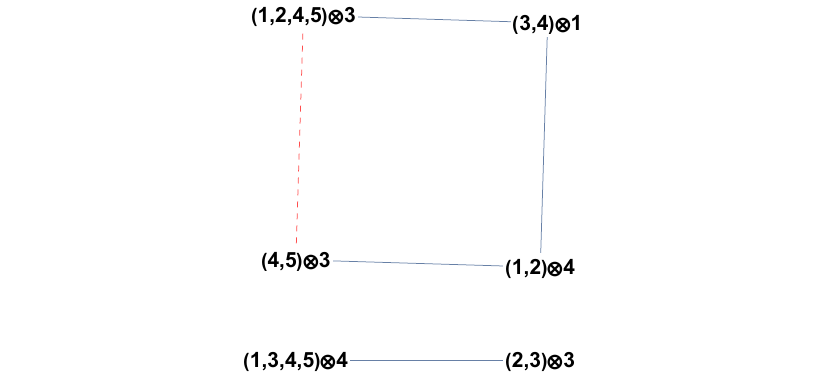}\\\hline
92&$\text{(1,2,3,4)$\otimes$1}\alb+\text{(1,2)$\otimes$4}\alb+\text{(1,3,4,5)$\otimes$4}\alb+\text{(1,4)$\otimes$2}\alb+\text{(2,3)$\otimes$3}\alb+\text{(4,5)$\otimes$3}$&  $A_1\alb+\ttt_1\alb+\u_{10}$&46&$\begin{array}{cccccccc}
 0 & 1 & 2 & 0 & \oplus  & 2 & 3 & 0 \\
  &  & 0\\
\end{array}$&\includegraphics[scale=\imgscale]{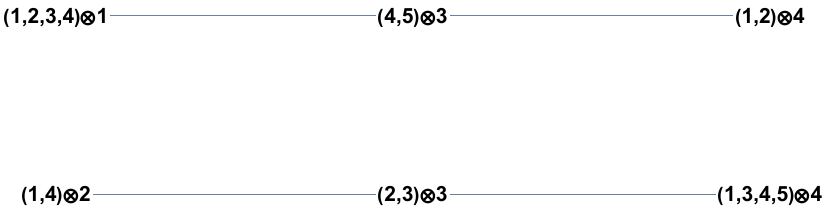}\\\hline
93&$\text{(1,2,4,5)$\otimes$3}\alb+\text{(1,2)$\otimes$4}\alb+\text{(1,3,4,5)$\otimes$4}\alb+\text{(1,3)$\otimes$1}\alb+\text{(2,4)$\otimes$1}\alb+\text{()$\otimes$3}$& $A_1\alb+\ttt_1\alb+\u_{10}$ &46&$\begin{array}{cccccccc}
 0 & 0 & 0 & 3 & \oplus  & 0 & 0 & 5 \\
  &  & 0\\
\end{array}$&\includegraphics[scale=\imgscale]{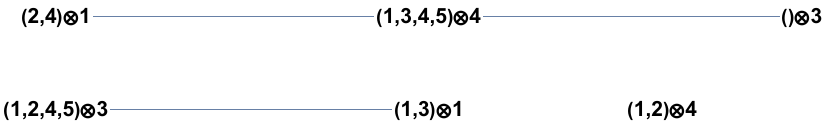}\\\hline
94&$\text{(1,2,4,5)$\otimes$3}\alb+\text{(1,2)$\otimes$4}\alb+\text{(1,3,4,5)$\otimes$4}\alb+\text{(1,4)$\otimes$2}\alb+\text{(3,4)$\otimes$1}\alb+\text{()$\otimes$3}$&  $\ttt_2\alb+\u_{12}$ &46&$\begin{array}{cccccccc}
 0 & 1 & 0 & 3 & \oplus  & 0 & 1 & 3 \\
  &  & 0\\
\end{array}$&\includegraphics[scale=\imgscale]{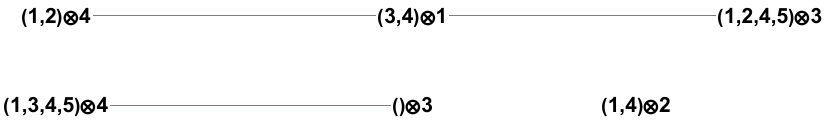}\\\hline
95&$\text{(1,2,3,4)$\otimes$2}\alb+\text{(1,3)$\otimes$4}\alb+\text{(1,5)$\otimes$3}\alb+\text{(2,4)$\otimes$4}\alb+\text{(3,4)$\otimes$1}$&    $A_1\alb+\ttt_2\alb+\u_{9}$&46&$\begin{array}{cccccccc}
 1 & 1 & 0 & 4 & \oplus  & 3 & 1 & 0 \\
  &  & 1\\
\end{array}$&\includegraphics[scale=\imgscale]{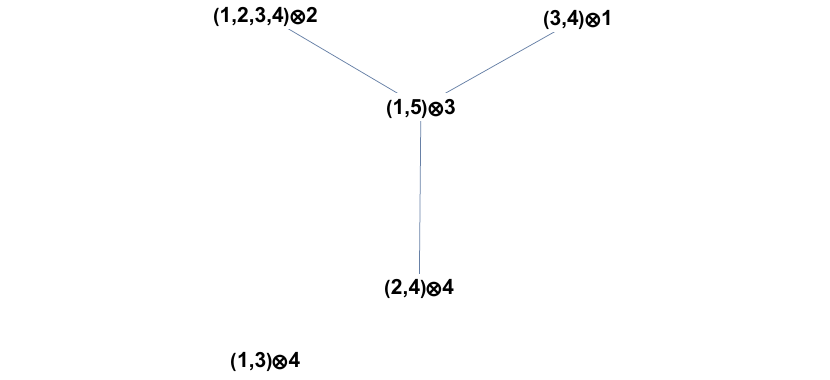}\\\hline
96&$\text{(1,2)$\otimes$3}\alb+\text{(1,3,4,5)$\otimes$4}\alb+\text{(1,3)$\otimes$1}\alb+\text{(2,4)$\otimes$1}\alb+\text{(3,4)$\otimes$2}$&   $A_1\alb+\ttt_2\alb+\u_{9}$&46&$\begin{array}{cccccccc}
 1 & 1 & 2 & 0 & \oplus  & 0 & 1 & 1 \\
  &  & 1\\
\end{array}$&\includegraphics[scale=\imgscale]{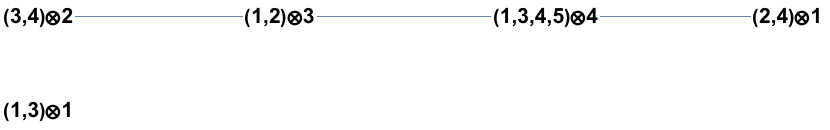}\\\hline
97&$\text{(1,2,4,5)$\otimes$3}\alb+\text{(1,2)$\otimes$4}\alb+\text{(1,3,4,5)$\otimes$3}\alb+\text{(1,4)$\otimes$2}\alb+\text{(3,4)$\otimes$1}\alb+\text{()$\otimes$3}$&  $\ttt_2\alb+\u_{13}$&45&$\begin{array}{cccccccc}
 1 & 0 & 1 & 2 & \oplus  & 1 & 0 & 3 \\
  &  & 0\\
\end{array}$&\includegraphics[scale=\imgscale]{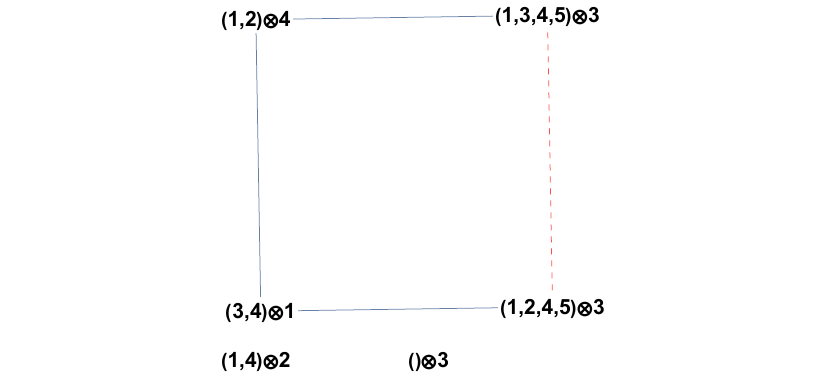}\\\hline
98&$\text{(1,2,3,4)$\otimes$1}\alb+\text{(1,2,4,5)$\otimes$3}\alb+\text{(1,3,4,5)$\otimes$4}\alb+\text{(1,4)$\otimes$2}\alb+\text{(2,3)$\otimes$3}\alb+\text{()$\otimes$4}$&  $2A_1\alb+\u_9$ &45&$\begin{array}{cccccccc}
 2 & 0 & 2 & 0 & \oplus  & 0 & 4 & 0 \\
  &  & 0\\
\end{array}$&\includegraphics[scale=\imgscale]{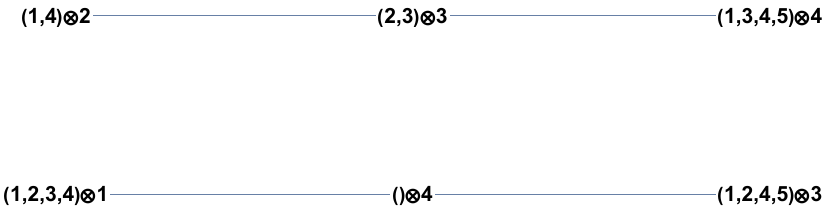}\\\hline
99&$\text{(1,2,4,5)$\otimes$3}\alb+\text{(1,3)$\otimes$4}\alb+\text{(1,4)$\otimes$2}\alb+\text{(2,3)$\otimes$3}\alb+\text{(2,4)$\otimes$4}\alb+\text{(3,4)$\otimes$1}$& $\ttt_2\alb+\u_{13}$ &45&$\begin{array}{cccccccc}
 1 & 1 & 1 & 0 & \oplus  & 2 & 1 & 1 \\
  &  & 1\\
\end{array}$&\includegraphics[scale=\imgscale]{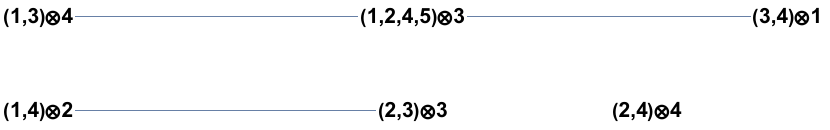}\\\hline
100&$\text{(1,2,4,5)$\otimes$3}\alb+\text{(1,2)$\otimes$4}\alb+\text{(1,3,4,5)$\otimes$4}\alb+\text{(3,4)$\otimes$1}\alb+\text{()$\otimes$3}$& $\ttt_3\alb+\u_{12}$&45&$\begin{array}{cccccccc}
 0 & 1 & 0 & 2 & \oplus  & 0 & 1 & 4 \\
  &  & 0\\
\end{array}$&\includegraphics[scale=\imgscale]{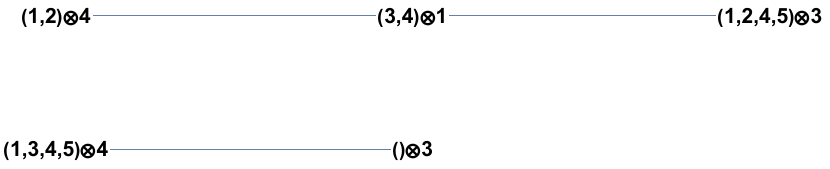}\\\hline
101&$\text{(1,2,4,5)$\otimes$3}\alb+\text{(1,3)$\otimes$4}\alb+\text{(2,4)$\otimes$1}\alb+\text{(3,4)$\otimes$2}$& $A_1\alb+\ttt_3\alb+\u_9$&45&$\begin{array}{cccccccc}
 0 & 2 & 2 & 0 & \oplus  & 0 & 2 & 0 \\
  &  & 0\\
\end{array}$&\includegraphics[scale=\imgscale]{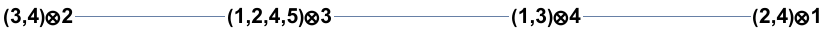}\\\hline
102&$\text{(1,2)$\otimes$4}\alb+\text{(1,5)$\otimes$3}\alb+\text{(2,3,4,5)$\otimes$3}\alb+\text{(3,4)$\otimes$1}$&  $A_2\alb+\ttt_2\alb+\u_5$ &45&$\begin{array}{cccccccc}
 2 & 0 & 0 & 0 & \oplus  & 2 & 0 & 6 \\
  &  & 0\\
\end{array}$&\includegraphics[scale=\imgscale]{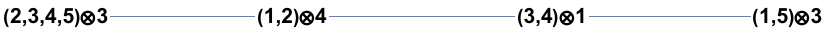}\\\hline
103&$\text{(1,2,4,5)$\otimes$3}\alb+\text{(1,3)$\otimes$4}\alb+\text{(1,4)$\otimes$2}\alb+\text{(2,4)$\otimes$4}\alb+\text{(3,4)$\otimes$1}\alb+\text{()$\otimes$3}$& $\ttt_2\alb+\u_{14}$&44&$\begin{array}{cccccccc}
 0 & 1 & 0 & 2 & \oplus  & 1 & 1 & 2 \\
  &  & 1\\
\end{array}$&\includegraphics[scale=\imgscale]{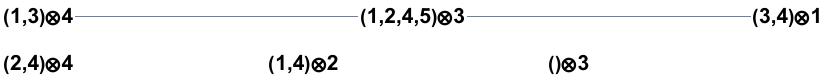}\\\hline
104&$\text{(1,2,3,4)$\otimes$2}\alb+\text{(1,5)$\otimes$3}\alb+\text{(2,3)$\otimes$3}\alb+\text{(2,4)$\otimes$4}\alb+\text{(3,4)$\otimes$1}$& $A_2\alb+\ttt_1\alb+\u_7$ &44&$\begin{array}{cccccccc}
 1 & 2 & 0 & 3 & \oplus  & 5 & 0 & 0 \\
  &  & 0\\
\end{array}$&\includegraphics[scale=\imgscale]{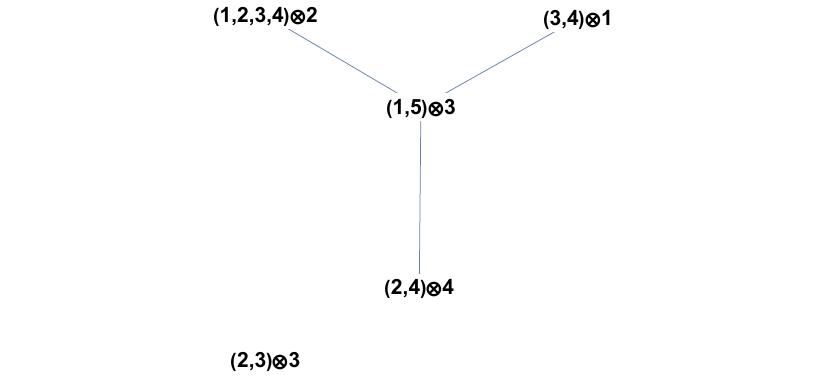}\\\hline
105&$\text{(1,2,4,5)$\otimes$3}\alb+\text{(1,2)$\otimes$4}\alb+\text{(1,3,4,5)$\otimes$3}\alb+\text{(3,4)$\otimes$1}\alb+\text{()$\otimes$3}$& $\ttt_3\alb+\u_{13}$&44&$\begin{array}{cccccccc}
 1 & 0 & 1 & 1 & \oplus  & 1 & 0 & 4 \\
  &  & 0\\
\end{array}$&\includegraphics[scale=\imgscale]{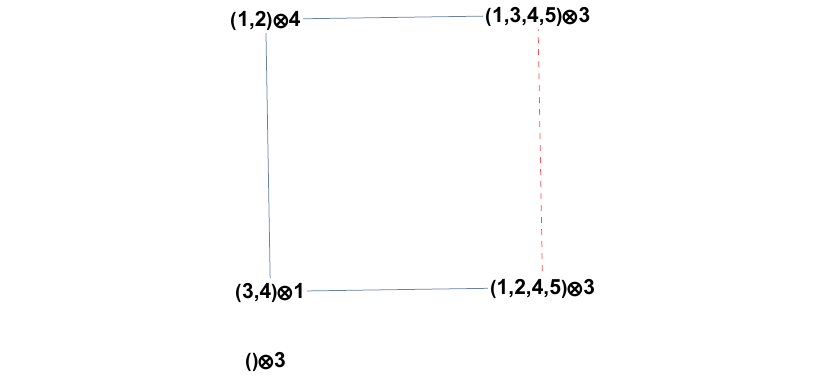}\\\hline
106&$\text{(1,2,4,5)$\otimes$3}\alb+\text{(1,2)$\otimes$4}\alb+\text{(1,3,4,5)$\otimes$3}\alb+\text{(1,4)$\otimes$2}\alb+\text{(3,4)$\otimes$1}$& $\ttt_3\alb+\u_{13}$  &44&$\begin{array}{cccccccc}
 0 & 0 & 2 & 1 & \oplus  & 0 & 0 & 3 \\
  &  & 0\\
\end{array}$&\includegraphics[scale=\imgscale]{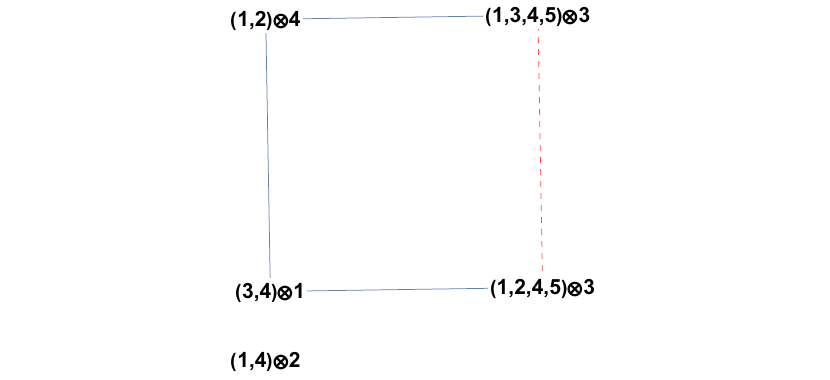}\\\hline
107&$\text{(1,3,4,5)$\otimes$4}\alb+\text{(1,4)$\otimes$1}\alb+\text{(1,5)$\otimes$3}\alb+\text{(2,3)$\otimes$3}\alb+\text{(2,4)$\otimes$4}$& $A_1\alb+\ttt_2\alb+\u_{11}$&44&$\begin{array}{cccccccc}
 1 & 0 & 0 & 1 & \oplus  & 1 & 2 & 3 \\
  &  & 1\\
\end{array}$&\includegraphics[scale=\imgscale]{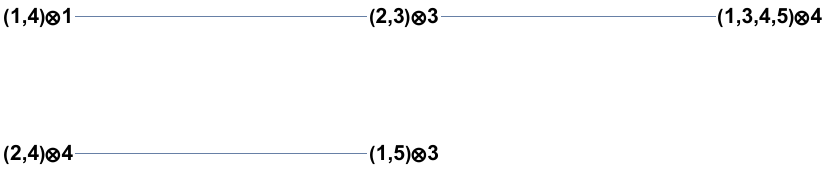}\\\hline
108&$\text{(1,3,4,5)$\otimes$3}\alb+\text{(1,3)$\otimes$1}\alb+\text{(1,4)$\otimes$2}\alb+\text{(2,3)$\otimes$3}\alb+\text{(2,4)$\otimes$4}$& $A_1\alb+\ttt_2\alb+\u_{11}$ &44&$\begin{array}{cccccccc}
 2 & 1 & 0 & 0 & \oplus  & 1 & 1 & 1 \\
  &  & 2\\
\end{array}$&\includegraphics[scale=\imgscale]{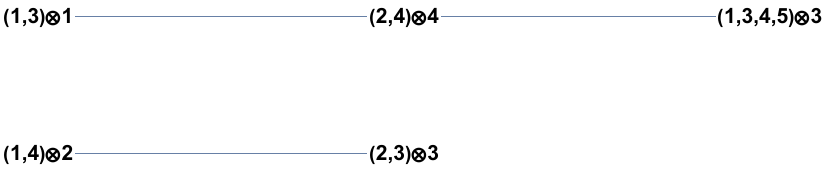}\\\hline
109&$\text{(1,2,4,5)$\otimes$3}\alb+\text{(1,3)$\otimes$4}\alb+\text{(1,4)$\otimes$2}\alb+\text{(2,3)$\otimes$3}\alb+\text{(3,4)$\otimes$1}$& $A_1\alb+\ttt_2\alb+\u_{12}$&43&$\begin{array}{cccccccc}
 1 & 0 & 1 & 0 & \oplus  & 3 & 0 & 1 \\
  &  & 2\\
\end{array}$&\includegraphics[scale=\imgscale]{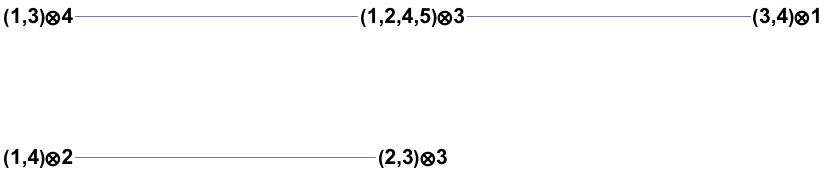}\\\hline
110&$\text{(1,2,4,5)$\otimes$3}\alb+\text{(1,3)$\otimes$4}\alb+\text{(2,4)$\otimes$4}\alb+\text{(3,4)$\otimes$1}\alb+\text{()$\otimes$3}$& $\ttt_3\alb+\u_{14}$&43&$\begin{array}{cccccccc}
 0 & 1 & 0 & 1 & \oplus  & 1 & 1 & 3 \\
  &  & 1\\
\end{array}$&\includegraphics[scale=\imgscale]{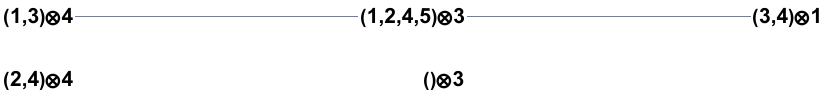}\\\hline
111&$\text{(1,2)$\otimes$3}\alb+\text{(1,3,4,5)$\otimes$4}\alb+\text{(1,4)$\otimes$2}\alb+\text{(2,4)$\otimes$4}\alb+\text{(3,4)$\otimes$1}$& $\ttt_3\alb+\u_{14}$ &43&$\begin{array}{cccccccc}
 1 & 0 & 1 & 1 & \oplus  & 0 & 1 & 2 \\
  &  & 1\\
\end{array}$&\includegraphics[scale=\imgscale]{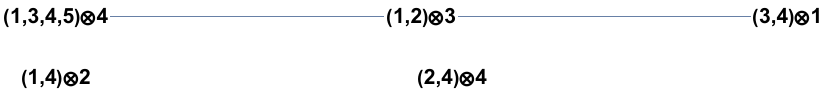}\\\hline
112&$\text{(1,2,3,4)$\otimes$1}\alb+\text{(1,3,4,5)$\otimes$3}\alb+\text{(1,3)$\otimes$4}\alb+\text{(1,4)$\otimes$2}\alb+\text{(2,4)$\otimes$4}\alb+\text{()$\otimes$3}$& $\ttt_2\alb+\u_{16}$&42&$\begin{array}{cccccccc}
 1 & 1 & 0 & 1 & \oplus  & 1 & 1 & 1 \\
  &  & 1\\
\end{array}$&\includegraphics[scale=\imgscale]{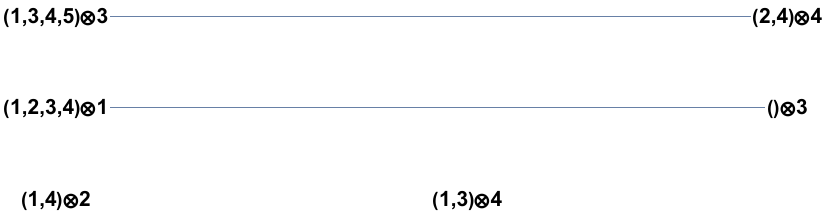}\\\hline
113&$\text{(3,4)$\otimes$1}\alb+\text{(1,2,4,5)$\otimes$3}\alb+\text{(1,3,4,5)$\otimes$3}\alb+\text{(1,2)$\otimes$4}$& $A_1\alb+\ttt_3\alb+\u_{12}$&42&$\begin{array}{cccccccc}
 0 & 0 & 2 & 0 & \oplus  & 0 & 0 & 4 \\
  &  & 0\\
\end{array}$&\includegraphics[scale=\imgscale]{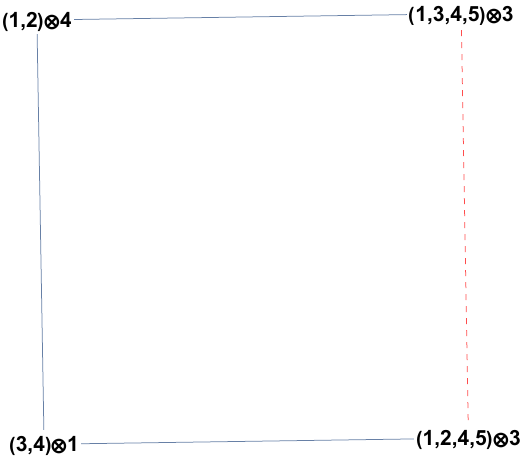}\\\hline
114&$\text{(1,2,3,4)$\otimes$2}\alb+\text{(1,5)$\otimes$3}\alb+\text{(2,4)$\otimes$4}\alb+\text{(3,4)$\otimes$1}$& $A_1\alb+A_2\alb+\ttt_1\alb+\u_6$ &42&$\begin{array}{cccccccc}
 0 & 2 & 0 & 4 & \oplus  & 4 & 0 & 0 \\
  &  & 0\\
\end{array}$&\includegraphics[scale=\imgscale]{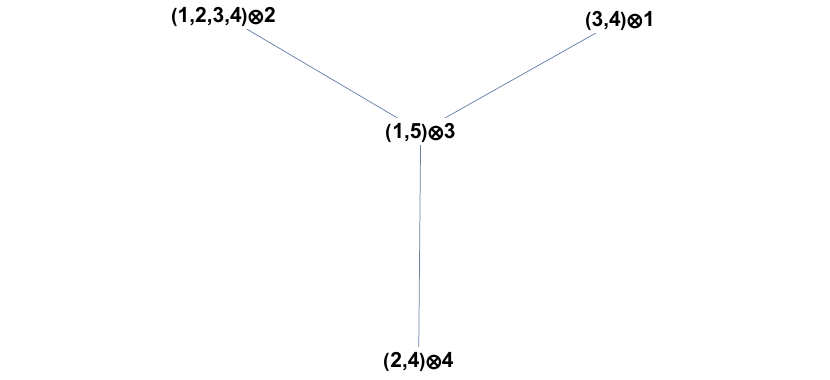}\\\hline
115&$\text{(1,2,4,5)$\otimes$3}\alb+\text{(1,3)$\otimes$4}\alb+\text{(1,4)$\otimes$2}\alb+\text{(3,4)$\otimes$1}\alb+\text{()$\otimes$3}$& $2A_1\alb+\ttt_1\alb+\u_{12}$&41&$\begin{array}{cccccccc}
 0 & 0 & 0 & 2 & \oplus  & 2 & 0 & 2 \\
  &  & 2\\
\end{array}$&\includegraphics[scale=\imgscale]{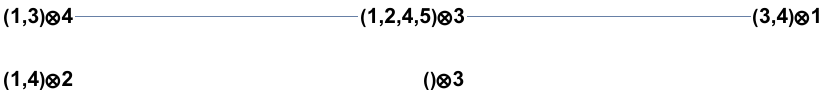}\\\hline
116&$\text{(1,2)$\otimes$3}\alb+\text{(1,3,4,5)$\otimes$4}\alb+\text{(2,4)$\otimes$4}\alb+\text{(3,4)$\otimes$1}$& $A_1\alb+\ttt_3\alb+\u_{13}$&41&$\begin{array}{cccccccc}
 1 & 0 & 1 & 0 & \oplus  & 0 & 1 & 3 \\
  &  & 1\\
\end{array}$&\includegraphics[scale=\imgscale]{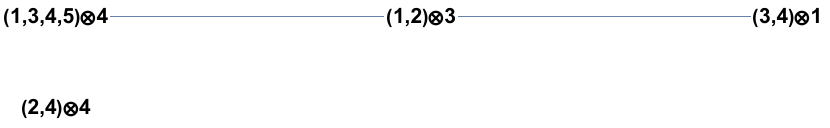}\\\hline
117&$\text{(1,2,3,4)$\otimes$1}\alb+\text{(1,2,4,5)$\otimes$3}\alb+\text{(1,4)$\otimes$2}\alb+\text{(3,4)$\otimes$4}\alb+\text{()$\otimes$3}$& $A_1\alb+\ttt_2\alb+\u_{15}$&40&$\begin{array}{cccccccc}
 0 & 2 & 0 & 1 & \oplus  & 2 & 0 & 1 \\
  &  & 0\\
\end{array}$&\includegraphics[scale=\imgscale]{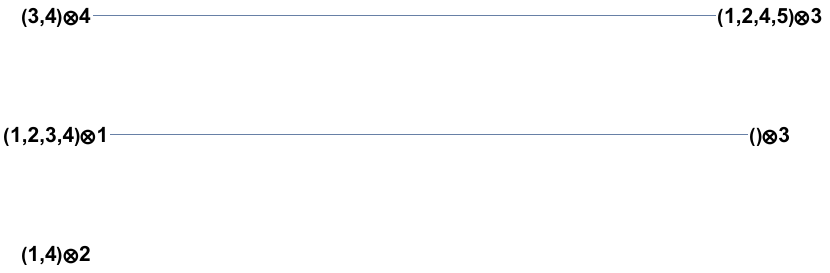}\\\hline
118&$\text{(1,3,4,5)$\otimes$3}\alb+\text{(1,3)$\otimes$4}\alb+\text{(1,4)$\otimes$1}\alb+\text{(2,3)$\otimes$3}\alb+\text{(2,4)$\otimes$4}$& $A_1\alb+\ttt_2\alb+\u_{15}$&40&$\begin{array}{cccccccc}
 1 & 1 & 0 & 0 & \oplus  & 1 & 1 & 2 \\
  &  & 1\\
\end{array}$&\includegraphics[scale=\imgscale]{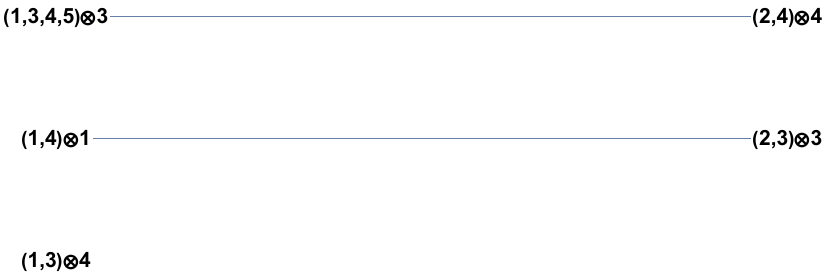}\\\hline
119&$\text{(1,2,4,5)$\otimes$3}\alb+\text{(1,3)$\otimes$4}\alb+\text{(3,4)$\otimes$1}\alb+\text{()$\otimes$3}$& $2A_1\alb+\ttt_2\alb+\u_{12}$&40&$\begin{array}{cccccccc}
 0 & 0 & 0 & 1 & \oplus  & 2 & 0 & 3 \\
  &  & 2\\
\end{array}$&\includegraphics[scale=\imgscale]{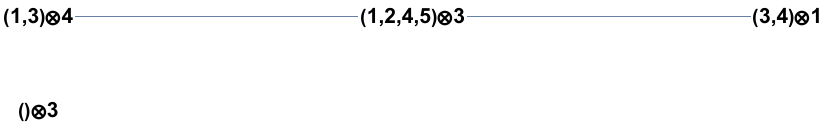}\\\hline
120&$\text{(1,3,4,5)$\otimes$4}\alb+\text{(1,3)$\otimes$1}\alb+\text{(1,4)$\otimes$2}\alb+\text{(2,4)$\otimes$3}$& $2A_1\alb+\ttt_2\alb+\u_{12}$&40&$\begin{array}{cccccccc}
 0 & 1 & 0 & 1 & \oplus  & 1 & 0 & 2 \\
  &  & 2\\
\end{array}$&\includegraphics[scale=\imgscale]{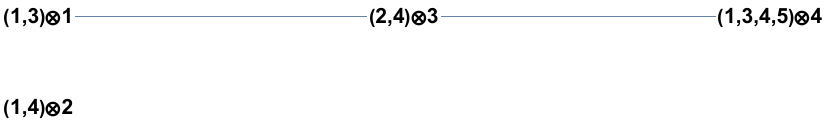}\\\hline
121&$\text{(1,2)$\otimes$4}\alb+\text{(1,4)$\otimes$2}\alb+\text{(2,3)$\otimes$3}\alb+\text{(3,4)$\otimes$1}$& $4A_1\alb+\u_8$&40&$\begin{array}{cccccccc}
 4 & 0 & 0 & 0 & \oplus  & 0 & 0 & 0 \\
  &  & 0\\
\end{array}$&\includegraphics[scale=\imgscale]{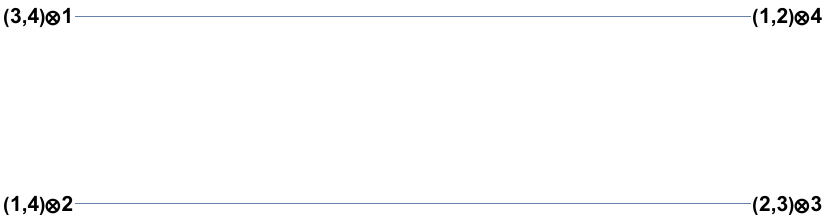}\\\hline
122&$\text{(1,2,3,4)$\otimes$1}\alb+\text{(1,2)$\otimes$3}\alb+\text{(1,3,4,5)$\otimes$3}\alb+\text{(1,3)$\otimes$4}\alb+\text{(1,4)$\otimes$2}\alb+\text{(2,4)$\otimes$4}$& $2A_1\alb+\u_{15}$&39&$\begin{array}{cccccccc}
 0 & 0 & 2 & 0 & \oplus  & 0 & 2 & 0 \\
  &  & 0\\
\end{array}$&\includegraphics[scale=\imgscale]{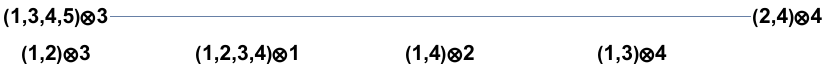}\\\hline
123&$\text{(1,2,3,4)$\otimes$1}\alb+\text{(1,3)$\otimes$4}\alb+\text{(1,4)$\otimes$2}\alb+\text{(2,4)$\otimes$4}\alb+\text{()$\otimes$3}$& $2A_1\alb+\ttt_1\alb+\u_{14}$&39&$\begin{array}{cccccccc}
 3 & 0 & 0 & 1 & \oplus  & 0 & 0 & 1 \\
  &  & 0\\
\end{array}$&\includegraphics[scale=\imgscale]{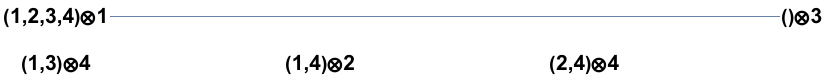}\\\hline
124&$\text{(1,2)$\otimes$3}\alb+\text{(1,3,4,5)$\otimes$3}\alb+\text{(1,3)$\otimes$4}\alb+\text{(1,4)$\otimes$1}\alb+\text{(2,4)$\otimes$4}$& $A_1\alb+\ttt_2\alb+\u_{17}$&38&$\begin{array}{cccccccc}
 0 & 0 & 1 & 1 & \oplus  & 0 & 2 & 1 \\
  &  & 0\\
\end{array}$&\includegraphics[scale=\imgscale]{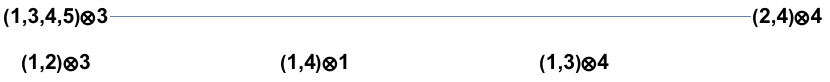}\\\hline
125&$\text{(1,2,3,4)$\otimes$1}\alb+\text{(1,2)$\otimes$3}\alb+\text{(1,3,4,5)$\otimes$3}\alb+\text{(1,4)$\otimes$2}\alb+\text{(3,4)$\otimes$4}$& $A_1\alb+\ttt_2\alb+\u_{17}$&38&$\begin{array}{cccccccc}
 1 & 0 & 1 & 0 & \oplus  & 1 & 1 & 0 \\
  &  & 1\\
\end{array}$&\includegraphics[scale=\imgscale]{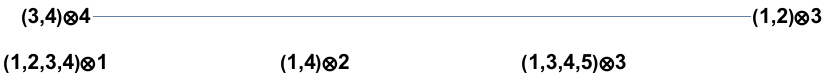}\\\hline
126&$\text{(1,2,3,4)$\otimes$4}\alb+\text{(1,4)$\otimes$1}\alb+\text{(2,3)$\otimes$3}\alb+\text{(4,5)$\otimes$3}$& $2A_1\alb+\ttt_2\alb+\u_{14}$&38&$\begin{array}{cccccccc}
 0 & 2 & 0 & 0 & \oplus  & 2 & 0 & 2 \\
  &  & 0\\
\end{array}$&\includegraphics[scale=\imgscale]{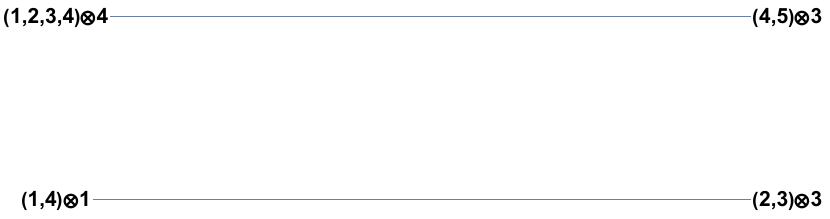}\\\hline
127&$\text{(1,2,3,4)$\otimes$1}\alb+\text{(1,2)$\otimes$3}\alb+\text{(1,4)$\otimes$2}\alb+\text{(3,4)$\otimes$4}$& $A_1\alb+\ttt_3\alb+\u_{17}$&37&$\begin{array}{cccccccc}
 2 & 0 & 1 & 0 & \oplus  & 0 & 1 & 0 \\
  &  & 0\\
\end{array}$&\includegraphics[scale=\imgscale]{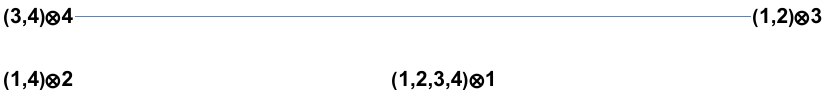}\\\hline
128&$\text{(1,3,4,5)$\otimes$4}\alb+\text{(1,3)$\otimes$1}\alb+\text{(2,4)$\otimes$3}$& $A_1\alb+A_2\alb+\ttt_2\alb+\u_{10}$&37&$\begin{array}{cccccccc}
 0 & 1 & 0 & 0 & \oplus  & 1 & 0 & 3 \\
  &  & 2\\
\end{array}$&\includegraphics[scale=\imgscale]{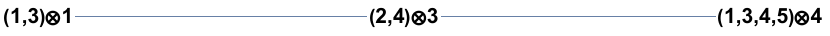}\\\hline
129&$\text{(1,2,3,4)$\otimes$4}\alb+\text{(1,3,4,5)$\otimes$3}\alb+\text{(1,4)$\otimes$1}\alb+\text{()$\otimes$3}$& $A_1\alb+\ttt_3\alb+u_{18}$&36&$\begin{array}{cccccccc}
 1 & 0 & 0 & 1 & \oplus  & 1 & 1 & 1 \\
  &  & 1\\
\end{array}$&\includegraphics[scale=\imgscale]{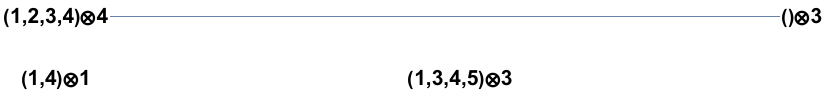}\\\hline
130&$\text{(1,3)$\otimes$4}\alb+\text{(1,5)$\otimes$3}\alb+\text{(2,3,4,5)$\otimes$3}\alb+\text{(2,4)$\otimes$4}$& $2A_1\alb+G_2\alb+\u_4$ &36&$\begin{array}{cccccccc}
 0 & 0 & 0 & 0 & \oplus  & 0 & 4 & 0 \\
  &  & 0\\
\end{array}$&\includegraphics[scale=\imgscale]{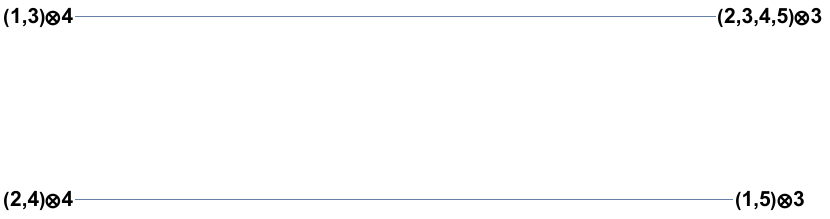}\\\hline
131&$\text{(1,3)$\otimes$4}\alb+\text{(1,4)$\otimes$1}\alb+\text{(2,3)$\otimes$3}\alb+\text{(2,4)$\otimes$4}$& $A_1\alb+B_2\alb+\ttt_1\alb+\u_{11}$&35&$\begin{array}{cccccccc}
 3 & 0 & 0 & 0 & \oplus  & 0 & 0 & 2 \\
  &  & 0\\
\end{array}$&\includegraphics[scale=\imgscale]{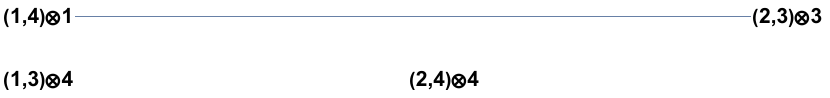}\\\hline
132&$\text{(1,2)$\otimes$3}\alb+\text{(1,3,4,5)$\otimes$3}\alb+\text{(1,3)$\otimes$4}\alb+\text{(2,4)$\otimes$4}$& $3A_1\alb+\ttt_1\alb+\u_{15}$&35&$\begin{array}{cccccccc}
 0 & 0 & 1 & 0 & \oplus  & 0 & 3 & 0 \\
  &  & 0\\
\end{array}$&\includegraphics[scale=\imgscale]{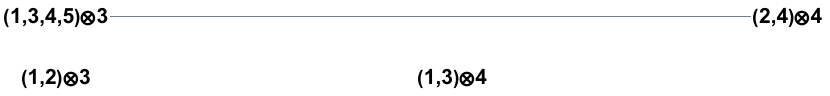}\\\hline
133&$\text{(1,2,3,4)$\otimes$1}\alb+\text{(1,3)$\otimes$3}\alb+\text{(1,4)$\otimes$2}\alb+\text{(2,4)$\otimes$3}\alb+\text{(3,4)$\otimes$4}$& $A_2\alb+\ttt_1\alb+\u_{17}$&34&$\begin{array}{cccccccc}
 1 & 1 & 0 & 0 & \oplus  & 1 & 0 & 0 \\
  &  & 1\\
\end{array}$&\includegraphics[scale=\imgscale]{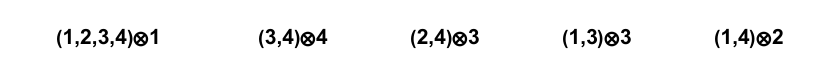}\\\hline
134&$\text{(1,2,3,4)$\otimes$4}\alb+\text{(1,4)$\otimes$1}\alb+\text{()$\otimes$3}$& $2A_1\alb+\ttt_3\alb+\u_{17}$ &34&$\begin{array}{cccccccc}
 2 & 0 & 0 & 1 & \oplus  & 0 & 1 & 1 \\
  &  & 0\\
\end{array}$&\includegraphics[scale=\imgscale]{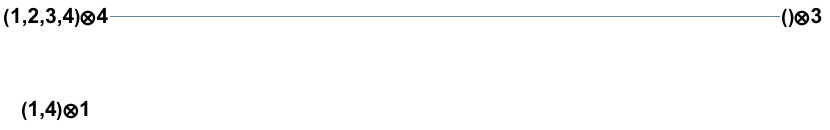}\\\hline
135&$\text{(1,2,3,4)$\otimes$4}\alb+\text{(1,3)$\otimes$3}\alb+\text{(1,4)$\otimes$1}\alb+\text{(2,4)$\otimes$3}$& $2A_1\alb+\ttt_2\alb+\u_{20}$ &32&$\begin{array}{cccccccc}
 1 & 0 & 1 & 0 & \oplus  & 1 & 0 & 1 \\
  &  & 0\\
\end{array}$&\includegraphics[scale=\imgscale]{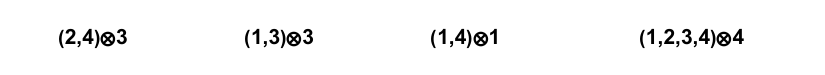}\\\hline
136&$\text{(1,3,4,5)$\otimes$3}\alb+\text{(1,4)$\otimes$4}\alb+\text{(2,3)$\otimes$3}$& $A_1\alb+A_2\alb+\ttt_2\alb+\u_{15}$&32&$\begin{array}{cccccccc}
 1 & 0 & 0 & 0 & \oplus  & 1 & 2 & 0 \\
  &  & 1\\
\end{array}$&\includegraphics[scale=\imgscale]{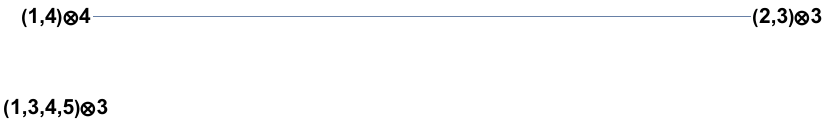}\\\hline
137&$\text{(1,2,3,4)$\otimes$1}\alb+\text{(1,3)$\otimes$3}\alb+\text{(1,4)$\otimes$2}\alb+\text{(3,4)$\otimes$4}$& $A_3\alb+\ttt_1\alb+\u_{14}$&30&$\begin{array}{cccccccc}
 1 & 0 & 0 & 0 & \oplus  & 0 & 0 & 0 \\
  &  & 2\\
\end{array}$&\includegraphics[scale=\imgscale]{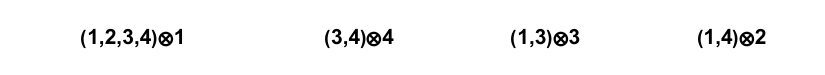}\\\hline
138&$\text{(1,2,3,4)$\otimes$1}\alb+\text{(1,4)$\otimes$2}\alb+\text{(2,4)$\otimes$3}\alb+\text{(3,4)$\otimes$4}$& $A_1\alb+A_3\alb+\u_{13}$ &29&$\begin{array}{cccccccc}
 0 & 2 & 0 & 0 & \oplus  & 0 & 0 & 0 \\
  &  & 0\\
\end{array}$&\includegraphics[scale=\imgscale]{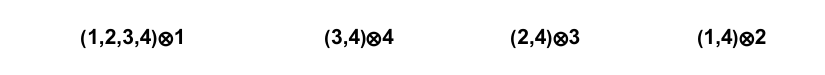}\\\hline
139&$\text{(1,2,3,4)$\otimes$4}\alb+\text{(1,4)$\otimes$1}\alb+\text{(3,4)$\otimes$3}$& $A_1\alb+A_2\alb+\ttt_2\alb+\u_{19}$&28&$\begin{array}{cccccccc}
 0 & 1 & 0 & 0 & \oplus  & 0 & 0 & 1 \\
  &  & 1\\
\end{array}$&\includegraphics[scale=\imgscale]{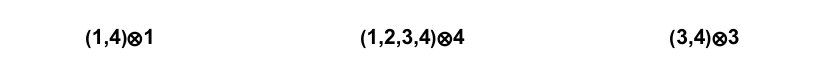}\\\hline
140&$\text{(1,4)$\otimes$4}\alb+\text{(2,3)$\otimes$3}$& $A_1\alb+A_3\alb+\ttt_2\alb+\u_{12}$&28&$\begin{array}{cccccccc}
 2 & 0 & 0 & 0 & \oplus  & 0 & 2 & 0 \\
  &  & 0\\
\end{array}$&\includegraphics[scale=\imgscale]{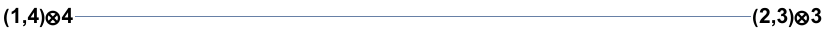}\\\hline
141&$\text{(1,3)$\otimes$3}\alb+\text{(1,4)$\otimes$4}\alb+\text{(2,4)$\otimes$3}$& $A_1\alb+B_2\alb+\ttt_2\alb+\u_{18}$ &27&$\begin{array}{cccccccc}
 1 & 0 & 0 & 1 & \oplus  & 1 & 1 & 0 \\
  &  & 0\\
\end{array}$&\includegraphics[scale=\imgscale]{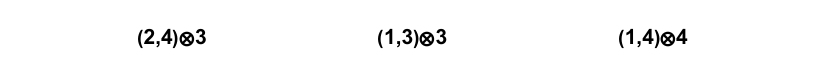}\\\hline
142&$\text{(1,2,3,4)$\otimes$3}\alb+\text{(1,4)$\otimes$4}$&  $3A_1\alb+A_2\alb+\ttt_1\alb+\u_{19}$&23&$\begin{array}{cccccccc}
 0 & 0 & 1 & 0 & \oplus  & 0 & 1 & 0 \\
  &  & 0\\
\end{array}$&\includegraphics[scale=\imgscale]{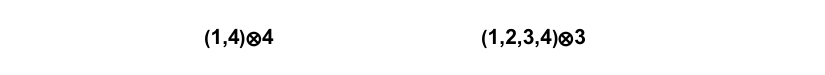}\\\hline
143&$\text{(1,3)$\otimes$3}\alb+\text{(2,4)$\otimes$3}$& $A_2\alb+B_3\alb+\ttt_1\alb+\u_{11}$&19&$\begin{array}{cccccccc}
 1 & 0 & 0 & 0 & \oplus  & 2 & 0 & 0 \\
  &  & 0\\
\end{array}$&\includegraphics[scale=\imgscale]{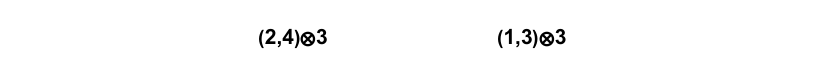}\\\hline
144&$\text{(1,4)$\otimes$3}$&  $A_2\alb+A_4\alb+\ttt_1\alb+\u_{13}$&14&$\begin{array}{cccccccc}
 0 & 0 & 0 & 1 & \oplus  & 1 & 0 & 0 \\
  &  & 0\\
\end{array}$&\includegraphics[scale=\imgscale]{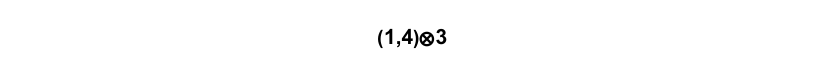}\\\hline
\end{longtable}
\end{tiny}
}

\ 

\begin{center}
\includegraphics[scale=.51]{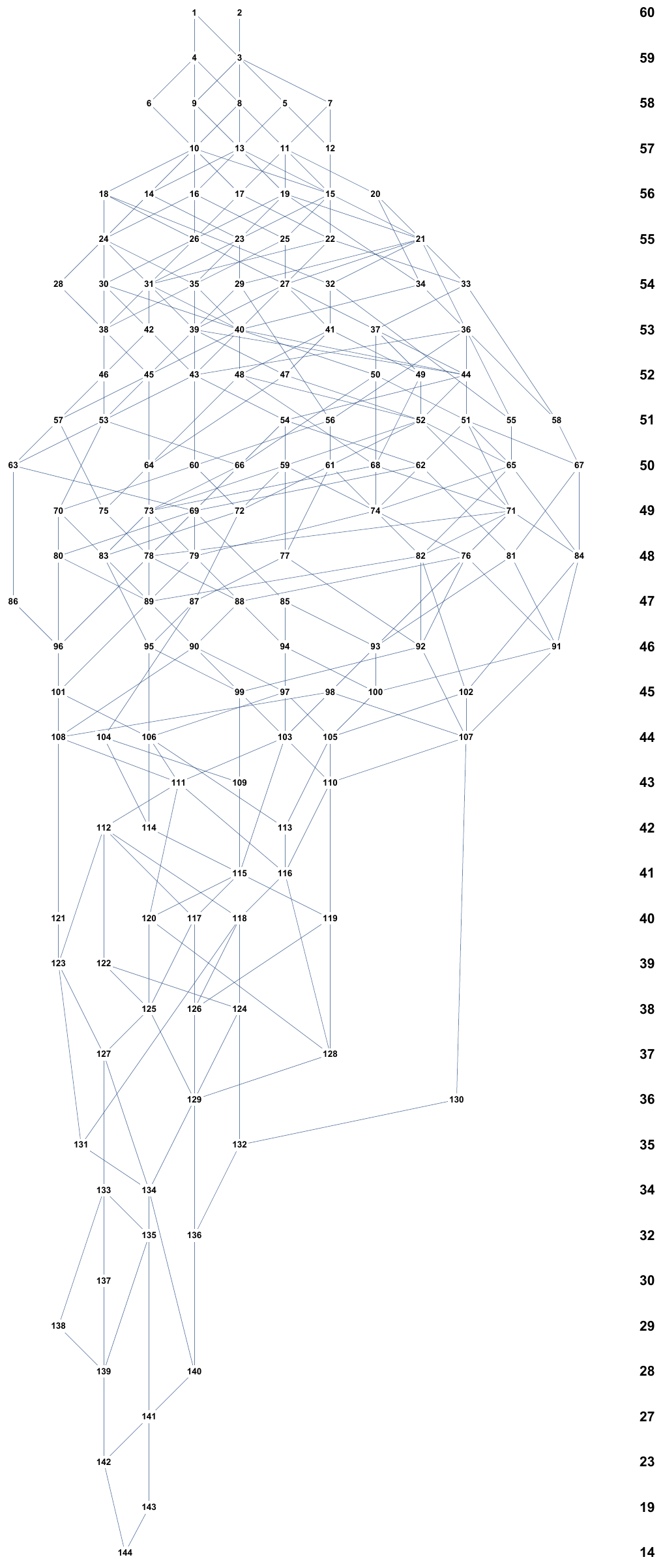}
\end{center}
\captionof{figure}{Hasse diagram of the closure ordering of the nilpotent orbits. The dimensions of the orbits are displayed on the right.}\label{fig}

\end{document}